\documentclass[12pt]{article}

\usepackage{latexsym,multicol,graphicx,color,comment}

\usepackage{amssymb}
\usepackage[centertags]{amsmath}
\usepackage{verbatim}

\usepackage[toc, page]{appendix}

\makeatletter
\renewcommand*{\@fnsymbol}[1]{\ensuremath{\ifcase#1\or *\or \star\or ***\or
   \mathsection\or \mathparagraph\or \|\or **\or \dagger\dagger
   \or \ddagger\ddagger \else\@ctrerr\fi}}
\makeatother

\usepackage{fancyhdr}
\newtheorem{theorem}{Theorem}
\newtheorem{lemma}{Lemma}

\newtheorem{definition}{Definition}

\numberwithin{equation}{subsection}

\newcommand{\dm}[1]{{\color{red}#1}}
\renewcommand{\dm}[1]{}

\newcommand{\ci}[1]{_{{}_{\scriptstyle #1}}}

\renewcommand{\phi}{\varphi} \renewcommand{\emptyset}{\varnothing}

\newcommand{\tame}{tame}

\newcommand{\e}{\varepsilon}  
  
  \newcommand{\la}{\lambda}
  \newcommand{\cF}{\mathcal F}
\newcommand{\cE}{\mathcal E}  
  
  \newcommand{\cP}{\mathcal P}
 \newcommand{\cH}{\mathcal H}

\newcommand{\BB}{\mathfrak B}

\newcommand{\vol}{\operatorname{vol}}

\def\done{{1\hskip-2.5pt{\rm l}}}

\renewcommand{\le}{\leqslant} \renewcommand{\ge}{\geqslant}

\newcommand{\cgf}{continuous Gaussian  function}
\newcommand{\pge}{parametric Gaussian ensemble}
\newcommand{\cccm}{connected component counting measure}

\newcommand{\bR}{\mathbb R} \newcommand{\bC}{\mathbb C} \newcommand{\bZ}{\mathbb Z}
 \newcommand{\bT}{\mathbb T} \newcommand{\bN}{\mathbb N}
 \newcommand{\bB}{\mathbb B} 
\newcommand{\bS}{\mathbb S} \newcommand{\bA}{\mathbb A}

\newcommand{\rd}{{\rm d}}

\newcommand{\spt}{\operatorname{spt}}

\begin{document}


\title
{Asymptotic laws for the spatial distribution\\ and the number
of connected components\\ of zero sets of Gaussian random functions}

\author{Fedor Nazarov
\thanks{Supported by grants No.\,2006136, 2012037 of the United States - Israel
Binational Science Foundation and by U.S. National Science Foundation Grants
DMS-0800243, DMS-1265623} \\
\and
\phantom{s}
\and
Mikhail Sodin
\thanks{Supported by grants No.\,2006136, 2012037 of the United States - Israel
Binational Science Foundation and by grant No.\,166/11
of the Israel Science Foundation of the Israel Academy
of Sciences and Humanities} }

\date{}
\maketitle

\vspace{-0.6cm}

\begin{tabular}{cccccc}\phantom{nushka} &
\begin{minipage}[c]{0,35\textwidth}
Dept. of Math. Sciences, \\
Kent State University, \\
Kent OH 44242, USA \\
nazarov@math.kent.edu
\end{minipage}
& \phantom{.} &
\begin{minipage}[c]{0,35\textwidth}
School of Math. Sciences\\ Tel Aviv University\\ Tel Aviv
69978, Israel \\
sodin@post.tau.ac.il
\end{minipage}
&\phantom{nush}
\end{tabular}

\vspace{1cm}

\rightline{\large\it In memory of Volodya Matsaev}

\vspace{0.6cm}

\begin{abstract}
We study the  asymptotic laws for the spatial distribution and the
number of connected components of zero sets of smooth Gaussian random functions
of several real variables. The primary examples are various Gaussian
ensembles of real-valued polynomials (algebraic or trigonometric) of large degree on
the sphere or torus, and translation-invariant smooth Gaussian functions on  the
Euclidean space restricted to large domains.
\end{abstract}

\tableofcontents

\section{Introduction and the main results}\label{sect:Intro}

The result we present has two main versions. The first one treats zero sets of smooth
Gaussian functions on the Euclidean space $\bR^m$ with translation-invariant distributions.
The second version deals with parametric ensembles of smooth Gaussian functions in open
domains in $\bR^m$. We also show how to translate the second version
to parametric ensembles of smooth Gaussian functions on smooth manifolds without boundary.

In Appendix~\ref{App-A}, all parts of the theory of smooth Gaussian functions needed for understanding
this work are developed from scratch. Appendix~\ref{App-B} contains the proof of the Fomin-Grenander-Maruyama
theorem in the multidimensional setting. None of the results in these Appendices is our own work.

\subsection{The translation invariant case}\label{subsect:Intro-transl-inv-case}

Suppose $F\colon \bR^m\to\bR$ is a continuous
Gaussian random function
with translation-invariant distribution (meaning that for every $v\in\bR^m$, the \cgf s $F$ and $F(\cdot+v)$
have the same distribution).
Then the covariance kernel
\[
K(x,y) = \cE \{ F(x)F(y)\}
\]
of $F$ depends only on the difference $x-y$ and can be
written in the form $K(x,y)=k(x-y)$ where $k:\bR^m\to\bR$ is a symmetric positive definite function.
By Bochner's theorem\footnote{
See \cite[\S~20]{Bochner} for the original proof or~\cite{Komatsu} for a clear self-contained exposition.},
$k$ can be written as the Fourier integral
$$
k(x)=(\cF\rho)(x)=\int_{\bR^m}e^{2\pi i(x\cdot\lambda)}\,\rd\rho(\lambda)
$$
of some finite symmetric positive Borel measure $\rho$ on $\bR^m$, which is called
{\em the spectral measure} of $F$.

We denote by $Z(F)=F^{-1}\{0\}$ the (random) zero set of $F$. Let $S$ be any bounded open convex set in $\bR^m$
containing the origin. By $S(R)$ we denote the set $\{x\in\bR^m\colon x/R\in S\}$.
By $N\ci S(R; F)$ we denote the number of the connected components of $Z(F)$ that are contained
in $S(R)$. When $S$ is the unit ball $B=\{x\colon |x|< 1\}$, we will write
$N(R; F)$ instead of $N\ci B(R; F)$.

We say that a finite complex-valued measure $\mu$ on $\bR^m$ is
{\em Hermitian} if for each Borel set $E\subset\bR^m$, we have
$\mu (-E) = \overline{\mu (E)}$. By $\cF\mu$ we denote the
Fourier integral of the measure $\mu$, and
by $\operatorname{spt}(\mu)$ we denote the closed support of $\mu$.

\begin{theorem}\label{thm:Euclid-version}\mbox{}
Suppose that the spectral measure $\rho$ of a continuous  Gaussian translation-invariant
function $F$ satisfies the following conditions:

\smallskip\par\noindent $(\rho1)$
\[
\int_{\bR^m} |\la|^4\, \rd\rho(\la) < \infty;
\]

\smallskip\par\noindent $(\rho2)$ $\rho$ has no atoms;


\smallskip\par\noindent $(\rho3)$ $\rho$ is not supported on a linear hyperplane.

\smallskip\par\noindent Then there exists a constant $\nu \ge 0$ such that for every bounded open convex set $S\subset\bR^m$
containing the origin,
\begin{equation}
\label{eq:a.s.}
\lim_{R\to\infty} \frac{N\ci S(R; F)}{\vol S(R)} = \nu  \ {\rm almost\ surely\qquad and} \qquad \lim_{R\to\infty} \cE \Bigl| \frac{N\ci S (R; F)}{\vol S(R)} - \nu \Bigr| = 0\,.
\end{equation}
\noindent Furthermore, $\nu > 0$ provided that

\smallskip\par\noindent $(\rho4)$
there exist a finite compactly supported Hermitian
measure $\mu$ with $\spt(\mu) \subset \spt(\rho) $ and a bounded domain $D\subset\bR^m$ such that
$\cF\mu\big|_{\partial D}<0$ and $ (\cF\mu) (u_0) > 0 $ for some $u_0\in D$.
\end{theorem}

\subsubsection{R\^{o}le of conditions $(\rho1)-(\rho3)$}
\label{subsubsect:assumptions_in_thm1}

Condition $(\rho1)$ guarantees that $F\in C^{2-}(\bR^m) \stackrel{\rm def}=\bigcap_{\alpha\in (0, 1)} C^{1+\alpha}(\bR^m)$.
Condition $(\rho3)$ says that the distribution of the gradient $\nabla F$ is non-degenerate.
Together conditions $(\rho1)$ and $(\rho3)$
allow us to think of the zero set $Z(F)$ as a collection of pairwise
disjoint smooth hypersurfaces that partition $\bR^m$ into ``nodal domains''.

The translation invariance allows us to consider the probability distribution measure generated
by $F$ on an appropriate space of functions as an invariant measure with respect to the action
of the abelian group $\bR^m$ by translations $( \tau_v g )( \cdot ) = g( \cdot + v )$. Condition ($\rho$2) ensures
that this action is ergodic, which in turn implies that the limit $\nu$ is non-random.

\subsubsection{Condition ($\rho$4)}\label{subsubsect:rho3}
Condition
$(\rho4)$ is essentially equivalent to the possibility to deterministically create at
least one bounded connected component of the zero set $Z(F)$. The measures not satisfying $(\rho4)$
have to be very degenerate. In particular, the support of any measure not satisfying
($\rho$4) has to be contained in a quadratic hypersurface in $\bR^m$.
We prove this, as well as some other observations pertaining to condition ($\rho$4), in Appendix~\ref{App-C}.

On the other hand, the Fourier transform of the Lebesgue surface measure on the sphere
centered at the origin is radial and sign changing. So if $\spt(\rho)$ is a sphere in $\bR^m$
centered at the origin then ($\rho$4) is still satisfied.

These observations suffice to check condition ($\rho$4) in most interesting examples.

\subsubsection{What can be said about the constant $\nu$?}\label{subsubsect:nu}
Unfortunately, the proof of Theorem~\ref{thm:Euclid-version} does not provide
much information about the value of the constant $\nu $.
There is a huge discrepancy between the lower
bounds that can be extracted from the ``barrier construction'' introduced in~\cite{NS}
and the upper bounds obtained by computing the mean number of
special points in the nodal domains or in the zero set
(cf. Nastasescu's undergraduate thesis~\cite{Nast}).

It is worth noting that the limiting constant $\bar\nu$ equals the expectation
$ \cE \bigl\{ {\rm vol}(G_0)^{-1}  \bigr\}  $, where $G_0$ is the connected component
of $\bR^m\setminus Z(F)$ containing the origin (or any other given point in $\bR^m$).
The random variable $ {\rm vol}(G_0) $ is, perhaps, even more mysterious than
$N(R; F)$. Our theorem shows that $\cP \bigl\{ {\rm vol}(G_0) < +\infty \bigr\} > 0 $,
but we still do not even know how to prove that this probability is $1$, not mentioning
any efficient tail estimate for its distribution.

\subsubsection{Further remarks about Theorem~\ref{thm:Euclid-version}}\label{subsubsect:remarks-to-thm1}

Theorem~\ref{thm:Euclid-version} can be viewed as a version of the ``law of large
numbers'' for the ``connected component process'' on $\bR^m$ associated with
the Gaussian function $F$. In most applications, one does not need as strong
convergence as is guaranteed by Theorem~\ref{thm:Euclid-version} and just the
convergence in probability (which is equivalent to the convergence in distribution for constant limits)
is enough.

Note also that the value of the intensity $\nu (F)$ is completely determined  by the covariance
kernel  $k(x-y)$ of $F$, or, which is the same, by the spectral measure $\rho$.

\medskip
Our last remark concerns a non-degenerate linear change of variables.
Let $T\colon \bR^m\to\bR^m$ be a non-degenerate linear operator and let
$\widetilde{F}(x)=F(Tx)$. Then $\widetilde F$ is also a Gaussian translation-invariant
function. Moreover, for every $S\subset\bR^m$ and $R>0$, we have
$N_{TS}(R; F)=N_S(R; \widetilde F)$, whence,
\[
\frac{\cE N_S(R; \widetilde F)}{\vol S(R)} = |\det T|\, \frac{\cE N_{TS}(R; F)}{\vol (TS)(R)}\,.
\]
Thus, if the intensity $\nu (F)$ exists, then so does $\nu(\widetilde F)$, and we have the
relation \[ \nu(\widetilde F)=|\det T|\, \nu (F)\,.\]

\subsection{Parametric Gaussian ensembles}\label{subsect:Intro-parametric-ensembles}

\begin{definition}[\pge]\label{def:parametric}
{\rm
A \pge\ $(f\ci L)$ on an open set $U\subset\bR^m$ (or on an $m$-dimensional manifold $X$ without
boundary) is any family $(f\ci L)$ of \cgf s on $U$ (on $X$ respectively)
indexed by some countable set of numbers $L\ge 1$ accumulating only at $+\infty$.
}
\end{definition}

Many interesting \pge s (in particular, two examples considered below in Section~\ref{sect:Examples})
arise from the following construction.
Let $X$ be a smooth compact $m$-dimensional manifold without boundary.
Let $\cH_L$ be a sequence of real finite-dimensional
Hilbert spaces of continuous functions on $X$ indexed by some scaling
parameter $L\ge 1$ so that $\lim_{L\to\infty}\operatorname{dim}\cH_L=\infty$.
Since, for every $x\in X$,
the point evaluation $\cH_L\ni f\mapsto f(x)$ is a continuous linear functional on $\cH_L$,
there is a unique function $K_L^{x}\in\cH_L$ such that $f(x)=\langle f,K_L^{x}\rangle$.
The function $K_L(x,y)=K_L^{x}(y)$ is called {\em the reproducing kernel} of the
space $\cH_L$. Since, $K_L^{x}\in\cH_L$, we have $K_L^{x}(y)=\langle K_L^{x},K_L^{y}\rangle$,
so $K_L(x,y)$ is symmetric.
Now let $\bigl\{ e_k \bigr\}$ be an orthonormal basis in $\cH_L$.
Then for every $f\in\cH_L$, we have $f=\sum_k \langle f,e_k\rangle e_k$ in $\cH_L$ and, therefore, pointwise.
Thus
\[
K_L(x, y) = K_L^{x}(y)=\sum_k e_k(x)e_k(y).
\]

Consider the \cgf
\[
f\ci L (x) = \sum_k \xi_k e_k(x), \qquad x\in X\,,
\]
where $ \xi_k $ are independent standard real Gaussian random variables. The covariance
kernel of the Gaussian function $f\ci L$ equals
\[
\cE \bigl\{ f\ci L(x) f\ci L(y) \bigr\} = \sum_k e_k(x) e_k(y) = K_L(x,y)\,,
\]
so it does not depend on the choice of the orthonormal basis $\bigl\{ e_k \bigr\}$
and coincides with the
reproducing kernel of $\cH_L$. It follows that the distribution of $f\ci L$ also
does not depend on the choice of the basis and is completely determined by the space $\cH_L$
itself. We shall call this \cgf\ $f\ci L$ {\em the \cgf\ generated by $\cH_L$}.

\subsubsection{Normalization}\label{subsubsect:parametric-normalization}

We say that a \cgf\ $f$ on $U$ with the covariance kernel $K$ is {\em normalized} if
\[ \cE \bigl\{f(x)^2 \bigr\}= K(x, x) = 1 \qquad {\rm for\ all\ } x\in U\,. \]
If the random Gaussian function $f$ is not normalized
but non-degenerate (that is, $\cE \bigl\{f(x)^2\bigr\}>0$, or, what is the same,
$\cP\{f(x)=0\}=0$ for every $x\in U$), we can just replace
$f$ by $\widetilde f(x)=\frac{f(x)}{\sqrt{K(x,x)}}$, which will correspond to replacing
the covariance kernel $K(x,y)$ by
\begin{equation}\label{eq:norm1}
\widetilde{K}(x, y) = \frac{K(x, y)}{\sqrt{K(x, x) \cdot K(y, y)}}\,,
\end{equation}
without affecting the zero set $Z(f)$ in any way.

Note that if we allow $f$ to degenerate at some points uncontrollably, then the zero set of $f$
may contain deterministic pieces of arbitrarily complicated structure and our talk about
the asymptotic behavior of the number of nodal components of $f$ may easily become totally meaningless.
Thus,
\begin{itemize}
\item {\em we will always assume that all \cgf s and all \pge s in this paper are normalized}.
\end{itemize}

Note that in many basic examples, including the ones we consider below
in Section~\ref{sect:Examples}, the function $x\mapsto K_L(x,x)$ is constant, so the normalization
of $K$ reduces to the division by that constant.

\subsubsection{Scaling and translation-invariant local limits}\label{subsubsect:tr-inv-local-limits}

Let $U$ be an open set in $\bR^m$ and let $(f\ci L)$ be a \pge\ on $U$.
Let $K\ci L$ be the covariance kernel of $f\ci L$.

We define the scaled covariance kernel $K_{x, L}$ at a point $x\in U$ by
\[
K_{x, L}(u, v) = K_L\bigl( x+\frac{u}L, x+\frac{v}L\bigr).
\]
Note that $K_{x, L}$ is the covariance kernel of the scaled Gaussian function
\[
f_{x, L}(u) = f\ci L\bigl( x+\frac{u}L \bigr),
\]
i.e., $K_{x, L}(u, v) = \cE\bigl\{ f_{x, L}(u) f_{x, L}(v)\bigr\}$. Note also
that if $x\in U$ is fixed and $L\to\infty$, the sets
$U_{x,L}=\{u\in\bR^m\colon x+\frac uL\in U\}$ exhaust $\bR^m$.

\begin{definition}[translation-invariant limit]\label{def:loc_tr-inv}
{\rm Let $(f\ci L)$ be a \pge\ on an open set $U\subset\bR^m$ and let
$K\ci L$ be the covariance kernel of $f\ci L$. Let $x\in U$. We say that
the scaled covariance kernels $K_{x,L}$ have a
{\em translation-invariant
limit} if
there exists a continuous function
$k_x\colon \bR^m\to \bR$ such that, for each $u,v\in\bR^m$,
\begin{equation}\label{eq:loc_tr-inv_limit}
\lim_{L\to\infty}K_{x, L} (u, v) = k_x(u-v)\,.
\end{equation}
We say that the \pge\ $(f\ci L)$ has a translation invariant
limit at the point $x$ if there exists a translation invariant \cgf\ $F_x$ on $\bR^m$ such
that, for every finite point set $\mathcal U\in\bR^m$, the finite-dimensional Gaussian vectors
$f_{x, L}|_\mathcal U$ converge to $F_x|_\mathcal U$ in distribution.
}
\end{definition}
We call the function $F_x$ {\em the local limiting function}
and its spectral measure $\rho_x$ {\em the local limiting spectral measure} of the \pge\ $(f\ci L)$
at the point $x$.
If a \pge\ $(f\ci L)$ on $U$ has a translation invariant
limit $F_x$ at some point $x\in U$, then the scaled covariance kernels $K_{x,L}$ have a
translation invariant limit as well and the limiting kernel $k_x(u-v)$ is the
covariance kernel of $F_x$. On the other hand, without any additional assumptions,
covariance kernels $K_{x, L}(u, v)$ may have a translation invariant limit
$k_x(u-v)$ that corresponds to no continuous Gaussian function $F$. However,
within the set-up considered in this paper, these notions become equivalent.

\medskip

It is natural to believe that if a \pge\ $(f\ci L)$ on $U$ has a translation invariant limit
at every point $x\in U$, then for large $L$, we can count the nodal components of $f\ci L$ in some
open set $V\subset U$ by partitioning $V$ into nice sets $V_j$ of size larger than $1/L$, choosing some
points $x_j\in V_j$, approximating
the number of nodal components of $f\ci L$ in each set $V_j$ by the number of nodal components of
$F_{x_j}$ in $(V_j)_{x_j,L}=\bigl\{v\in\bR^m\colon x_j + L^{-1}v\in V_j  \bigr\}$, and adding all these counts up. If we are lucky enough, the nodal components
of $F_x$ may have asymptotic intensity $\bar\nu(x)=\nu(F_x)$ and then the total count we get will be
typically close to
\[ \sum_j\bar\nu(x_j)\vol (V_j)_{x_j, L} = L^m \sum_j\bar\nu(x_j)\vol V_j\,.\]
If we are even
luckier, the quantity $\bar\nu(x)$ may depend on $x$ in a nice enough way for the Riemann sums $\sum_j\bar\nu(x_j)\vol V_j$
to converge to $\int_{V}\bar\nu\,\rd \vol$.

The formalization of this intuitive argument requires some accuracy, especially because the
standard integral calculus nowadays is Lebesgue, not Riemann. The classical
form of a convergence statement for integrals in the Lebesgue  language is that of the dominated convergence theorem,
whose general structure is
\begin{itemize}\item
{\em Given a sequence of nice objects that converge in some fairly
weak and easy to check sense to some limiting object, and assuming that our pre-limiting objects are uniformly
controlled in some way, the limiting object is nice as well, and some integral functional of
the limiting object is the limit of the integral functionals of the pre-limit objects.}
\end{itemize}
Our Theorem~\ref{thm:main_thm} will be exactly of this structure.

\medskip
We have already introduced in Definition~\ref{def:loc_tr-inv} the modes of convergence we will
be using. Now it is time to define ``controllability''.

\subsubsection{Uniform smoothness of covariance kernels $K\ci{L}$}\label{subsubsect:unif-smoothness}

The control we want to impose will be two-fold.
First, we will need to restrict the typical speed of oscillation of the \cgf s $f\ci L$. Some restriction
of this type is inevitable because fast oscillating \cgf s like the Brownian motion on $\bR^1$ change sign
infinitely many times near every their zero and they still have fairly decent moduli of continuity
on the H\"older scale. The control we will impose will guarantee that $f_{x,L}\in C^{2-}(U)=\bigcap_{\tau\in (0, 1)} C^{1+\tau}(U)$
on every compact subset of $U$.

For $k\ge 1$, by $C^{k, k}(U\times U)$ we denote the class of functions $g\colon U\times U\to \bR$ for which
all partial derivatives $\partial_x^\alpha\, \partial_y^\beta\, g(x, y)$, $|\alpha|, |\beta|\le k$ (taken
in any order) exist and are continuous\footnote{
in which case, they do not depend on the order}.
For $L\ge 1$, a compact set $Q\subset U$, and $g\in C^{k, k}(U\times U)$, we put
\[
\| g \|_{L, Q, k} \stackrel{\rm def}= \max_{|\alpha|, |\beta|\le k} \ \max_{x, y\in Q}\
L^{-(|\alpha|+|\beta|)} \bigl| \partial_x^\alpha\, \partial_y^\beta\, g(x, y) \bigr|.
\]
When $L=1$, we will write $ \| g \|_{Q, k} $ instead of $ \| g \|_{1, Q, k} $.

If the covariance kernel $K$ of a continuous Gaussian function $f$ on $U$ belongs to  $C^{k, k}(U\times U)$,
then the semi-norms $\| K \|_{L, Q, k}$ can be computed on ``the diagonal'' $\alpha=\beta$ and $x=y$:
\[
\| K \|_{L, Q, k} = \max_{|\alpha|\le k}\ \max_{x\in Q}\ L^{-2|\alpha|} \bigl| \partial_x^\alpha\, \partial_y^\alpha\, K(x, y)|_{y=x}\, \bigr|.
\]
A na\"{\i}ve explanation to this fact comes from the Cauchy-Schwarz inequality combined with the formula
\[
\partial_x^\alpha\, \partial_y^\beta\, K(x, y) =
\cE \bigl\{ \partial_x^\alpha f(x) \cdot \partial_y^\beta f(y) \bigr\},
\]
which is true in the case when the derivatives on the RHS exist and are continuous random functions.
The proof of this fact for the general case will be given in Appendix~\ref{A_subsect:remarks_Kolmogorov}\,.

The uniform smoothness of the kernels
$K_{L}$ with $k\ge 1$ is more than enough to erase any distinction between the existence of a translation
invariant limit of the kernels $K_{x,L}$ and the existence of a translation
invariant limit at $x$ of the \pge\ $(f\ci L)$.

\subsubsection{Local uniform non-degeneracy of the \pge~$(f\ci L)$}\label{subsubsect:unif-nondegeneracy}

Our second restriction will be of the opposite character. While the local uniform smoothness
guarantees that the \cgf s $f\ci L$ do not change too fast or in a too rough way, the condition
we discuss in this section will ensure that $f\ci L$ cannot change too slowly or in a too
predictable way in any direction.
Without any such restriction, there will be nothing that would prevent long regular components to prevail
and, with our methods,
we will either not be able to say anything at all in such case, or will just conclude that some
limit is $0$, which would merely mean that the particular scaling we have chosen is a wrong one for the problem.
With all this in mind, let us pass to the formal definitions.

\medskip
Let $K\in C^{1,1}(U\times U)$ and let $C_x$ be the matrix with the entries
\[ C_x(i,j)=\partial_{x_i}\, \partial_{y_j}\, K(x, y)|_{y=x}\,, \qquad x\in U\,. \]
If $K$ is the covariance kernel of some $C^1$ Gaussian function $f$ on $U$, then $C_x$ is the covariance
matrix of the Gaussian random vector $\nabla f(x)$. Assuming that $\det C_x\ne 0$, we can say that the density
of the probability distribution of the random Gaussian vector $\nabla f(x)$ in $\bR^m$ is given by
$$
p(\xi)=\frac{1}{(2\pi)^{m/2}\sqrt{\det C_x}}\,e^{-\frac12 (C_x^{-1}\xi\cdot\xi)}\,.
$$
Since in this case $C_x^{-1}$ is positive definite, we have
\[ \max_\xi p(\xi)=p(0)=(2\pi)^{-m/2}(\det C_x)^{-1/2}\,.\]

\begin{definition}[local uniform non-degeneracy of $(f\ci L)$]\label{def:non-degen}
{\rm
We say that a \pge\ $(f\ci L)$ on some open set $U\subset\bR^m$ is locally
uniformly non-degenerate  if the corresponding kernels $K_L(x,y)$
are at least in $C^{1, 1}(U\times U)$ and
for every compact set $Q\subset U$,
$$
\varliminf_{L\to\infty} \, \inf_{x\in Q}\det C_{x,L}>0
$$
where $C_{x,L}$ is the matrix with the entries
$$
C_{x, L}(i, j) =
\partial_{u_i} \partial_{v_j} K_{x,L}(u, v)\bigl|_{u=v=0}=
L^{-2} \partial_{x_i} \partial_{y_j} K_{L}(x, y)\bigl|_{y=x}\,.
$$
}
\end{definition}

As the argument above shows, if $f\ci L$ is $C^1$-smooth, then our non-degeneracy condition
just means that, for every compact
set $Q\subset U$, there is a uniform upper bound for the densities of the distributions
of all Gaussian vectors $L^{-1}\nabla f\ci L(x)$ with $x\in Q$.

\medskip
Suppose that, for some $x\in U$,
the kernels $K_{x, L}$ have a translation invariant limit and the convergence holds in the
semi-norm $ \| \cdot \|_{Q, 1} $ for some compact set $Q$ containing $x$ in its interior. Then  the matrix $C_{x, L}$ converges to
the matrix $c_x$ with the entries
\[
c_x(i, j) = - \bigl( \partial_{u_i} \partial_{u_j} k_x\bigr)(0) = 4\pi^2\, \int_{\bR^m} \la_i \la_j\, \rd\rho_x(\la)
\]
and we see that in this case the limiting measure $\rho_x$ satisfies
\[
\inf_{\xi\in\bS^{m-1}}\, \int_{\bR^m} \bigl| \la \cdot \xi \bigr|^2\, \rd\rho_x(\la) > 0\,,
\]
which means that $\rho_x$ cannot be supported on any linear hyperplane $\bigl\{\la\colon \la\cdot\xi = 0 \bigr\}$, i.e.,
condition ($\rho3$) is satisfied.

\subsubsection{Controllability}\label{control}

Now we are ready to say what we mean by a locally uniformly controllable \pge\ $(f\ci L)$ on $U$.

\begin{definition}[locally uniform controllability]\label{def:control}
{\rm
The \pge\ $(f\ci L)$ on an open set $U\subset\bR^m$ is {\em locally uniformly controllable}
if it is locally
uniformly non-degenerate and the corresponding covariance kernels $K_L$ satisfy
\[
\varlimsup_{L\to\infty} \| K_L \|_{L, Q, 2} < \infty
\]
for every compact set $Q\subset U$.
}
\end{definition}

The above considerations combined with results presented in Appendix (see~\ref{A_subsect:remarks_Kolmogorov} and~\ref{A_subsect_tr-inv})
imply that
\begin{itemize} \item
{\em if the kernels $K_L$ are locally uniformly controllable and
if the scaled kernels $K_{x, L}$ have translation-invariant limits, then the limiting spectral
measure $\rho_x$ satisfies assumptions} ($\rho1$) {\em and} ($\rho3$) {\em of Theorem~\ref{thm:Euclid-version}}.
\end{itemize}

\subsubsection{Tameness}

\begin{definition}[\tame\ ensembles]\label{def:tame}
{\rm
The \pge\ $(f\ci L)$ on an open set $U\subset\bR^m$ is {\em tame}
if

\medskip\noindent
(i) it is locally uniformly controllable,

\medskip\noindent and there exists a Borel subset $U'\subset U$ of full Lebesgue measure
such that, for all $x\in U'$,

\medskip\noindent
(ii) the scaled
kernels $(K_{x, L})$ have translation invariant limits;

\smallskip\noindent
(iii) the limiting spectral measure $\rho_x$ has no atoms.
}
\end{definition}

\medskip
A tame \pge\ $ ( f\ci L )$ has a translation invariant limit at every point $x\in U'$. Moreover,
by Theorem~\ref{thm:Euclid-version}, the point intensity $\bar\nu (x)\stackrel{\rm def}= \nu (F_x)$
associated with the ensemble $(f\ci L)$ is well-defined on $U'$.

\subsection{The main result}\label{subsect:main_result}

Before we state our second main theorem, we will introduce one more object.
Let $U$ be an open set in $\bR^m$ and let $f$ be a \cgf\ on $U$.
We say that a (depending on the implicit probability variable $\omega$) Borel measure $n$ on $U$
is {\em a \cccm} of $f$ if $\spt(n)\subset Z(f)$ and the $n$-mass
of each connected component of $Z(f)$ equals $1$.
Note that
we do not require the dependence of $n$ on $\omega$ to be measurable in any sense
(for this reason, we do not call $n$ a random measure), so
in the statement of the next theorem we will have to use
``the upper expectation'' $\cE^*$ instead of the usual one $\cE$.

\begin{theorem}\label{thm:main_thm}
Suppose that $(f\ci L)$ is a \tame\ \pge\ on an open set $U\subset \bR^m$.
Then

\smallskip\par\noindent {\rm (i)}
the function $x\mapsto \bar\nu(x)$ is measurable and locally bounded in $U$;

\noindent and

\smallskip\par\noindent {\rm (ii)}
for every sequence of \cccm s $n\ci L$ of $f\ci L$ and for every compactly supported in $U$
continuous function $\phi$, we have
\[
\lim_{L\to\infty}\cE^* \Bigl\{ \Bigl| \frac1{L^m} \int \phi\, \rd n\ci L
- \int \phi \bar\nu\, {\rm d\, vol} \Bigr| \Bigr\} = 0\,.
\]
\end{theorem}

Note that the second statement of that theorem can be strengthened to
\[
\lim_{L\to\infty}\cE^* \Bigl\{ \Bigl| \frac1{L^m} \int \phi\, \rd n\ci L
- \int \phi \bar\nu\, {\rm d\, vol} \Bigr|^q \Bigr\} = 0
\]
for some $q=q(m)>1$ (which tends to $1$ as $m\to\infty$) without any essential
change in the proof but we are not aware of any application of this stronger result
for which the current version would not suffice as well.

\subsection{The manifold version of Theorem~\ref{thm:main_thm}}\label{subsect:manifold}

Theorem~\ref{thm:main_thm} can be transferred to \pge s on smooth manifolds without
boundary. Everywhere in this section $X$ is an $m$-dimensional $C^2$-manifold without boundary
(not necessarily compact) that can be covered by countably many charts, all charts being
assumed open and $C^2$-smooth, and $(f\ci L)$ is a \pge\ on $X$.
We start with two definitions.
\begin{definition}[\tame\ ensembles on manifolds]\label{def:tame_on_mflds}
{\rm
We say that a \pge\ $(f\ci L)$ on $X$ is {\em \tame} if, for every chart $\pi\colon U\to X$,
the \pge\ $(f\ci L \circ \pi)$ is \tame\ on $U$.
}
\end{definition}

This definition implies that for every chart $\pi\colon U\to X$, the \pge\ $(f\ci L\circ\pi)$
satisfies the assumptions of Theorem~\ref{thm:main_thm}. So the associated point intensity $\bar\nu_\pi$
belongs to $L^\infty_{\tt loc} (U)$.

\begin{definition}[Volumes compatible with smooth structure]
{\rm
We say that a locally finite Borel positive measure $\vol\ci X$ on $X$
is a volume {\em compatible with
the smooth structure} of $X$ if for every chart $\pi:U\to X$,
the measures $\pi_*\vol$ and  $\vol\ci X$ are mutually absolutely continuous and the corresponding Radon-Nikodym densities
are continuous on $\pi(U)$.
}
\end{definition}

Of course, the main example we have in mind giving this definition is that of a smooth Riemannian manifold $X$ and
the volume generated by the Riemannian metric on $X$.

We also note that despite the manifold
$X$ may be endowed with no measure, the words ``almost every $x\in X$'' still have meaning because all push-forward
measures $\pi_*\vol$ corresponding to various charts $\pi:U\to X$ of $X$ are mutually absolutely continuous wherever
they can be compared to each other.

\medskip At last, we can state the manifold version of Theorem~\ref{thm:main_thm}.

\begin{theorem}\label{thm:manifold}
Suppose that $(f\ci L)$ is a \tame\ \pge\ on $X$.
Then

\medskip\par\noindent {\rm (i)}
there exists a locally finite Borel non-negative measure $n_\infty$ on $X$ such that for every choice of \cccm s
$n\ci L$ of $f\ci L$ and every function $\phi\in C_0(X)$,
\[
\lim_{L\to\infty}\cE^* \Bigl\{ \Bigl| \frac1{L^m} \int \phi\, \rd n\ci L
- \int \phi \, \rd\, n_\infty \Bigr| \Bigr\} = 0\,;
\]

\medskip\par\noindent {\rm (ii)}
for every chart $\pi:U\to X$, the measure $n_\infty$ coincides on $\pi(U)$ with the push-forward
$\pi_*(\bar\nu_\pi\,\vol)$ where $\bar\nu_\pi$ is the point intensity associated with the
\pge\ $f\ci L\circ\pi$;

\medskip\par\noindent {\rm (iii)}
if $\vol\ci X$
is some volume measure compatible with
the smooth structure of $X$, then $n_\infty$ is absolutely continuous
with respect to $\vol\ci X$, and there exists a set $X'\subset X$ of full $\vol\ci X$ such that,
for every $x\in X'$, the quantity
\[
\mathfrak n(x) = \bar\nu_\pi (\pi^{-1}(x))\, \frac{\rd \pi_*\vol}{\rd \vol\ci X}(x)
\]
is well-defined and does not depend on the choice of the chart $\pi\colon U\to X$ with $x\in \pi (U)$. Moreover,
$$
\rd n_\infty  = \mathfrak n \, \rd \vol\ci X.
$$
\end{theorem}

The point of part (iii) is that,
for $\vol\ci X$-almost all $x\in X$, it allows one to compute the Radon-Nikodym derivative $\frac{\rd n_\infty}{\rd \vol\ci X}(x)$
using {\em any} chart containing $x$. In particular, nothing prevents us from choosing for each point its own individual chart.

\subsubsection{How to verify tameness?}\label{subsubsect:tameness}
Theorem~\ref{thm:manifold}, as stated,
has an essential shortcoming: it may be somewhat unpleasant to verify tameness of $(f\ci L)$
because formally it requires one to estimate various quantities in the local
coordinates given by $\pi$ for every chart $\pi\colon U\to X$, however weird or ugly. The next two observations (both
of purely technical nature) allow one to substantially reduce this workload. Recall that an atlas on $X$
is any family of charts $\bA=\bigl\{\pi_\alpha:U_\alpha\to X\bigr\}_\alpha$ such that $\bigcup_{\alpha}\pi_\alpha(U_\alpha)=X$.
Here is our first observation:

\begin{itemize}
\item{\em Suppose $\bA$ is an atlas on $X$ and that, for every chart $\pi_\alpha\in\bA$,
$(f\ci L \circ \pi_\alpha)$ is tame on $U_\alpha$. Then $(f\ci L)$ is tame on $X$}.
\end{itemize}

The possibility to check the tameness for the charts from any atlas
of our choice is quite a relief. However, one unpleasant thing still remains. It may (and often does)
happen that for every point $x\in X$ there is one ``preferred'' chart $\pi_x:U_x\to X$ covering $x$
such that the computations in this chart are a piece of cake in any infinitesimal neighborhood
of $x$ but not quite so even a bit away from $x$. In this case we would strongly prefer to compute
all quantities and check all conditions at $x$ using its preferred chart $\pi_x$. However, we are still
formally required to run the computations concurrently on any compact subset of any given chart using
the local coordinates given by that particular chart. Our next observation takes care of this
difficulty.

\begin{definition}
{\rm
We say that an atlas $\bA$ of $X$ has {\em uniformly bounded distortions} if there exists a constant $A>0$ such that
all partial derivatives of orders $\le 2$ of all coordinate functions of all
transition maps between the charts of $\bA$ are bounded by $A$.
}
\end{definition}
Note that this definition doesn't require $X$
to be uniform in any sense; rather it requires that the charts in $\bA$ be small enough so that
$X$ doesn't show any non-trivial structure within the union of each chart with all charts it intersects, and
that the chart scalings be more or less consistent with each other within small regions.
Our second observation says that
\begin{itemize}
\item
{\em If the atlas $\bA$ has uniformly bounded distortions, then to check the tameness of $(f\ci L)$ on $X$, it suffices
to check the relevant conditions and uniform bounds (on compact subsets of $X$) for the related quantities
computed in the charts $(U_x, \pi_x)$ at the points $\pi_x^{-1}(x)$ only}.
\end{itemize}

These two observations may be not obvious and we will explain them more in Section~\ref{sect:mflds}.

\medskip
Note that in our examples, we will  deal with compact manifolds admitting a transitive group
$\mathcal G$ of diffeomorphisms leaving the \pge\ $(f\ci L)$ under consideration invariant (meaning that
for each ${\mathbf g}\in \mathcal G$ and each $L$, the \cgf s $f_L$ and $f_L\circ {\mathbf g}$ have the same
distribution).
In such situation, all one needs is to find one chart $\pi:U\to X$ such that the atlas consisting of the charts
${\mathbf g}\circ \pi$,  ${\mathbf g}\in \mathcal G$, has uniformly bounded distortions. Then one may fix his/her favorite
point $x=\pi(u)$ in that chart, and establish all the required bounds and conditions at this single point for this single
chart. All passages about ``almost every $x$'' and suprema and infima over $Q$ in all conditions can be ignored in
such setup because all the related objects and quantities do not depend on $x$ at all.

\subsection{The final remarks about Theorems~\ref{thm:main_thm} and~\ref{thm:manifold}}\label{subsect:remarks-Thms2_and_3}

\subsubsection{}
Note that the particular choice of the counting measures $n_L$ plays no r\^{o}le.
The reason is that, for large $L$,  with high probability the overwhelming part of $n_L$ comes from components of
arbitrarily small diameter. Such components
can be viewed as single points at the macroscopic level.

\subsubsection{}
If the manifold $X$ is compact, we can apply the conclusion of Theorem~\ref{thm:manifold}
to $\phi\equiv 1$ and to obtain the asymptotics $(n_\infty(X)+o(1))L^m$ for the typical (and the mean)
total number of nodal components of $f\ci L$ on $X$ as $L\to\infty$. Of course, this asymptotic law
is really useful only when $n_\infty(X)>0$. Finding an asymptotic formula (or even a decent
estimate) for the variance of the total number of nodal components in such regimes remains
an open problem.

\subsubsection{}
The proof of Theorem~\ref{thm:main_thm}  also shows that the value $\bar\nu (x)$
can be recovered as a double-scaling limit. In
Lemma~\ref{lemma:local} we show that, for almost every $x\in U$ and for each $\e>0$,
we have
\[
\lim_{R\to\infty}\, \varlimsup_{L\to\infty} \cP \Bigl\{ \Big|
\frac{N\bigl( x, R/L; f\ci L \bigr)}{\vol B(R)} - \bar\nu (x) \Big| > \e  \Bigr\} =0
\]
where $ N\bigl( x, \tfrac{R}{L}; f\ci L \bigr) = N(R, f_{x, L})$ is the number of the connected components of
the zero set $Z(f\ci L)$ contained in the open ball centered at $x$ of radius $R/L$.

\subsubsection{}
A few words should be said about the measurability issues. While we prove every measurability result
that is necessary for the completeness of the formal exposition, when possible, we circumvent
this discussion by using upper integral and upper expectation instead of the usual ones.
Note that the Borel measurability of similar quantities
has been discussed in detail in Rozenshein's Master Thesis~\cite[Section~5]{YR}.

\subsection{Pertinent works}\label{subsect:pertinent-works}

\subsubsection{}
The earliest non-trivial lower bound for the mean number of connected components is,
probably, due to
Malevich. In~\cite{Mal}, she considered a $C^2$-smooth translation-invariant Gaussian random
function $F$ on $\bR^2$ with {\em positive} covariance kernel decaying at a certain rate at infinity. She
proved that $ \cE N(R; F)/R^2 $ is bounded from below and from above by two positive constants. Her
proof of the lower bound uses Slepian's inequality and probably cannot be immediately extended to models with covariance
kernels that change their signs.

\subsubsection{}
Several years ago, Bogomolny and Schmit~\cite{BS} proposed a bond percolation model
for the description of the zero set of the translation-invariant Gaussian function $F$ on $\bR^2$
whose spectral measure is the Lebesgue measure on the unit circumference. This model completely
ignores slowly decaying correlations between values of the random function at different points and
is very far from being rigorous. The predictions of Bogomolny and Schmit were checked by computational
experiments carried out by Nastasescu~\cite{Nast}, Konrad~\cite{Konrad}, and Beliaev and Kereta~\cite{BK}.
The observed value of the constant $\nu$ was very close to but still noticeably less than the
Bogomolny and Schmit prediction. It would be very interesting to reveal a hidden ``universality law'' that
provides the rigorous foundation for the work done by Bogomolny and Schmit.
Note also that it is not clear whether or to what degree their approach can be extended to make reasonably accurate
predictions about the behavior of nodal components of translation-invariant Gaussian functions corresponding
to other spectral measures in $\bR^2$ or in dimensions $m>2$.

\subsubsection{}
In~\cite{NS}, we showed that for the Gaussian ensemble of spherical harmonics of large
degree $L$ on the two-dimensional sphere, the total number $N(f_L)$ of connected components of $Z(f_L)$ satisfies
\[ \cP \left\{ \left| L^{-2} N(f\ci L) -
\upsilon \right| >\e\right\} < C(\e) e^{-c(\e) \dim \cH_L}\,,
\]
with some $\upsilon > 0$. The
limiting function for this ensemble is the one considered by Bogomolny and Schmit.
The case of higher dimension (in a slightly different setting) was treated by Rozenshein in~\cite{YR}.
The exponential concentration of $N(f\ci L)/L^2$ is interesting since this model
has slowly decaying correlations.

We were unable to prove the exponential
concentration for other ensembles considered here. The difficulty is caused by the
small components, which do not exist when $f\ci L$ is an eigenfunction of the Laplacian. Even in the
univariate case, the question about the exponential concentration in
Theorem~\ref{thm:Euclid-version} remains open; cf. Tsirelson's lecture
notes~\cite{Tsirelson}.

Some lower bounds for
the number of connected components of the zero set and for other similar quantities were obtained in different settings
by Bourgain and Rudnick~\cite{BR}, Fyodorov, Lerario, Lundberg~\cite{FLL}, Gayet and Welschinger~\cite{GW1, GW2, GW3}, Lerario and Lundberg~\cite{LL}
using the ``barrier construction'' from~\cite{NS}.

\subsubsection{}
Certain versions of main results of this  work were
presented at the St. Petersburg Summer School in Probability and Statistical Physics
(June, 2012) and appeared in the lecture notes~\cite{S}.

\subsubsection{}
There have been several works of interest relying on ideas and techniques developed in this paper,
among which those by Bourgain~\cite{Bourgain}, Canzani and Sarnak~\cite{CS}, Kurlberg and Wigman~\cite{KW} and
Sarnak and Wigman~\cite{SW} deserve special attention of the reader.

\bigskip\par\noindent{\bf Acknowledgments}. On many occasions, Boris Tsirelson helped us by providing
information and references concerning Gaussian measures and measurability. Leonid Polterovich and
Ze\'ev Rudnick have read parts of a preliminary version of this work and made several valuable comments,
which we took into account. We have had encouraging discussions of this work with Andrei
Okounkov, Peter Sarnak, and Jean-Yves Welschinger. Alex Barnett and Maria Nastasescu showed us the
beautiful and inspiring simulations. We thank them all.

\section{Examples}\label{sect:Examples}

Here, we point out two examples illustrating
Theorem~\ref{thm:manifold}. In our examples, the manifold $X$ has a natural
Riemannian metric and a transitive group of isometric diffeomorphisms that
leaves the distribution of $(f\ci L)$ invariant. As discussed near the
end of Section~\ref{subsubsect:tameness}, this will allow us to check
the conditions of Theorem~\ref{thm:manifold} at just one point $x\in X$ with
respect to a natural local chart associated with this point and to conclude
that the limiting measure $n_\infty$ on the manifold $X$ is a constant multiple
of the Riemannian volume on $X$. Moreover, since in our examples
the kernels $K_{x,L}(u,v)$ converge to $k(u-v)$ uniformly with all derivatives
on compact subsets of $\bR^{m}\times \bR^m$, the uniform smoothness
and non-degeneracy of the kernels $K_{x,L}(u, v)$ can be derived from the corresponding
properties of the limiting kernel $k(u-v)$. Passing to the limit in our
examples is an elementary exercise in Taylor calculus and complex analysis.
This list
of examples may be continued (see~\cite{CS, KW, YR, SW}) but the two ones we included into this paper
should be already enough to convey its main message, which is

\begin{itemize}
\item
{\em Under not unreasonably unfavorable conditions, establishing the
asymptotics for the number of nodal domains for \pge s is about as easy (or,
if the reader prefers, as hard) as establishing the convergence of the
scaled kernels and investigating the resulting limiting processes.}
\end{itemize}

\subsection{Trigonometric ensemble}\label{subsect:trig-ensemble}
Here $\cH_{n}$  is the subspace of $L^2(\bT^m)$ that
consists of real-valued trigonometric polynomials
\[
{\rm Re}\, \sum_{\nu\in\bZ^m\colon |\nu|_\infty\le n} c_\nu e^{2\pi {\rm i}(\nu \cdot
x)}
\]
in $m$ variables of degree $\le n$ in each variable.
A straightforward computation shows that the corresponding normalized covariance kernel
coincides with the product of $m$ Dirichlet's kernels:
\[
K_{n} (x,y) = \prod_{j=1}^m \frac{\sin \left[ \pi (2n+1)(x_j-y_j)
\right]}{(2n+1)\sin \left[ \pi (x_j-y_j) \right]}\,.
\]
In this case, it is natural to choose the degree $n$ as the scaling parameter $L$.
The scaled kernels $K_{x,n}(u, v) = K_{n}(x+n^{-1}u, x+n^{-1}v)$
do not depend on the choice of the point $x\in\bT^m$.
They extend analytically from $\bR^m\times \bR^m$ to $\bC^m \times \bC^m$ and the extensions converge uniformly on compact subsets of $\bC^m \times \bC^m$
to
\[
\prod_{j=1}^m \frac{\sin 2\pi (u_j-v_j)}{2\pi (u_j-v_j)}\,.
\]
This implies the convergence with all derivatives on all compact subsets
of $\bR^m \times \bR^m$. The limiting spectral measure $\rho$ is the normalized
Lebesgue measure on the cube $[-1, 1]^m\subset \bR^m$.

\subsection{Kostlan's ensemble}\label{subsect:Kostlan}
In this case, $\cH_n$ is the space of the homogeneous real-valued polynomials
of degree $n$ in $m+1$ variables restricted to the unit sphere $\bS^m$. The scalar
product in $\cH_n$ is given by
\begin{equation}\label{eq:sc_product}
\langle f, g\rangle = \sum_{|J|=n} { {n}\choose{J} }^{-1} f_J
g_J
\end{equation}
where
\[
f(X) = \sum_{|J|=n} f_J X^J, \quad g(X) = \sum_{|J|=n} g_J X^J, \qquad
X^J = x_0^{j_0}x_1^{j_1} x_2^{j_2} \ldots x_{m}^{j_{m}},
\]
and
\[ J=(j_0, j_1, j_2, \ldots ,
j_{m}), \quad |J|=j_0+j_1+j_2+ \ldots +j_{m}, \quad { {n}\choose{J} } = \frac{n!}{j_0!j_1!
j_2! \ldots j_{m}!}\,.
\]
The form of the scalar product~\eqref{eq:sc_product} comes from the
complexification: extending the homogeneous polynomials $f$ and $g$ to $\bC^{m+1}$,
one can show that
\[
\langle f , g \rangle_{\cH_n} =  c(n, m) \int_{\bC^{m+1}} f(Z) \overline{g(Z)} e^{-|Z|^2}\, {\rm
d}\vol (Z)\,,
\]
i.e., $\langle f, g \rangle_{\cH_n}$ coincides (up to a positive factor) with the scalar
product in the Fock-Bargmann space (or any other weighted $L^2$-space of entire
functions with fast decaying radial weight).

It is known that the complexified Kostlan ensemble is {\em the only unitarily
invariant} Gaussian ensemble of homogeneous polynomials.
On the other hand, there are many other
orthogonally invariant Gaussian ensembles, all of them having been classified by
Kostlan~\cite{Kostlan2} (see~\cite[Section~2]{FLL} for some details).

The normalized covariance kernel of Kostlan's ensemble equals
$ (x \cdot y)^n $.
Take $x=(0,\, \ldots \, , 0, 1)$ (the ``North Pole'') and consider the local chart
$\pi (u) = \bigl( u, \sqrt{1-|u|^2} \bigr)$ where $u$ runs over a small neighbourhood of the
origin in $\bR^m$.
Then
\begin{align*}
\pi(u)\cdot\pi(v)&=\sum_{j=1}^m u_j v_j + \Bigl( 1-\sum_{j=1}^m u_j^2 \Bigr)^{\frac12} \cdot
\Bigl(1-\sum_{j=1}^m v_j^2\Bigr)^{\frac12}
\\
&=1-\frac12\sum_{j=1}^m(u_j-v_j)^2+O\bigl(|u|^4+|v|^4\bigr)\text{ as }u,v\to 0\,.
\end{align*}

This suggests that the correct scaling in this case is $L=\sqrt{n}$ and the limiting
covariance kernel is
\begin{align*}
\lim_{n\to\infty} \Bigl( \pi(n^{-\frac12}u) \cdot \pi(n^{-\frac12}v) ) \Bigr)^n
&= \lim_{n\to\infty} \Bigl(  1-(2n)^{-1}\sum_{j=1}^m(u_j-v_j)^2+n^{-2}O\bigl(|u|^4+|v|^4\bigr)\Bigr)^n
\\
&= \exp\Bigl\{-\frac12 \sum_{j=1}^m(u_j-v_j)^2\Bigr\}\,.
\end{align*}
The justification of the local uniform convergence with all derivatives is similar to that in the previous example,
and we skip it. The limiting spectral measure is the Gaussian measure on $\bR^m$ with the density $c_m
e^{-2\pi^2 |\la|^2} $.

An interesting feature of this example
is a very rapid off diagonal decay of the covariance kernel.

\section{Notation}\label{sect:Notation}

We denote by $B(x, r)$ the open ball of radius
$r$ centered at $x$, $\bar B(x, r)$ denotes the corresponding closed ball.
$B(r)$ always denotes the open ball of radius $r$ centered at the origin.

For a closed set $\Gamma\subset \bR^m$, we denote by $N(x, r; \Gamma)$
the number of the connected components of $\Gamma$ that are contained
in the open ball $B(x, r)$, and by $N^*(x, r; \Gamma)$ the number of
the connected components of $\Gamma $ that intersect the closed ball $\bar B(x, r)$.
If $\Gamma=Z(f)$ is the zero set of a continuous function $f$, we will abuse
the notation slightly and write $N(x, r; f)$ instead of $N(x, r; Z(f))$.
For a bounded open convex set $S$ and $R>0$, we denote by $N_S(R; \Gamma)$ the number of
connected components of $\Gamma$ that are contained in $S(R)=\{u\colon R^{-1}u\in S\}$.

\smallskip Throughout the paper, we denote by $c$ and $C$ various positive constants, which may
depend on the dimension $m$ and on the parameters of the Gaussian process or
ensemble under consideration (the parameters in the conditions of
Theorems~\ref{thm:Euclid-version} and~\ref{thm:main_thm}) but on nothing else.
The values of these constants may
vary from line to line. Usually, the constants denoted by $C$ should be thought of as large,
and the constants denoted by $c$ as small.
The notation $a \lesssim b$ means that $a \le C\cdot b$.

\medskip
Quite frequently, we will use the smoothness class $C^{2-}(U)$ ($U\subset\bR^m$ is an open set), which
we define as
\[
C^{2-}(U) = \bigcap_{0<\beta<1}\, C^{1+\beta}(U)\,.
\]
Recall that to check that $g\in C^{1+\beta}(U)$ it suffices to show that
$g\in C^1(U)$ and the first order partial derivatives $\partial_{x_i} g$ are
$\beta$-H\"older functions on any closed ball $\bar B\subset U$.

\section{Lemmata}\label{sect:lemmata}

In this section, we present several lemmas needed for the proofs of
Theorems~\ref{thm:Euclid-version} and~\ref{thm:main_thm}.

\subsection{Some integral geometry}\label{subsect:lemma-integr-geom}

The first result is taken from~\cite[Claim~5.1]{NS} where it appears in a slightly different form.
\begin{lemma}\label{lemma:integr-geom}
Suppose $ \Gamma\subset\bR^m $ is a closed set and $S \supset B(1)$ is a bounded open convex set.
Then, for $0<r<R$,
\[
\int_{S(R-r)} \frac{N(u, r; \Gamma)}{\vol B(r)}\, {\rm d} \vol (u) \le N_S(R; \Gamma)
\le \int_{S(R+r)} \frac{N^*(u, r; \Gamma)}{\vol B(r)}\, {\rm d}\vol (u)\,.
\]
\end{lemma}

\medskip
Note that $N(u, r; \Gamma)$ is lower semicontinuous as a function of $u$.
Proving the Lebesgue measurability of $u\mapsto N^*(u, r; \Gamma)$ without
additional assumptions on $\Gamma$ may be somewhat nontrivial. However, we will
apply this lemma only in the case when the set of connected components of
$\Gamma$ is countable. Also, replacing the integral on the RHS by the upper
Lebesgue integral will not affect the argument in any way. So, we will not
dwell on this particular measurability.

\medskip\noindent{\em Proof:} For a connected component $\gamma$ of $\Gamma$, we put
\[
G_*(\gamma) = \bigcap_{y\in\gamma} B(y, r), \quad G^*(\gamma) = \bigcup_{y\in\gamma} \bar B(y,
r)\,.
\]
Note that since $\gamma$ is closed, $G_*(\gamma)$ is open and $G^*(\gamma)$ is closed.
Also, for any $y\in\gamma$, $ G_*(\gamma) \subset B(y, r) \subset G^*(\gamma) $.
Hence,
\begin{align*}
\int_{S(R-r)} N(u, r; \Gamma)\, \rd\vol (u)
&= \int_{S(R-r)} \Bigl( \sum_{\gamma\colon \gamma\subset B(u, r)} 1 \Bigr) \rd\vol (u)
\\
&\le \int_{S(R-r)} \Bigl( \sum_{\gamma\colon \gamma\subset S(R),\, u\in G_*(\gamma)} 1 \Bigr)
\rd\vol (u) \\[7pt]
&= \sum_{\gamma\subset S(R)} \vol \bigl( G_*(\gamma)\cap S(R-r) \bigr) \\[10pt]
&\le N_S(R; \Gamma) \cdot \vol B(r)\,,
\end{align*}
proving the left inequality.

On the other hand,
\begin{align*}
\int_{S(R+r)} N^*(u, r; \Gamma)\, \rd\vol (u)
&= \int_{S(R+r)} \Bigl( \sum_{\gamma\colon u\in G^*(\gamma)} 1 \Bigr)\, \rd\vol(u) \\[7pt]
&= \sum_\gamma \vol \bigl( G^*(\gamma) \cap S(R+r) \bigr)\,.
\end{align*}
Since for every connected component $\gamma$ having a common point $y$ with $S(R)$, we have
$B(y, r) \subset G^*(\gamma)\cap S(R+r)$, the last sum is at least
$N_S(R; \Gamma) \cdot \vol B(r)$, so the right inequality holds as well.
\hfill $\Box$

\subsection{Stability of components of the zero set under small perturbations}\label{subsect:lemma-stability}

If zero is not a critical value of a smooth function then the zero set of this function
is stable under small
perturbations. The following lemma, which quantifies this general principle, is taken
from~\cite[Claim~4.2]{NS} where it was proven in the two-dimensional case. The proof of the general
case needs no changes.

Denote by $V_{+t}$ the open $t$-neighbourhood of a set $V\subset\bR^m$.
\begin{lemma}\label{lemma:Claim4.2}
Fix $\alpha, \beta >0$. Let $F$ be a $C^1$-smooth
function on an open ball $B\subset \bR^m$ such that at every point $u\in B$, either
$|F(u)|>\alpha$, or $|\nabla F(u)|>\beta$. Then each component $\gamma$ of the zero set $Z(F)$ with
${\rm dist}(\gamma, \partial B) > \alpha/\beta$ is contained in an open ``annulus''
$A_\gamma\subset \gamma_{+\alpha/\beta}$ bounded by two smooth connected hypersufaces such that
$ F=+\alpha$ on one boundary component of $A_\gamma$, and $F=-\alpha$ on the other one. Furthermore,
the ``annuli'' $A_\gamma$ are pairwise disjoint.
\end{lemma}
As an immediate corollary, we obtain
\begin{lemma}\label{lemma:Cor4.3}
Under the assumptions of the previous lemma, suppose that $G\in
C(B)$ with $\sup_B |G|<\alpha$. Then each component $\gamma$ of $Z(F)$ with ${\rm dist}(\gamma,
\partial B) > \alpha/\beta$ generates a component $\widetilde\gamma$ of the zero set $Z(F+G)$ such
that $\widetilde\gamma\subset\gamma_{+\alpha/\beta}$
and different components $\gamma_1\ne \gamma_2$ of $Z(F)$ generate different components
$\widetilde\gamma_1 \ne \widetilde\gamma_2$ of $Z(F+G)$.
\end{lemma}

\subsection{Statistical independence of $g$ and $\nabla g$}\label{subsect:lemma-independence}
Quite often, we will use the following
well-known fact:

\begin{lemma}\label{lemma:indep}
Suppose $U\subset\bR^n$ is an open set and
$g\colon U\to \bR$ is a Gaussian $C^1$-function on $U$ that has constant
variance. Then $g(u)$ and its gradient $\nabla g(u)$ are independent
for every $u\in U$.
\end{lemma}

\par\noindent{\em Proof:} Denote by $g_{u_i}$ the partial derivative  $\partial_{u_i} g$. The
covariance kernel $K(u,v)=\cE\bigl\{ g(u)g(v) \bigr\}$ is a $C^1$-function,
and $ \cE\bigl\{ g_{u_i}(u) g(u) \bigr\} = K_{u_i}(u, v)\big|_{v=u} $.
Since the function $u'\mapsto K(u', u)$ attains its maximal value at
$u'=u$ and is $C^1$-smooth, we have $K_{u_i}(u, v)\big|_{v=u}=0$.
Therefore, $ \cE\bigl\{ g_{u_i}(u) g(u) \bigr\} = 0$.
Since $g(u)$ and $\nabla g(u)$ are jointly Gaussian, this orthogonality implies
their independence. \hfill $\Box$

\medskip We will be using the following corollary:

\begin{lemma}\label{lemma:non-degen}
Suppose $F\colon \bR^m\to \bR$ is a Gaussian random function
with translation-invariant distribution whose spectral measure $\rho$ satisfies conditions
($\rho$1) and ($\rho$3). Then the distribution of the Gaussian vector $\bigl( F(u), \nabla F(u)
\bigr)$ does not degenerate.
\end{lemma}
\noindent{\em Proof of Lemma~\ref{lemma:non-degen}:}
By Lemma~\ref{lemma:indep},
$F(u)$ and $\nabla F(u)$ are independent. Hence, it suffices to show that
the distribution of $\nabla F(u)$ does not degenerate. If it degenerates, then there exists a
non-zero vector $v\in\bR^m$ such that
\[ 0 = \cE \bigl\{  ( v \cdot \nabla F )^2 \bigr\} = 4\pi^2\, \int_{\bR^m}
( v \cdot \la )^2 \, {\rm d}\rho (\la)\,,
\]
which is impossible since, due to condition
($\rho$3), the spectral measure $\rho$ cannot be supported on a linear hyperplane.
\hfill $\Box$

\section{Quantitative versions of Bulinskaya's lemma}\label{sect:Bulinskaya}

\subsection{Preliminaries}\label{subsect:Bulinsk-prelim}

The purpose of this part is to show that certain ``bad events'' have negligibly small probability.
The particular bad events we want to get rid of are the event that the random Gaussian
function and its gradient are simultaneously small at some point and the event that $Z(f)$ has too
many connected components.

\medskip
Everywhere in this part, $B_R\subset \bR^m$ is a fixed ball of large radius $R>1$, $S=\partial B_R$,
$U$ is an open neighbourhood of $B_{R+1}$,
and $f$ is a continuous Gaussian function on $U$
with the covariance kernel $K$. As usual, we will assume that
the function $f$ is normalized, that is, $\cE |f(x)|^2 = K(x, x) =1 $, $x\in U$.
We will impose certain bounds on the smoothness and non-degeneracy. These bounds are normalized versions
of estimates used in the definition of controllability of parametric Gaussian ensembles. Namely, we will assume
that

\smallskip\par\noindent (i) the kernel $K$ is $C^{2,2}(U\times U)$-smooth and
\[
\max_{|\alpha|\le 2}\, \max_{x\in \bar B_{R+1}}\, \bigl| \partial_x^\alpha\, \partial_y^\alpha K(x, y)|_{y=x} \bigr| \le M  < \infty\,,
\]

\noindent and that

\smallskip\par\noindent (ii) the process $f$ is non-degenerate on $U$ and
\[
\inf_{x\in\bar B_{R+1}} \det C_x \ge \kappa >0\,,
\]
where $C_x$ is the covariance matrix of the Gaussian random vector $\nabla f(x)$, that is, the matrix with the entries
$C_x(i, j) = \partial_{x_i} \partial_{y_j} K(x, y)|_{y=x}$.

\begin{itemize}
\item
{\em Till the end of Section~\ref{sect:Bulinskaya}, the constants $M$ and $\kappa$ remain fixed and all the constants
that appear in the  conclusions of all results proven here may depend on $M$ and $\kappa$}.
\end{itemize}

\medskip As shown in Appendix~\ref{A_subsect:remarks_Kolmogorov}, the smoothness assumption (i) yields that, almost surely, the process $f$ is
$C^{2-}(U)$-smooth. We will be frequently using a quantitative version of this statement, which is also given in~\ref{A_subsect:remarks_Kolmogorov}.
For a closed ball $\bar B\subset U$, denote by $\| f \|_{\bar B, 1+\beta}$
the least $N$ such that
\[
\max_{\bar B} |f| \le N, \quad \max_{\bar B} |\nabla f| \le N, \quad{\rm  and}\quad  |\nabla f(x)-\nabla f(y)|
\le N|x-y|^\beta \ \  {\rm for\ }\  x, y \in\bar B.
\]
Then, for every $\beta<1$ and every $p<\infty$,
\[
\sup_{x\in \bar B_{R}} \cE \bigl\{ \| f \|_{\bar B, 1+\beta}^p \bigr\} \le C (\beta, p, M) < \infty\,.
\]

\subsection{The function $\Phi$}\label{subsect:Phi}

A prominent r\^{o}le in our approach will be played by the function
\[
\Phi (x) = |f(x)|^{-t} |\nabla f(x)|^{-tm}\,, \qquad x\in U,\ t\in (0, 1)\,,
\]
and by its spherical version
\[
\Phi_S (x) = |f(x)|^{-t} |\nabla_S f(x)|^{-t(m-1)}\,, \qquad x\in S=\partial B_R,\ t\in (0, 1)\,,
\]
where $\nabla _S f(x)$ is the projection of
the vector $\nabla f(x)$ to the tangent space to $S$ at the point $x\in S$. The main feature of this function
is that if $f$ and $\nabla f$ (or $\nabla_S f$) are very small at two points that are close to each other (in particular,
if they are small at the same point), then $\Phi$ ($\Phi_S$ correspondingly) is very large in a neighbourhood
of these two points. At the same time, since $f$ is normalized and $\nabla f$ is non-degenerate,
\begin{itemize}\item
{\em the moments $ \cE \bigl\{ \Phi^q (x) \bigr\} $ and $ \cE \bigl\{ \Phi_S^q (x) \bigr\} $
are bounded locally uniformly on $U$ and uniformly on $S$ whenever we fix $t<1<q$ so that $tq<1$.
Moreover,
if $t$ and $q$ satisfying
this restrictions are fixed, the suprema
$ \sup_{\bar B_{R+1}} \cE \bigl\{ \Phi^q  \bigr\} $ and  $ \sup_S \cE \bigl\{ \Phi_S^q  \bigr\} $
are bounded by constants depending only on $\kappa$}.
\end{itemize}
\subsection{Almost surely, zero is not a critical value of $f$}\label{subsect:Bulinskaya-original-form}

As a warm up, we prove a useful qualitative result that goes back to Bulinskaya.
\begin{lemma}\label{lemma:Bulinskaya}\mbox{} Almost surely, the following assertions hold:

\smallskip\par\noindent {\rm (i)} zero is not a critical value of $f$;

\smallskip\par\noindent {\rm (ii)}
there is no point $z\in S \cap Z(f)$ at which $\nabla_S f(z)=0$.
\end{lemma}
\noindent{\em Proof}:
In the first case, we use the function $ \Phi $.
Fix a compact set $Q\subset U$ and take a positive $\delta<{\rm dist}(Q, \partial U)$.
Consider the event
\[
\Omega_Q = \bigl\{\exists z\in Q\colon {\rm such\ that\ } f(z)= 0, \ \nabla f(z)=0 \bigr\}
\]
and take a ball $\bar B\subset U$ centered at $z$ of radius
less than $\delta$.
Since the function $\nabla f(x)$ is $\beta$-H\"older with every $\beta<1$, we have, for
all $x\in\bar B$,
\[
|f(x)| \lesssim |x-z|, \quad |\nabla f(x)| \lesssim |x-z|^\beta\,,
\]
whence $\Phi(x) \gtrsim |x-z|^{-t(1+\beta m)}$. Hence,  choosing $t$ and $\beta$ so close to $1$ that
$t(1+\beta m) > m$, we see that
\[
\int_{Q_{+\delta}} \Phi\,\rd\vol \ge \int_B \Phi\,\rd\vol =+\infty\,.
\]
Recalling that $\cE \{ \Phi (x) \}$ is uniformly bounded on $\bar Q_{+\delta}$ and using Fubini's theorem,
we conclude that the event $\Omega_Q$ has zero probability. It remains to note that $U$ can be covered by countably
many compact subsets.

Similarly, in the second case we take $ \Phi_S $.
As above, the expectation $\cE \{ \Phi_S (x)\}$ is uniformly bounded on
$S$. Suppose that, for some $z\in S \cap Z(f)$,
$\nabla_S f(z) = 0$, that is, the gradient $\nabla f(z)$ is orthogonal to the sphere $S$.
Then, for $x\in S$, we have
\[
|f(x)| \lesssim |x-z|, \qquad |\nabla_S f(x)| \lesssim |x-z|^\beta + R^{-1}|x-z| \lesssim  |x-z|^\beta\,,
\]
and, thereby, $\Phi_S(x) \gtrsim |x-z|^{-t(1+\beta(m-1))}$. Therefore, choosing $t$ and $\beta$ so close to $1$ that $\nolinebreak t(1+\beta(m-1)) >m-1$, we
get
\[
\int_S \Phi\,\rd\vol_{S} = +\infty,
\]
and conclude that the event we consider has zero probability. \hfill $\Box$

\subsection{With probability close to one, $f$ and $\nabla f$ cannot be simultaneously small}\label{subsect:Bulinskaya-quantitative}

Here, we prove a quantitative version of Lemma~\ref{lemma:Bulinskaya}.
\begin{lemma}\label{lemma:prob-stable}
Given $\delta>0$, there exists $\tau>0$ (possibly, depending on $R$)
such that
\[
\cP\bigl\{ \min_{x\in \bar B_R}\, \max\{ |f(x)|, |\nabla f(x)| \} < \tau \bigr\} < \delta\,.
\]
\end{lemma}
\noindent{\em Proof}: Denote by $\Omega_\tau$ the event
\[
\bigl\{ \exists z\in \bar B_R\colon |f(z)|, |\nabla f(z)|<\tau \bigr\}
\]
and put
\[
W = 1 + \| f \|_{\bar B_{R+1}, 1+\beta}\,.
\]
The parameter $\beta\in (0, 1)$ will be specified later. If the event $\Omega_\tau$ occurs, then
in the ball $B = B(z, \tau)$ with $\tau\in (0, 1)$, we have
\[
|f(x)| \le \tau + \tau \| f \|_{\bar B_{R+1}, 1+\beta} = W \tau\,,
\]
and
\[
| \nabla f(x) | \le \tau + \tau^\beta \| f \|_{\bar B_{R+1}, 1+\beta} < W \tau^\beta\,.
\]
Then, on $\Omega_\tau$,
\[
\Phi (x) \ge \tau^{-t(1+\beta m)} W^{-t(1+m)}\,, \qquad {\rm for}\ x\in B
\]
and
\[
\int_{B_{R+1}} \Phi\,\rd\vol \ge \int_{B} \Phi\, \rd\vol \ge c \tau^{m-t(1+\beta m)} W^{-t(1+m)}\,.
\]
Therefore,
\begin{align*}
\cP\bigl\{ \Omega_\tau \bigr\} &\le C \tau^{t(1+\beta m)-m}\, \vol (B_{R+1})\, \cE \Bigl\{ W^{t(1+m)}\, \frac1{\vol (B_{R+1})}\, \int_{B_{R+1}} \Phi\,\rd\vol \Bigr\}\\
&\le  C \tau^{t(1+\beta m)-m} \vol (B_{R+1}) \Bigl( \cE \bigl\{  W^{pt(1+m)} \bigr\} \Bigr)^{\frac1p}\,
\Bigl( \frac1{\vol (B_{R+1})}\, \int_{B_{R+1}} \cE \bigl\{ \Phi^q \} \rd\vol \Bigr)^{\frac1q}\,,
\end{align*}
with $ \frac1p + \frac1q = 1$.
The only restriction we have is $\beta<1<q<\frac1t$. So we can take $\beta$ and $t$ so close to $1$ that the
exponent $t(1+\beta m) - m = t-m(1-\beta)$ remains positive. This completes the proof. \hfill $\Box$

\subsection{General principle for estimating the number of connected components}\label{subsect:general-principle}

Our next aim is to estimate how many connected components of various kinds $Z(f)$ may have. We start with
``an abstract scheme'', which our estimates will be based on.

\medskip
Let $(X, \mu)$ be a measure space with $0<\mu (X)<\infty$, and let $X=\bigcup_j X_j$ be a cover of $X$ with
bounded covering number $C_0$ (that is, for every $x\in X$, $\#\bigl\{ j\colon x\in X_j \bigr\}\le C_0$).
Let $(\Omega, \cP)$ be a probability space, and let $\bigl\{ (Y_i(\omega), z_i(\omega) ) \bigr\}_{1\le i \le N(\omega)}$ be
disjoint subsets of $X$ with marked points $z_i\in Y_i$ depending on the parameter $\omega\in\Omega$.
Our aim is to estimate the cardinality $N(\omega)$ of the collection $\{ Y_i \}$.

\begin{lemma}\label{lemma:abstract}
Let $\Phi\colon X\to\bR_+$ be a random function such that, for some $q>1$,
\[
\sup_X \cE\{\Phi^q\} < \infty\,.
\]
Let $\{ W_j \}$ be non-negative random variables such that, for any $p<\infty$,
\[
\sup_j \cE \{ W_j^p \} < \infty\,.
\]
Suppose that, for every pair $(i, j)$ with  $z_i\in X_j$, we have
\[
\int_{X_j\cap Y_i} \Phi\,\rd\mu  \ge \rho\, \mu (X_j\cap Y_i)^{-\sigma} W_j^{-\eta}
\]
with some $\rho, \sigma, \eta>0$.
Then
\[
\cE^* \{ N^q \}\le (C_0 C(\rho, \sigma))^q \, \mu (X)^q\, \Bigl[\, \sup_{X} \cE \bigl\{ \Phi^q \bigr\} \Bigr]^{\frac1{1+\sigma}} \cdot
\Bigl[ \sup_j \cE \bigl\{ W_j^{\frac{q\eta}{\sigma}} \bigr\} \Bigr]^{\frac{\sigma}{1+\sigma}}\,.
\]
\end{lemma}

\noindent{\em Proof of Lemma~\ref{lemma:abstract}}:
Let $ N_j $  be the number of $i$'s such that $ z_i\in X_j $.
Then, for at least $\frac12 N_j$ indices $i$ with this property, we have
$ \mu (X_j \cap Y_i) \le \frac2{N_j}\, \mu (X_j) $, and therefore,
\[
\int_{X_j\cap Y_i} \Phi\,\rd\mu \ge c(\rho, \sigma) N_j^\sigma \mu (X_j)^{-\sigma} W_j^{-\eta}\,,
\]
whence,
\[
\int_{X_j} \Phi\,\rd\mu \ge c(\rho, \sigma) N_j^{1+\sigma} \mu (X_j)^{-\sigma} W_j^{-\eta}\,.
\]
Applying H\"older's inequality with the exponents $1+\sigma$ and $\frac{1+\sigma}\sigma$, we get
\begin{align*}
N \le \sum_j N_j &= \sum_j N_j \mu (X_j)^{-\frac{\sigma}{1+\sigma}} W_j^{-\frac{\eta}{1+\sigma}} \cdot
\mu (X_j)^{\frac{\sigma}{1+\sigma}} W_j^{\frac{\eta}{1+\sigma}} \\
&\le \Bigl( \sum_j N_j^{1+\sigma} \mu (X_j)^{-\sigma} W_j^{-\eta} \Bigr)^{\frac1{1+\sigma}}
\cdot \Bigl( \sum_j \mu (X_j) W_j^{\frac{\eta}\sigma} \Bigr)^{\frac{\sigma}{1+\sigma}} \\
&\le C(\rho, \sigma) \Bigl( \int_{X} \Phi\,\rd\mu \Bigr)^{\frac1{1+\sigma}}
\cdot  \Bigl( \sum_j \mu (X_j) W_j^{\frac{\eta}\sigma} \Bigr)^{\frac{\sigma}{1+\sigma}}\,.
\end{align*}
Then
\[
\cE^* \bigl\{ N^q \bigr\} \le C(\rho, \sigma)^q\, \cE \Bigl\{ \Bigl( \int_{X} \Phi\,\rd\mu \Bigr)^{\frac{q}{1+\sigma}} \cdot
\Bigl( \sum_j \mu (X_j) W_j^{\frac{\eta}\sigma} \Bigr)^{\frac{q\sigma}{1+\sigma}} \Bigr\}
\]
and, applying H\"older's inequality with the same exponents again, we obtain
\[
\cE^* \bigl\{ N^q \bigr\}
\le C(\rho, \sigma)^q  \Bigl[ \cE \Bigl\{ \Bigl( \int_{X} \Phi\,\rd\mu \Bigr)^q \Bigr\} \Bigr]^{\frac1{1+\sigma}} \cdot
\Bigl[ \cE \Bigl\{ \Bigl( \sum_j \mu (X_j) W_j^{\frac{\eta}\sigma} \Bigr)^q \Bigr\} \Bigr]^{\frac{\sigma}{1+\sigma}}\,.
\]
At last, using H\"older's inequality with the exponents $\tfrac{q}{q-1}$ and $q$, we get
\[
\cE \Bigl\{ \Bigl( \int_{X} \Phi\,\rd\mu \Bigr)^q \Bigr\}
\le \mu (X)^{q-1}\, \int_X \cE\bigl\{ \Phi^q \bigr\}\, {\rm d}\mu \le \mu (X)^q\, \sup_X \cE\bigl\{ \Phi^q \bigr\}
\]
and
\begin{align*}
\cE \Bigl\{ \Bigl( \sum_j \mu (X_j) W_j^{\frac{\eta}\sigma} \Bigr)^q \Bigr\}
&= \cE \Bigl\{ \Bigl( \sum_j \mu (X_j)^{1-\frac1q} \mu(X_j)^\frac1q\,  W_j^{\frac{\eta}\sigma} \Bigr)^q \Bigr\} \\
&\le \Bigl[ \sum_j \mu (X_j) \Bigr]^{q-1}   \cdot \Bigl[ \sum_j \mu (X_j)\, \cE \bigl\{  W_j^{\frac{q\eta}\sigma}  \bigr\} \Bigr] \\
&\le \bigl( C_0 \mu (X) \bigr)^q \, \sup_j \cE \bigl\{  W_j^{\frac{q\eta}\sigma}  \bigr\}\,.
\end{align*}
Finally,
\begin{multline*}
\cE^*\{ N^q \} \le C(\rho, \sigma)^q
\Bigl[ \mu (X)^{q}\, \sup_{X} \cE\{ \Phi^q \} \Bigr]^{\frac1{1+\sigma}} \cdot
\Bigl[ \bigl( C_0 \mu (X) \bigr)^q\, \sup_j \cE\bigl\{ W_j^{\frac{q\eta}{\sigma}} \bigr\} \Bigr]^{\frac{\sigma}{1+\sigma}} \\
\le (C_0 C(\rho, \sigma))^q\,  \mu (X)^q\, \Bigl[\, \sup_{x\in X} \cE \bigl\{ \Phi (x)^q \bigr\} \Bigr]^{\frac1{1+\sigma}} \cdot
\Bigl[ \sup_j \cE \bigl\{ W_j^{\frac{q\eta}{\sigma}} \bigr\} \Bigr]^{\frac{\sigma}{1+\sigma}}\,,
\end{multline*}
completing the proof. \hfill $\Box$

\subsection{Components on a sphere}\label{subsect:components-on-sphere}

For the sphere $S =  \partial B_{R}$, we denote by $\mathfrak N(S; f)$ the number of connected components of
$S \setminus Z(f)$.

\begin{lemma}\label{lemma:sphere}
There are positive constants $C<\infty$ and $q>1$ such that
\[
\cE^* \{ \mathfrak N^q(S; f) \} \le C R^{q(m-1)}.
\]
\end{lemma}

\noindent{\em Proof of Lemma~\ref{lemma:sphere}}:
We cover the sphere $S$ with bounded covering
number by closed spherical caps $X_j$ of Euclidean radius $1$, and denote by
$\bar B_j$ the closed $m$-dimensional Euclidean balls of radius $1$ having the same centers as $X_j$.
The total number of the caps in the cover is $\lesssim R^{m-1}$.
By $Y_i$ we denote the connected components of $S\setminus Z(f)$. In each domain $Y_i$ we fix a point $z_i$ where the gradient
$\nabla f (z_i)$ is directed normally to $S$, that is, $\nabla_S f(z_i)=0$. The number of $i$'s such that, for some $j$,
$X_j \subset Y_i$ is $\lesssim R^{m-1}$. Thus, in what follows, we consider only those $i$'s for which $Y_i$ does not contain
any $X_j$.

In order to apply Lemma~\ref{lemma:abstract} with the function $\Phi_S$ and with
$W_j = \| f  \|_{\bar B_j, 1+\beta}$
we need to establish the lower bounds for the integrals
\[
\int_{X_j\cap Y_i} \Phi_S \, \rd\vol_S (x) \qquad {\rm with\ } \Phi_S = |f|^{-t}\, |\nabla_s f|^{-t(m-1)}\,,
\]
assuming that $z_i\in X_j$.

Since the sets $Y_i \cap X_j$ and $ X_j\setminus Y_i $ aren't empty, the closed set $\partial Y_i \cap X_j $ is not empty too.
Denote $\rho_i = {\rm dist}(z_i, \partial Y_i \cap X_j)\le 2$ and take a closest to $z_i$ point $p\in \partial Y_i \cap X_j$.
By $V_i$ we denote the spherical cap centered at $z_i$ such that $p\in\partial V_i$.
Note that, by the construction, $\vol_S (Y_i \cap X_j) \ge \vol_S (V_i \cap X_j) \gtrsim \rho_i^{m-1}$.
Since $f(p)=0$,  we have
\[
|f(x)| \lesssim \rho_i\, \| f  \|_{\bar B_j, 1+\beta} = \rho_i W_j\,, \qquad x\in V_i \cap X_j\,.
\]
Furthermore, since $\nabla_S f(z_i)=0$,
\[
|\nabla_S f (x)| \lesssim \rho_i^\beta\,   \| f  \|_{\bar B_j, 1+\beta}
+ \frac{\rho_i}R  \| f  \|_{\bar B_j, 1+\beta}
\lesssim  \rho_i^\beta W_j\,,  \qquad x\in V_i \cap X_j\,.
\]
Hence, on $V_i\cap X_j$ we have
\[
\Phi_S  \gtrsim \rho_i^{-t(1+\beta(m-1))}\, W_j^{-tm}
\gtrsim ( \vol_S (V_i \cap X_j) )^{-t(\frac1{m-1}+\beta)}\, W_j^{-tm}\,,
\]
and
\[
\int_{Y_i \cap X_j} \Phi_S \, \rd\vol_S (x) \ge
\int_{V_i \cap X_j} \Phi_S \, \rd\vol_S (x)
\gtrsim ( \vol_S ( V_i \cap X_j ) )^{1-t(\frac1{m-1}+\beta)} W_j^{-tm}\,.
\]
Now, fixing the parameters $t$ and $\beta$ so close to $1$ that
$ t(\frac1{m-1}+\beta) >1$, we see that the RHS is
$ \ge ( \vol_S ( Y_i \cap X_j ) )^{1-t(\frac1{m-1}+\beta)} W_j^{-tm} $.
At last, applying Lemma~\ref{lemma:abstract}, we complete the proof. \hfill $\Box$

\subsection{Regular components}\label{subsect:regular-components}

\begin{definition}\label{def:small-components}
{\rm
We call a connected component $G$ of the set $U\setminus Z(f)$ {\em regular}, if $G$ is compactly supported in $U$ and
$\vol (G) < \vol (B(1))$.
By $N_{\tt reg}(B_R; f)$ we denote the number of regular connected components $G$ compactly contained in $B_R$.
}
\end{definition}

\begin{lemma}\label{lemma:small-components}
There exist constants $q>1$ and $C<\infty$ such that
\[
\cE^* \bigl\{ N_{\tt reg }^q (B_R; f) \bigr\} \le C  R^{qm}\,.
\]
\end{lemma}

\noindent{\em Proof:} The proof of this lemma follows closely that of Lemma~\ref{lemma:sphere}.
Cover the ball $\bar B_R$ by closed balls $X_j$ of radius $1$ with centers in $B_R$ keeping the covering number bounded,
and put $ X = \bigcup_j X_j$. Then $\bar B_R \subset X \subset B_{R+1}$.
Denote by $\{Y_i\}$ the set of regular nodal domains of $f$ that are contained in $B_R$. In each domain $Y_i$ choose a point
$z_i$ with $\nabla f(z_i)=0$. In order to apply Lemma~\ref{lemma:abstract} with the function $\Phi = |f|^{-t} |\nabla f|^{-tm} $ and
with $W_j =  \| f  \|_{\bar B_j, 1+\beta} $, we need to estimate from below the integrals
$\displaystyle \int_{X_j \cap Y_i} \Phi\, \rd\vol $
assuming that $z_i\in X_j$.

Since $\vol (Y_i) < \vol (X_j) $, we note again that
$\partial Y_i \cap X_j \ne \emptyset$.
Put $\rho_i = {\rm dist}(z_i, \partial Y_i \cap X_j) \le 2$, take a closest to $z_i$ point
$p\subset \partial Y_i \cap X_j$, and denote $V_i = B(z_i, \rho_i)$. By the construction,
\[
\vol (Y_i \cap X_j) \ge \vol (V_i \cap X_j) \gtrsim \rho_i^m\,.
\]
Since $f(p)=0$ and $\nabla f(z_i)=0$, we have
\[
|f(x)|  \le \rho_i\,  \| f  \|_{\bar B_j, 1+\beta} = \rho_i W_j, \qquad  x\in V_i\cap X_j\,, \\
\]
and
\[
| \nabla f(x) | \lesssim \rho_i^\beta\,  \| f  \|_{\bar B_j, 1+\beta} = \rho_i^\beta W_j, \qquad x\in V_i\cap X_j\,.
\]
Hence, on $V_i \cap X_j$,
\[
\Phi \gtrsim \rho_i^{-t(1+\beta m)} W_j^{-t(1+m)} \gtrsim
(\vol (V_i \cap X_j) )^{-t(\beta+\frac1m)}\, W_j^{-t(m+1)}\,,
\]
and
\[
\int_{Y_i \cap X_j} \Phi\, \rd\vol \ge \int_{V_i\cap X_j} \Phi\, \rd\vol
\gtrsim (\vol (V_i \cap X_j) )^{1-t(\beta+\frac1m)}\, W_j^{-t(m+1)}\,.
\]
Fixing the parameters $t$ and $\beta$ so close to $1$ that $t( \beta + \frac1m ) > 1$, we get
\[
\int_{Y_j \cap X_j} \Phi\, \rd\vol \gtrsim ( \vol (Y_i \cap X_j ) )^{1-t(\beta+\frac1m)}\, W_j^{-t(m+1)}\,.
\]
Finally, Lemma~\ref{lemma:abstract} ends the job. \hfill $\Box$

\subsection{The moment estimate for the total number of connected components}\label{subsubsect:moment_estimate}

If the function $f$ is $C^1$-smooth and $0$ is not a critical value, then we can bound the number of connected
components $\gamma$ of $Z(f)$ contained in $B_R$ by the number of connected components $G$ of $U\setminus Z(f)$ compactly
contained in $B_R$. All we need for that is to note that each $\gamma\subset B_R$ is the outer boundary\footnote{
i.e., the part of the boundary of $G$ that bounds the unbounded connected component of $\bR^m\setminus G$ as well
}
of some $G$ compactly supported in $B_R$ and no two different connected components $\gamma\subset B_R$ of $Z(f)$ can serve as
the outer boundary of the same connected component $G$ of $U\setminus Z(f)$ simultaneously.

Thus, combining the estimate of Lemma~\ref{lemma:small-components} with the trivial bound
\[
\# \bigl\{ G\colon \bar G \subset B_R, \ \vol(G)\ge \vol(B(1)) \bigr\} \le R^m\,,
\]
we conclude that, for some $q>1$,
\[
\cE^* \bigl\{ N(B_R; f)^q \bigr\} \lesssim R^{mq}\,.
\]
If, in addition, $f$ is non-degenerate on $S=\partial B_R$ in the sense that $f$ and $\nabla_S f$ do not
vanish simultaneously anywhere on $S$ (due to Lemma~\ref{lemma:Bulinskaya} this event has probability $1$),
then, arguing in a similar way,  we can estimate the number of connected components of $Z(f)$ intersecting
$S$ by the number of  connected components of $S\setminus Z(f)$. Thus, the result of Lemma~\ref{lemma:sphere}
can be viewed as an upper bound for the $q$-th moment of the number of connected components of $Z(f)$ intersecting
$S$.

\medskip
We will use these observations several times when referring to Lemmas~\ref{lemma:small-components} and~\ref{lemma:sphere}
as if they were about the connected components of $Z(f)$ rather than about those of $U\setminus Z(f)$ and $S\setminus Z(f)$.

\section{Proof of Theorem~\ref{thm:Euclid-version}}\label{sect:proof-Thm1}

\subsection{Preliminaries}\label{subsect:proof-thm1-preliminaries}
We need several basic notions from the ergodic theory. Suppose  $ \bigl( \Omega, \mathfrak
S, \cP \bigr) $ is a probability space on which $ \bR^m$ acts by measure-preserving
transformations $\tau_v $, $v\in\bR^m$. This means that for each $v\in\bR^m$,
$\tau_v\colon \Omega\to\Omega$ is a $\mathfrak S$-measurable transformation, $\tau_u \circ \tau_v =
\tau_{u+v}$, $\tau_{-v}=\tau_v^{-1}$, and for each $v\in\bR^m$ and each $A\in\mathfrak S$,
we have $\cP (\tau_v A) = \cP (A)$.

The following version of Wiener's ergodic theorem suffices for our purposes:

\medskip\par\noindent{\bf Wiener's ergodic theorem:}
{\em\ Suppose $\bigl( \Omega, \mathfrak S, \cP \bigr) $
is a probability space on which $ \bR^m$ acts by
measure-preserving transformations $\tau_v$, $v\in\bR^m$. Suppose that $\Phi\in
L^1(\cP)$, and that the function $(v, \omega)\mapsto \Phi\circ\tau_v$
is measurable with respect to the product $\sigma$-algebra $\mathfrak B(\bR^m) \times\mathfrak S$,
where $\mathfrak B(\bR^m)$ is the Borel $\sigma$-algebra generated by open sets in $\bR^m$.
Suppose that $S\subset\bR^m$ is a
bounded open convex set containing the origin.
Then the limit
\[
\lim_{R\to\infty} \frac1{\vol S(R)}\, \int_{S(R)} \Phi (\tau_v \omega) \, \rd\vol(v)
= \bar\Phi (\omega)
\]
exists with probability $1$ and in $L^1(\cP)$. The limiting random
variable $\bar\Phi$ is $\tau$-invariant {\rm (}i.e., for each $v\in\bR^m$,
$\bar\Phi\circ\tau_v = \bar\Phi${\rm )}, and does not depend on the choice of the convex set $S$.}

\medskip This is a special case of a theorem proven in Becker~\cite[Theorems~2 and 3]{Becker}.
Note that Becker's formulation of this theorem deals with rather general increasing families $( U_R )$ of open
sets in $\bR^m$ satisfying two conditions:

\smallskip\noindent (A) {\em the Hardy-Littlewood maximal operator associated with
the family $ ( U_R )$ is of weak type $(1,1)$},

\smallskip\noindent and

\smallskip\noindent (B) {\em for each $t\in\bR^m$},
\[
\lim_{R\to\infty} \vol ( (t+U_R) \triangle U_R) / \vol (U_R) = 0\,,
\]
where $\triangle$ denotes the symmetric difference.

\smallskip
In the case when $S$ is the unit ball, condition (A) reduces to
the classical Hardy-Littlewood maximal theorem, after which it remains to note that the maximal function
associated with the family $S(R)$ is dominated (up to a constant factor) by the one corresponding to the unit ball.
The verification of condition (B) is straightforward.

Note that Becker's presentation does not formally contain the claim that the limiting random variable $\Phi$ does not depend on the family
$ (U_R) $ but in our situation
it can be easily established by applying Becker's theorem to a family $ U_R $ containing
arbitrarily large homothetic images of two bounded convex sets  $S'$ and $S''$.

\medskip Next, recall that the action of $\bR^m$ is called
{\em ergodic}\footnote{a.k.a. metric-transitive} if for
every set $A\in\mathfrak S$ satisfying $\cP( (\tau_v A) \triangle A) = 0$,
either $\cP (A)=0$, or $\cP (A)=1$. In the ergodic case, the limiting random variable
$\bar\Phi$ is a constant function. Due to the $L^1 (\cP)$-convergence,
the value of this constant equals the expectation of
$\Phi$: $\bar\Phi = \cE \bigl\{ \Phi \bigr\}$.

\medskip
Let $X\subset C(\bR^m)$ be an invariant set of continuous functions
(i.e., $G\in X$ implies $G\circ\tau_v\in X$ for all $v\in\bR^m$). Let $\mathfrak S$
be the minimal $\sigma$-algebra on $X$ containing all ``intervals''
$I(u; a, b)=\bigl\{ G\in X\colon G(u)\in [a, b) \bigr\}$. Let
$\gamma$ be a Gaussian probability measure on $(X, \mathfrak S)$
meaning that for every finitely many points $u_1, \ldots , u_k~\in~\bR^m$, the
push-forward of $\gamma$ by the mapping $G\mapsto [ G(u_1), \ldots , G(u_k) ]$ is a (centered)
Gaussian, possibly degenerate, measure on $\bR^k$. If $\gamma$ is invariant under the introduced
action of $\bR^m$ on $X$,
then \[ \bR^m \times X\ni(u, G) \mapsto G(u)\in\bR \] is a translation-invariant
Gaussian function on the probability space $(X, \mathfrak S, \gamma)$ with continuous
trajectories and continuous covariance kernel and we can talk about its spectral
measure $\rho$.

\medskip\noindent{\bf Fomin-Grenander-Maruyama theorem:}
{\em Suppose that $\rho$ has no atoms. Then the action of
$\bR^m$ on $(X, \mathfrak S, \gamma)$ by translations is ergodic.
}

\medskip\noindent For the reader's convenience, we remind
the proof of this theorem\footnote{
The full version of the Fomin-Grenander-Maruyama theorem states that the continuity of the
spectral measure $\rho$ is {\em necessary and sufficient} for the ergodicity of the
action of
$\bR^m$ on $(X, \mathfrak S, \gamma)$ by translations. We will use (and prove) only the sufficiency
part. The proof we present follows
the argument
for the univariate case given in~\cite[Section~5.10]{Grenander}.
}
in Appendix\!~\ref{App-B}.

\medskip Now, let $F$ be a Gaussian function on $\bR^m$ satisfying the assumptions
of Theorem~\ref{thm:Euclid-version}. By the moment assumption $(\rho1)$, with probability
$1$ it is $C^{2-}$-smooth. Hence, it generates a Gaussian measure $\gamma_F$ on
$(C^1(\bR^m), \mathfrak B(C^1(\bR^m)))$ where $\mathfrak B(C^1(\bR^m))$ is the Borel $\sigma$-algebra
generated by open sets in $C^1(\bR^m)$. In what follows, it will be convenient to pass
from $C^1(\bR^m)$ to its subset
\[
C^1_*(\bR^m) = \bigl\{ G\in C^1(\bR^m)\colon |G| + |\nabla G| \ne 0 \bigr\},
\]
which consists of functions for which $0$ is not a critical value.
Note that $C^1_*(\bR^m)$ is a Borel subset of $C^1(\bR^m)$ and, by
the first statement in Lemma~\ref{lemma:Bulinskaya}, \[ \gamma_F\bigl( C^1(\bR^m) \setminus C^1_*(\bR^m) \bigr) = 0.\]
Furthermore, $\bR^m$ acts on $\bigl( C^1_* (\bR^m), \mathfrak B(C_*^1(\bR^m)), \gamma_F \bigr)$
by translations and, since the distribution of $F$ is translation invariant,
the action is measure-preserving.
Thus, Wiener's theorem applies in this setting. To apply the Fomin-Grenander-Maruyama
theorem, we only need to note that the Borel $\sigma$-algebra $\mathfrak B(C_*^1(\bR^m))$
coincides with the $\sigma$-algebra  $\mathfrak S$
generated by the intervals $I(u; a, b)$ (see~Appendix~\ref{A_subsect_C^k}).

\medskip
We conclude that
\begin{itemize}\item
{\em under the assumption $(\rho1)$ of Theorem~\ref{thm:Euclid-version}, for any
random variable $\Phi\in L^1(\gamma_F)$ such that the function
$(v, G)\mapsto \Phi(\tau_v G)$ is
measurable, the ergodic averages
\[
( A_R^S \Phi )(G) \stackrel{\rm def}= \frac1{\vol S(R)}\, \int_{S(R)} \Phi(\tau_v G)\,\rd\vol (v)
\]
converge to a $\tau$-invariant limit $\bar\Phi$ with probability $1$, as well as in
$L^1(\gamma_F)$, as $R\to\infty$. Moreover, under assumption $(\rho2)$, we have}
$ \bar\Phi = \cE \bigl\{ \Phi \bigr\} $.
\end{itemize}

\medskip
We split the proof of Theorem~\ref{thm:Euclid-version} into two
parts: first, we prove the convergence of $ (\vol S(R))^{-1}N_S(R; F) $ to
a limit $\nu$. Then, assuming condition ($\rho$4), we show that
this limit is positive.

\subsection{Existence of the limit}\label{subsect:proof-thm1-existence-limit}

\subsubsection{The sandwich estimate for $ N_S(R; G)/\vol S(R) $}\label{subsubsect:proof-thm1-sandwich}
Without loss of generality, we assume that $S \supset B(1)$. Then, the integral-geometric Lemma~\ref{lemma:integr-geom} provides us with the ``sandwich estimate'':
\begin{multline*}
\frac1{\vol S(R)}\int_{S(R-r)} \frac{N(v, r; G)}{\vol B(r)}\, {\rm d}\vol (v) \le
\frac{N_S(R; G)}{\vol S(R)} \\
\le \frac1{\vol S(R)}\, \int_{S(R+r)} \frac{N^*(v, r; G)}{\vol B(r)}\,
\rd\vol (v)\,.
\end{multline*}
The difference $ N^*(v, r; G) - N(v, r; G) = N^*(r; \tau_v G) - N(r; \tau_v G)$ is bounded
by $ \mathfrak N_{\#} (r; \tau_v G) $, where
\[
\mathfrak N_{\#} (r; G) \stackrel{\rm def}=
\begin{cases}
\mathfrak N(\partial B(r); G) & \text{ if } G \text{ is non-degenerate on } \partial B(r), \\
+ \infty & \text{ otherwise},
\end{cases}
\]
and $ \mathfrak N(\partial B(r); G) $ is the number of connected components of $\partial B(r)\setminus Z(G)$.
Recall that we say that $G$ is non-degenerate on the sphere $\partial B(r)$ if $G$ and $\nabla_{\partial B(r)} G$
do not vanish simultaneously anywhere on $\partial B(r)$.

We introduce the functionals
\[
\Phi_r(G) \stackrel{\rm def}= \frac{N(r; G)}{\vol B(r)}, \qquad
\Psi_r(G) \stackrel{\rm def}= \frac{\mathfrak N_\#(r; G)}{\vol B(r)}\,.
\]
Then the sandwich estimate takes the form
\begin{multline}\label{eq:sandwich}
\Bigl( 1-\frac{r}R \Bigr)^m (A_{R-r}^S \Phi_r)(G)
\le \frac{N_S(R; G)}{\vol S(R)}
\\
\le \Bigl( 1+\frac{r}R \Bigr)^m \bigl[ (A_{R+r}^S \Phi_r)(G) + ( A_{R+r}^S \Psi_r )(G) \bigr].
\end{multline}

\subsubsection{Checking measurability}\label{subsubsect:proof-thm1-measurability}

We need to check that, given $r>0$, the functions
\[
(v, G) \mapsto \Phi_r (\tau_v G), \quad (v, G) \mapsto \Psi_r (\tau_v G)
\]
are measurable with respect to the product $\sigma$-algebra
$\mathfrak B(\bR^m) \times \mathfrak B(C^1_*(\bR^m)) $.
The function $(v, G)\mapsto \tau_v G$ is
a measurable (even continuous) map
\[
\bigl( \bR^m \times C^1_*(\bR^m), \mathfrak B(\bR^m) \times \mathfrak B(C^1_*(\bR^m) \bigr)
\to \bigl( C^1_*(\bR^m), \mathfrak B(C^1_*(\bR^m)) \bigr)\,.
\]
Since the composition of measurable functions is measurable, it remains to show that, given $r>0$,
the functions $G\mapsto N(r, G)$ and $G\mapsto \mathfrak N_\#(r, G)$ are measurable as maps
from $\bigl( C^1_*(\bR^m),  \mathfrak B(C^1_*(\bR^m)) \bigr)$ to
$\bigl( [0, +\infty], \mathfrak B( [0, +\infty] ) \bigr)$.

The measurability of the map $G\mapsto N(r, G)$ follows from its
lower semicontinuity  on $C^1_*(\bR^m)$. To see that $G\mapsto \mathfrak N_\#(r, G)$
is measurable, first, consider the set ${\rm Degen}(r)$ of functions $G\in C^1_* (\bR^m)$
for which there exists a point $x\in \partial B(r)$ such that $\nabla G(x)$ is orthogonal
to the tangent space to $\partial B(r)$ at $x$.
This set is closed in $C^1_* (\bR^m)$ with respect to the $C^1$-topology
and, therefore, is $\mathfrak B(C^1_*(\bR^m))$-measurable. On the other hand, our map $G\mapsto \mathfrak N_\#(r, G)$
is lower semi-continuous on $C^1_*(\bR^m)\setminus {\rm Degen}(r)$.

\subsubsection{Integrability}\label{subsubsect:proof-thm1-integrability}
Next, we note that, for every fixed  $r>0$, the functions
$\Phi_r$ and $\Psi_r$ on $C^1_*(\bR^m)$ are $\gamma_F$-integrable.
This readily follows from Lemma~\ref{lemma:small-components} and Lemma~\ref{lemma:sphere}
correspondingly.

\subsubsection{Proof of convergence}\label{subsubsect:proof-thm1-convergence}
By the sandwich estimate~\eqref{eq:sandwich}, for every function $G\in C^1_*(\bR^m)$, we have
\begin{multline}\label{eq:mult*}
\Bigl| \frac{N_S(R; G)}{\vol S(R)} - (A_R^S\Phi_r)(G)\Bigr| \le
\Bigl| \Bigl( 1 - \frac{r}R \Bigr)^m (A_{R-r}^S\Phi_r)(G)  - (A_R^S \Phi_r)(G) \Bigr| \\ +
\Bigl| \Bigl( 1 + \frac{r}R \Bigr)^m (A_{R+r}^S\Phi_r)(G)  - (A_R^S \Phi_r)(G) \Bigr|
+ \Bigl( 1 + \frac{r}R \Bigr)^m (A_{R+r}^S\Psi_r)(G)\,.
\end{multline}

By the Wiener ergodic theorem, there exist $\tau$-invariant functions $\bar\Phi_r$ and $\bar\Psi_r$
such that
\[
\lim_{R\to\infty} A_R^S \Phi_r = \bar\Phi_r \quad {\rm and} \quad \lim_{R\to\infty} A_R^S\Psi_r = \bar\Psi_r
\]
both $\gamma_F$-almost everywhere and in $L^1(\gamma_F)$. Letting $R\to\infty$ on both sides of~\eqref{eq:mult*}, we get
\begin{equation}\label{eq:A}
\varlimsup_{R\to\infty} \Bigl| \frac{N_S(R; G)}{\vol S(R)} - (A_R^S\Phi_r)(G)\Bigr| \le
\bar\Psi_r (G) \qquad \text{for\ } \gamma_F\text{-almost every\ } G,
\end{equation}
and
\begin{equation}\label{eq:B}
\varlimsup_{R\to\infty} \int \Bigl| \frac{N_S(R; G)}{\vol S(R)} - (A_R^S\Phi_r)(G)\Bigr|\, \rd\gamma_F (G)
\le \int \bar\Psi_r\, \rd\gamma_F =  \frac{\cE \{ \mathfrak N_\# (r; F)\}}{\vol B(r)}\,.
\end{equation}
By Lemma~\ref{lemma:sphere}, the RHS of~\eqref{eq:B} is $ \lesssim r^{-1}$ for $ r\ge 1$.
So taking a sequence $r_k\uparrow\infty$, we observe that
\[
\lim_{k\to \infty} \int \bar\Psi_{r_k}\, \rd\gamma_F = 0\,.
\]
and, consequently,
\[
\inf_k \bar\Psi_{r_k} = 0\, \qquad \gamma_F\text{-almost everywhere}\,.
\]
Since $A_R^S \Phi_r(G)$ converge to $\bar\Phi_r$ for $\gamma_F$-almost every $G$,
the second observation together with~\eqref{eq:A} imply that  $  (\vol S(R))^{-1}\, N_S(R; G) $ is Cauchy
for $\gamma_F$-almost every $G$. Similarly, the convergence of
$A_R^S \Phi_r(G)$ to $\bar\Phi_r$  in $L^1(\gamma_F)$ together with the first
observation and~\eqref{eq:B} imply that  $  (\vol S(R))^{-1}\, N_S(R; G) $ is Cauchy in $L^1(\gamma_F)$.
Thus, the limit
\[
\nu \stackrel{\rm def}= \lim_{R\to\infty} \frac{N_S(R; G)}{\vol S(R)}
\]
exists $\gamma_F$-almost everywhere and in $L^1(\gamma_F)$. It follows from~\eqref{eq:sandwich} that, for every $r>0$,
\[
\bar\Phi_r \le \nu \le \bar\Phi_r + \bar\Psi_r \qquad \gamma_f-{\rm almost\ everywhere}\,.
\]

If, in addition, the action of $\bR^m$ on $\bigl( C^1_*(\bR^m), \mathfrak B(C^1_*(\bR^m)), \gamma_F \bigr)$ is ergodic,
then
$\bar\Phi_r = \cE\{ \Phi_r \}$, $\bar\Psi_r = \cE\{ \Psi_r \}$. Therefore,
\begin{equation}\label{eq:***}
\cE\bigl\{ \bar\Phi_r \bigr\} \le \nu \le \cE\bigl\{ \bar\Phi_r \bigr\} + \cE\bigl\{ \bar\Psi_r \bigr\}
\qquad \gamma_f-{\rm almost\ everywhere}\,,
\end{equation}
whence, for every $r>0$,
$\gamma_F$-essential oscillation of $\nu$ does not exceed $\cE\{ \Psi_r \}$.
Recalling that $\cE\{ \Psi_r \} \lesssim r^{-1}$ and letting $r\to\infty$,
we see that
$\nu$ is a (non-random) constant. This completes the proof of convergence in Theorem~\ref{thm:Euclid-version}.
\hfill $\Box$

\subsection{Positivity of $\nu$}\label{subsect:proof-thm1-positivity}
It remains to show that condition ($\rho4$)
yields the positivity of the limiting constant $\nu$. We prove that if assumption ($\rho$4)
holds, then $\cP\{N(r; F)>0\}>0$ when $r$ is sufficiently big. Since
the LHS of estimate~\eqref{eq:***} can be rewritten as
$\nu \ge \cE \{ N(r; F)\}/\vol B(r)$ for each $r>0$,
this will yield the positivity of $\nu$.

\subsubsection{A Gaussian lemma}\label{subsubsect:proof-thm1-Gaussian-lemma}

\begin{lemma}\label{lemma:gauss-funct-anal}
Let $\mu$ be a compactly supported Hermitian measure
with $\spt(\mu)\subset \spt(\rho)$. Then for each closed ball $\bar B\subset\bR^m$ and for each $\e>0$,
\[
\cP
\bigl\{ \| F - \widehat\mu \|_{C(\bar B)} < \e \bigr\} >0\,.
\]
\end{lemma}

\noindent{\em Proof of Lemma~\ref{lemma:gauss-funct-anal}}:
The part of the theory of continuous Gaussian functions developed in Appendix (\ref{A_subsect_top_support} and
~\ref{A_subsect_tr-inv}) yields
the statement of the lemma for all measures $\mu$ absolutely continuous with respect
to $\rho$ with density
\[
h\in L^2_{\tt H}(\rho) \stackrel{\rm def}= \bigl\{g\in L^2(\rho)\colon
g(-x)=\overline{g(x)}\,  \ {\rm for\ all\ } x\in\bR^m \bigr\}\,.
\]

In the general case, we can approximate the measure $\mu$ in the weak topology by measures
${\rm d}\mu = h\, {\rm d}\rho$ with $\operatorname{spt}(h)$ contained
in a fixed compact neighbourhood of $\operatorname{spt}(\mu)$.
Then it remains to recall that
for measures supported on a fixed compact set, the weak
convergence yields locally uniform convergence of their Fourier integrals.
\mbox{} \hfill $\Box$

\subsubsection{Proof of the positivity of $\nu(\rho)$}\label{subsubsect:proof-thm1-positivity}
We take a Hermitian compactly supported measure $\mu$ with $ \spt(\mu)
\subset \spt(\rho) $ and a bounded domain $D\subset\bR^m$ so that $\widehat{\mu}\big|_{\partial
D}<0$ and $\widehat{\mu}(u_0)>0$ for some $u_0\in D$. Choose $r$ so big that $\bar D \subset
B(r)$.  If $\e>0$ is sufficiently small, then $G(u_0)>0$ and $G\big|_{\partial D}<0$ for every function
$G$ satisfying $\| G-\widehat\mu \|_{C(\bar B(r))}<\e$.
Thus, for every such function $G$, the zero set $Z(G)$ has at least one connected component in $D$.
Applying Lemma~\ref{lemma:gauss-funct-anal}, we see that
\[
\cP \bigl\{ N(r; F)>0\bigr\}
\ge \cP \bigl\{ \| F-\widehat{\mu} \|_{C(\bar B(r))} < \e \bigr\} >0
\]
completing the proof of Theorem~\ref{thm:Euclid-version}. \mbox{}\hfill $\Box$

\section{Recovering the function $ \bar\nu $ by a double scaling limit}\label{sect:double-scaling-limit}

The proof of Theorem~\ref{thm:main_thm} will rely upon the following lemma, which is of independent
interest.
Let $(f_L)$ be a tame parametric Gaussian ensemble, that is, an ensemble satisfying the assumptions
of Theorem~\ref{thm:main_thm}.
As above, we put
$f_{x, L}(u)=f(x+L^{-1}u)$ and \[ K_{x, L}(u, v) = \cE \bigl\{ f_{x, L}(u)
f_{x, L}(v) \bigr\} = K_L(x+L^{-1}u, x+L^{-1}v)\,.\]
Till the end of this section, we
fix a point $x\in U$ so that
\[
\lim_{L\to\infty}\, K_{x, L}(u, v) = k_x(u-v) \qquad {\rm pointwise\ in\ } \bR^m\times \bR^m\,,
\]
where the Hermitian positive-definite  function $k_x$ is the
Fourier integral of a measure $\rho_x$ satisfying assumptions ($\rho1$)--($\rho3$).
By $F_x$ we denote the limiting Gaussian function on $\bR^m$, and put $\nu=\bar\nu (x) = \nu (F_x)$.

\begin{lemma}\label{lemma:local}
For every $\e>0$,
\[
\lim_{R\to\infty} \varlimsup_{L\to \infty}
\cP \Bigl\{ \Bigl| \frac{N(R; f_{x, L})}{\vol B(R)}  - \nu  \Bigr|>\e \Bigr\} = 0\,.
\]
\end{lemma}

\noindent{\em Proof\,}: Fix $R>2$ and $\e >0$. Our goal will be to show that, for every $t$,
\begin{align*}
\varlimsup_{L\to\infty} \cP \Bigl\{ N(R; f_{x, L}) > t \Bigr\}
&\le \cP \Bigl\{ N(R+1; F\ci x) > t \Bigr\}\,, \\
\varlimsup_{L\to\infty} \cP \Bigl\{ N(R; f_{x, L}) < t \Bigr\}
&\le \cP \Bigl\{ N(R-1; F\ci x)< t \Bigr\}\,.
\end{align*}
Applying these inequalities with $t=(\nu+\e)\vol B(R)$ and $t=(\nu-\e)\vol B(R)$ respectively, and
then combining the results with Theorem~\ref{thm:Euclid-version}, we get the
conclusion of Lemma~\ref{lemma:local}.
The proofs of these two relations are very similar, so we will present only the proof of the first one.

\medskip
We choose a big constant $M$ and a small constant $\kappa$ so that the kernels $k_x$ and $K_{x, L}$ (with $L\ge L_0$)
satisfy the ``$(M, \kappa)$-conditions'' introduced in the beginning of Section~\ref{subsect:Bulinsk-prelim}.
For the kernel $k_x$ this is possible due to conditions ($\rho$1) and ($\rho$3). For the scaled kernels $K_{x, L}$
this is possible due to the controllability of $(f\ci L)$.

\medskip
Given positive constants $A$ and $a$, we put
\[
E(A, a) = \Bigl\{
g\in C^1(B(R+1.1))\colon \| g \|_{C^1(\bar B(R+1))} \le A\,, \ \min_{\bar B(R+1)}\, \max\{ |g|, |\nabla g| \} \ge a  \Bigr\}\,.
\]
Introduce the events $ \Omega'_L = \bigl\{ f_{x, L}\notin E(A, a) \bigr\}$ and $\Omega'' = \bigl\{ F_x\notin E(A, a)\bigr\} $.
By Lemma~\ref{lemma:prob-stable}, the aforementioned ``$(M, \kappa)$-conditions'' imply that, for a given $\delta>0$,
we can make the probabilities of both events less than $\delta$ if we choose sufficiently big $A$ and
sufficiently small $a$.
We fix a finite $a/(2A)$-net $X$ in $\bar B(R+1)$ and denote by $E'\subset \bR^{|X|}$ the set of traces on $X$ of functions
$g\in E(A, a)$ satisfying $N(R; g) > t$. This is a bounded subset of $\bR^{|X|}$.
Note that if $g, h\in E(A, a)$ and $|g-h|< a/2$ on $X$, then $|g-h|< a$ everywhere on $\bar B(R+1)$, and by Lemma~\ref{lemma:Cor4.3} (applied with $\alpha=\beta=a$), $N(R+1; h) \ge N(R; g)$.

\medskip
We fix a function $\phi\in C^\infty_0(\bR^{|X|})$ satisfying $0\le \phi \le 1$
everywhere, $\phi\equiv 1$ on $E'$
and $\phi\equiv 0$ on $\bR^{|X|}\setminus E'_{+a/2}$
(as usual, by $E'_{+s}$ we denote the $s$-neighbourhood of $E'$), and
consider the finite dimensional Gaussian vectors $f_{x, L}|_X$ and $F_x|_X$.
First, we note that
\begin{multline*}
\bigl\{\omega\colon N(R; f_{x, L}) > t  \bigr\} \subset
\bigl\{\omega\colon f_{x, L}|_X \in E', \ f_{x, L}\in E(A,a) \bigr\}\, \cup\, \Omega_L' \\
\subset \bigl\{\omega\colon \phi (f_{x, L}|_X)=1 \}\, \cup\, \Omega_L'\,,
\end{multline*}
whence,
\[
\cP \bigl\{N(R; f_{x, L}) > t  \bigr\} < \cE \{ \phi (f_{x, L}|_X) \} + \delta\,.
\]

The pointwise convergence of the scaled kernels $K_{x, L}(u, v)$
to the limiting kernel $k_x(u-v)$ yields\footnote
{
If $\xi\ci L$ are Gaussian $n$-dimensional vectors and the entries of the covariance matrices $K\ci L$
of $\xi_L$ converge to the entries
of the covariance matrix $K$ of $\xi$, then
\begin{multline*}
\cE \bigl\{ \phi (\xi_L) \bigr\} = \cE \Bigl\{ \int_{\bR^n} \widehat{\phi}(\la) e^{2\pi{\rm i}\la\cdot \xi_L}\, \rd\la \Bigr\}
= \int_{\bR^n} \widehat{\phi}(\la) \cE \Bigl\{ e^{2\pi{\rm i}\la\cdot \xi_L} \Bigr\}\, \rd\la \\
= \int_{\bR^{n}} \widehat{\phi}(\la) e^{-\pi K_L\la\cdot\la}\, \rd\la \to
\int_{\bR^n} \widehat{\phi}(\la) e^{-\pi K\la\cdot\la}\, \rd\la =  \cE \bigl\{ \phi (\xi) \bigr\},
\end{multline*}
where the convergence holds by the dominated convergence theorem.
}
\[
\cE \{ \phi (f_{x, L}|_X) \} \stackrel{L\to\infty}\to \cE \{ \phi (F\ci x|_X) \} \le \cP \bigl\{\phi (F_x|_X)>0\bigr\}
\]
(in the inequality we used that $\phi \le 1$ everywhere).
Now,
\begin{multline*}
\{\omega\colon \phi (F\ci x|_X)>0 \} \subset  \bigl\{\omega\colon F\ci x|_X \in E'_{+a/2}\bigr\}
\subset \bigl\{\omega\colon F\ci x|_X \in E'_{+a/2}, \ F_x\in E(A, a) \bigr\} \cup \Omega'' \\
\subset
\bigl\{\omega\colon N(R+1; F_x) > t  \bigr\} \cup \Omega''\,.
\end{multline*}
In the last step we used
that, by our construction, if $F_x\in E(A, a)$ and $ F\ci x|_X \in E'_{+a/2} $,
then there is a function $g\in E(A, a)$ such that $N(R, g) > t$, and
$|F_x-g|<\frac12 a$ on $X$, whence,
$ N(R+1; F_x)\ge N(R; g) > t $. Hence,
\[
\cP \{ \phi (F\ci x|_X)>0 \} < \cP\bigl\{ N(R+1; F_x) > t  \bigr\} + \delta\,.
\]
Thus, for sufficiently large $L$, we have
\begin{multline*}
\cP \bigl\{N(R; f_{x, L}) > t  \bigr\} < \cE \{ \phi (f_{x, L}|_X) \} + \delta \\
< \cP \bigl\{\phi (F_x|_X)>0\bigr\} + 2\delta < \cP\bigl\{ N(R+1; F_x) > t  \bigr\} + 3\delta\,,
\end{multline*}
completing the argument. \hfill $\Box$

\section{Proof of Theorem~\ref{thm:main_thm}}\label{Sect:Proof_Thm2}

It remains to tie the ends together. Let $(f_L)$ be a \tame\ \pge\, on an open set $U\subset\bR^m$.
This implies that,
\begin{itemize}\item
{\em for every compact set $Q\subset U$,
there exist constants $M<\infty$ and $\kappa>0$  such that the covariance kernels of
the functions $f_{x, L}$ on ${\bar B(R+1)}$ satisfy the $(M, \kappa)$-conditions
from Section~\ref{subsect:Bulinsk-prelim} whenever $x\in Q$, $R>0$, and $L\ge L_0(Q, R)$}.
\end{itemize}

Fix a Borel set $U'\subset U$ of full volume on which the scaled functions $f_{x, L}$
have translation invariant limits $F_x$. Then, by Appendix~\ref{A_subsect_tr-inv},
\begin{itemize}\item
{\em the covariance kernels $k_x(u-v)$ of the limiting functions $F_x$ satisfy the $(M, \kappa)$-conditions
whenever $x\in Q \cap U'$}.
\end{itemize}

\subsection{$\bar\nu\in L^\infty_{\rm loc} (U)$}\label{subsect:proof-thm2-L-infty}

First, we show that $\bar\nu$ is locally uniformly bounded on $U'$ and then that it is measurable.

\subsubsection{Boundedness of $\bar\nu$ }
Recall that
\[
\bar\nu (x) =
\lim_{R\to\infty} \frac{\cE \{ N(R; F_x) \}}{\vol B(R)}\,, \qquad x\in U'\,.
\]
Given any compact set $Q\subset U$,
Lemma~\ref{lemma:small-components} implies that, for every $x\in U' \cap Q $, we have
$ \cE \{ N(R; F_x) \} \le C(Q)\vol B(R) $. Thus,
the function $\bar\nu$ is locally bounded on $U'$. \hfill $\Box$

\subsubsection{Measurability of $\bar\nu$}

Put
\[
\nu_{R, L} (x, \omega) = \frac{N (R; f_{x, L})}{\vol B(R)}\,.
\]
The function $\nu_{R, L}$ is defined on the set $U_{-(R+1)/L}\times \Omega'$, where
$U_{-r}=\bigl\{x\in U\colon {\rm dist}(x, \partial U)>r\bigr\}$
and $\Omega'=\bigl\{ \omega\in\Omega\colon f_L\in C^1_*(U) \bigr\}$, $\cP(\Omega\setminus\Omega')=0$.
It is measurable as a composition of a lower semicontinuous
mapping
\[
C^1_* (B(R+1)) \ni g \mapsto \frac{N (R; g)}{\vol B(R)} \in\bR\,,
\]
a continuous mapping
\[
U_{-(R+1)/L} \times C^1_*(U)   \ni (x, g) \mapsto g_{x, L}\big|_{B(R+1)} \in C^1_* (B(R+1))\,,
\]
and a measurable mapping
\[
U_{-(R+1)/L} \times \Omega' \ni (x, \omega) \mapsto (x, f_L) \in U_{-(R+1)/L} \times C^1_*(U)\,.
\]

Fix $x\in U'$.
By Lemma~\ref{lemma:small-components}, there exist $q>1$ and $C<\infty$ such that
\[
\int_{\Omega'} \nu_{R, L}^q\,\rd\cP < C
\]
for all sufficiently large $L$.
Given $\e>0$, put
\[
\Omega_\e (R, L, x) = \bigl\{\omega\in\Omega'\colon |\nu_{R, L} (x, \omega) - \bar\nu (x)|>\e \bigr\}\,.
\]
Then
\begin{multline*}
\int_{\Omega_\e} |\nu_{R, L}(x, \omega) - \bar\nu (x)|\, \rd\cP(\omega)
\le \int_{\Omega_\e} \nu_{R, L}\, \rd\cP + \bar\nu (x) \cP\{ \Omega_\e \} \\
\le (\cP\{\Omega_\e\})^{1-\frac{1}q}\, \Bigl( \int_{\Omega_\e} \nu_{R, L}^q\, \rd\cP \Bigr)^{\frac1{q}}
+ \bar\nu (x) \cP\{ \Omega_\e \} \le C (\cP\{\Omega_\e\})^{1-\frac{1}q}\,.
\end{multline*}
Therefore,
\[
\Bigl| \int_{\Omega'} \nu_{R, L}(x, \omega)\,\rd\cP (\omega) - \bar\nu (x)  \Bigr| \le
\int_{\Omega'} |\nu_{R, L}(x, \omega) - \bar\nu (x)|\, \rd\cP(\omega) \le \e
+  C (\cP\{\Omega_\e\})^{1-\frac{1}q}
\]
and
\[
\lim_{R\to\infty}\, \varlimsup_{L\to\infty}
\Bigl| \int_{\Omega'} \nu_{R, L}(x, \omega)\,\rd\cP (\omega) - \bar\nu (x)  \Bigr| \le
\e +  C \lim_{R\to\infty}\, \varlimsup_{L\to\infty} (\cP\{\Omega_\e\})^{1-\frac{1}q}\,.
\]
By Lemma~\ref{lemma:local}, the double limit on the RHS vanishes, so
\[
\lim_{R\to\infty}\, \varlimsup_{L\to\infty}
\Bigl| \int_{\Omega'} \nu_{R, L}(x, \omega)\,\rd\cP (\omega) - \bar\nu (x)  \Bigr| = 0\,.
\]
It follows from here that the function $\bar\nu (x)$ can be represented as, say,
\[
\bar\nu (x) = \lim_{R\to\infty}\, \varlimsup_{L\to\infty}\, \int_{\Omega'} \nu_{R, L}(x, \omega)\, \rd\cP(\omega)\,.
\]
Since the functions $ \nu_{R, L}(x, \omega) $ are non-negative and measurable in $(x, \omega)$, their integrals
with respect to $\omega$ over a fixed set $\Omega'$
are also measurable as functions of $x\in U'$.
Thus, the function $\bar\nu$ is also measurable. \hfill $\Box$

\subsection{Towards the proof of Theorem~\ref{thm:main_thm}: another sandwich estimate}\label{subsect:proof-thm2-another-sandwich}
Without loss of generality we assume
that the continuous compactly supported function $\phi$ in the assumptions of Theorem~\ref{thm:main_thm} is non-negative.
We denote $Q=\spt(\phi)$. Fix $\delta>0$ such that $Q_{+4\delta}\subset U$
and put $Q_1=Q_{+\delta}$, $Q_2=Q_{+2\delta}$. For $x\in Q_1$, let
\[
\phi_- (x) = \min_{\bar B(x, \delta)} \phi, \quad \phi_+ (x) = \max_{\bar B(x, \delta)} \phi\,.
\]
Note that
\[
\phi_-(x) \le \phi(y) \le \phi_+(x)
\]
whenever $x\in Q_1$, $y\in B(x, \delta)$.

Fix the parameters $D, R, L$ so that $1<D<R<\delta L$.
We have
\[
L^{-m} \int_U \phi\, \rd n_L = \int_{Q_1} \Bigl( \int_{B(x, R/L)} \frac{\phi (y)\, \rd n_L(y)}{\vol B(R) } \Bigr) \rd\vol (x)\,,
\]
whence,
\begin{multline}\label{eq:sandwich1}
\int_{Q_1} \phi_- (x)\, \frac{n_L(B(x, R/L))}{\vol B(R)}\, \rd\vol (x) \le L^{-m} \int_U \phi\, \rd n_L \\
\le \int_{Q_1} \phi_+(x)\, \frac{n_L(B(x, R/L))}{\vol B(R)}\, \rd\vol (x)\,.
\end{multline}

Since the total $n_L$-mass of each connected component of $Z(f_L)$ equals $1$, the LHS of~\eqref{eq:sandwich1} cannot be less than
\[
\int_{Q_1} \phi_-(x) \nu_{R, L}(x, \omega)\, \rd\vol (x),
\]
where, as above, $\nu_{R, L}(x, \omega) = ( \vol B(R) )^{-1} N(R; f_{x, L})$.

In order to estimate the RHS of~\eqref{eq:sandwich1}, we cover $Q_2$ by $\simeq \vol (Q_2) \bigl( \frac{L}{D} \bigr)^m $ open balls of
diameter $D/L$. Denote by $\bigl\{ S_j \bigr\} $ the collection of boundary spheres of these balls.
Due to the second statement in Lemma~\ref{lemma:Bulinskaya}, with probability $1$ there is no point $x$
such that, for some $j$, $x\in S_j \cap Z(f_L)$ and $\nabla_{S_j} f_L (x) = 0$. Under this non-degeneracy condition,
the number of connected components of $Z(f_L)$ that intersect the sphere $ S_j $ is bounded by the number $ \mathfrak N(S_j; f\ci L) $
of connected components of $S_j\setminus Z(f_L)$.
Denote by $n^*_L$ the part of the component counting measure $n_L$ supported on the connected components of
$Z(f_L)$ intersecting at least one of the spheres $S_j$.
Since every other component of $Z(f_L)$ intersecting a ball $B(x, R/L)$ centered at $x\in Q_1$
is contained in $B(x, (R+D)/L)$, we see that the RHS of~\eqref{eq:sandwich1} does not exceed
\[
\Bigl( \frac{R+D}R \Bigr)^m\,  \int_{Q_1} \phi_+(x) \nu_{R+D, L}(x, \omega) \, \rd\vol (x)
+ \int_{Q_1} \phi_+(x)\, \frac{n_L^*(B(x, R/L))}{\vol B(R)}\, \rd\vol (x)\,.
\]

By Fubini, the second integral on the RHS is bounded by $( \max_U \phi ) L^{-m} n_L^*(Q_2) $. In turn, $ n_L^*(Q_2) \le \sum_j \mathfrak N(S_j; f\ci L)  $
with probability $1$. Thus, for almost every $\omega$, we have
\begin{multline*}
\int_{Q_1} \phi_-(x) \nu_{R, L}(x, \omega)\, \rd\vol (x) \le L^{-m} \int_U \phi\, \rd n_L \\
\le \bigl( 1 + \tfrac{D}R \bigr)^m\,
\int_{Q_1} \phi_+(x)  \nu_{R+D, L}(x, \omega)\, \rd\vol (x)
+ ( \max_U \phi ) \, L^{-m}  \sum\nolimits_j \mathfrak N(S_j; f\ci L)\,.
\end{multline*}

\subsection{Completing the proof of Theorem~\ref{thm:main_thm}}\label{subsect:proof-thm2-completing-the-proof}
To juxtapose the integrals
\[
L^{-m} \int_U \phi\, \rd n_L \quad {\rm and} \quad \int_U \phi\bar\nu\, \rd\vol\,,
\]
we note that, since pointwise $\phi_+ \le \phi + \omega_\phi (\delta)$, where $\omega_\phi$ is the modulus of
continuity of $\phi$, we have
\begin{multline*}
\int_U \phi\bar\nu\, \rd\vol
\ge \int_U \phi_+ \bar\nu\, \rd\vol - \omega_\phi (\delta) \| \bar\nu \|_{L^\infty (Q_1)} \vol(Q_1)
\ge \bigl( 1 + \tfrac{D}R \bigr)^m\, \int_U \phi_+ \bar\nu\, \rd\vol \\[7pt]
- \bigl[ \bigl( 1 + \tfrac{D}R \bigr)^m - 1 \bigr] (\max_U \phi) \| \bar\nu \|_{L^\infty (Q_1)} \vol(Q_1)
- \omega_\phi (\delta) \| \bar\nu \|_{L^\infty (Q_1)} \vol(Q_1)\,,
\end{multline*}
whence, for almost every $\omega$,
\begin{multline*}
L^{-m} \int_U \phi\, \rd n_L - \int_U \phi\bar\nu\, \rd\vol
\le 2^m\, (\max_U \phi )\, \int_{Q_1} | \nu_{R+D, L}(x, \omega) - \bar\nu (x)|\,
\rd\vol (x) \\[7pt]
+ \bigl[ \bigl( 1 + \tfrac{D}R \bigr)^m - 1 \bigr] (\max_U \phi) \| \bar\nu \|_{L^\infty (Q_1)} \vol(Q_1)
+ \omega_\phi (\delta) \| \bar\nu \|_{L^\infty (Q_1)} \vol(Q_1) \\[7pt]
+ ( \max_U \phi ) \, L^{-m} \sum\nolimits_j  \mathfrak N(S_j; f\ci L)\,.
\end{multline*}
The matching lower bound is similar but somewhat simpler: for almost every $\omega$, we have
\begin{multline*}
L^{-m} \int_U \phi\, \rd n_L - \int_U \phi\bar\nu\, \rd\vol \\
\ge - (\max_U \phi )\, \int_{Q_1} | \nu_{R, L}(x, \omega) - \bar\nu (x)|\,
\rd\vol (x) -  \omega_\phi (\delta) \| \bar\nu \|_{L^\infty (Q_1)} \vol(Q_1)\,.
\end{multline*}
Gathering the upper and lower bounds and taking the upper expectation, we obtain
\begin{align*}
\cE^* \Bigl| L^{-m} \int_U \phi\, \rd n_L &-  \int_U \phi \bar\nu\,\rd\vol  \Bigr| \\
& \le 2^m ( \max_U \phi)\, \int_{Q_1}  \cE \{ | \nu_{R+D, L}(x)-\bar\nu (x) | +  | \nu_{R, L}(x)-\bar\nu (x) |\}\, \rd\vol (x) \\[7pt]
& \quad + ( \max_U \phi)\, L^{-m}  \sum\nolimits_j \cE^* \bigl\{ \mathfrak N(S_j; f\ci L) \bigr\} \\[7pt]
&  \quad + ( \max_U \phi)\, \| \bar\nu \|_{L^\infty (Q_1)} \vol (Q_1) \bigl[ \bigl( 1 + \tfrac{D}R \bigr)^m - 1 \bigr] \\[7pt]
& \quad  + 2 \omega_\phi(\delta) \| \bar\nu \|_{L^\infty(Q_1)}\vol (Q_1)\,.
\end{align*}
It remains to estimate the terms on the RHS.

\medskip
Fix $\e>0$ and choose $\delta$ so small that $\omega_f(\delta)<\e$. This takes care of the last term
on the RHS.
To treat the second term we use Lemma~\ref{lemma:sphere}, which yields
$ \cE^* \bigl\{ \mathfrak N(S_j; f\ci L) \bigr\} \lesssim D^{m-1}$ uniformly in $j$. Therefore,
\[
L^{-m}  \sum\nolimits_j \cE^* \bigl\{ \mathfrak N(S_j; f\ci L) \bigr\} \lesssim L^{-m} \cdot \vol (Q_2) (L/D)^m \cdot D^{m-1} \lesssim D^{-1}\vol (Q_2)\,.
\]
Let $U'\subset U$ be a Borel subset of full volume on which the scaled functions $f_{x, L}$ have translation invariant limits.
The functions $\bar\nu$ and $\cE \bigl\{ \nu_{R, L} \bigr\}$ are locally uniformly bounded on $U'$ by a constant independent of $R$ and $L$.
Let
\[ \eta_R (x) \stackrel{\rm def}= \varlimsup_{L\to\infty} \cE \bigl| \nu_{R, L}(x) - \bar\nu (x) \bigr| \,.\]
The function $\eta_R$ is uniformly bounded on $Q_1\cap U'$ by a constant independent of $R$.
Then, we obtain
\begin{multline*}
\varlimsup_{L\to\infty}
\cE^* \Bigl| L^{-m} \int_U \phi\, \rd n_L -  \int_U \phi \bar\nu\,\rd\vol  \Bigr| \\
\quad \le C(\phi, Q) \Bigl( \int_{Q_1} \bigl[ \eta_{R+D}(x) + \eta_R(x) \bigr]\,\rd\vol (x)
+  \| \bar\nu \|_{L^\infty(Q_1)} \bigl[ \bigl( 1 + \tfrac{D}R \bigr)^m - 1 + \e \bigr]
+ D^{-1} \Bigr).
\end{multline*}
For any $x\in U'$, we have
$  \eta_R (x) \to 0 $ as $R\to\infty$. Using the dominated convergence theorem, we get
\[
\varlimsup_{L\to\infty}
\cE^* \Bigl| L^{-m} \int_U \phi\, \rd n_L -  \int_U \phi \bar\nu\,\rd\vol  \Bigr| \\
\le C(\phi, Q) \bigl( \e \| \bar\nu \|_{L^\infty(Q_1)} + D^{-1} \bigr).
\]
Letting $\e\to 0$ and $D\to \infty$.
we finish off the proof of Theorem~\ref{thm:main_thm}. \hfill $\Box$

\section{The manifold case. Proof of Theorem~\ref{thm:manifold}}\label{sect:mflds}\label{sect:proof-thm3}

\subsection{Smooth Gaussian functions and their covariance kernels under $C^2$-changes of variable}\label{subsect:invariance-diff}

Suppose that $U$, $V$ are open  subsets of $\bR^m$ and $\psi\colon V\to U$ is a $C^2$-diffeomorphism.
Suppose that $f\colon U\to \bR$ is a continuous Gaussian function on $U$ with a $C^{2,2}$ covariance kernel
$K(x, y)$. Then $f\circ\psi$ is a continuous Gaussian function on $V$ with the reproducing kernel
$\widetilde{K}(x, y)=K(\psi(x), \psi(y))$. Note that for every pair of the multi-indices $\alpha, \beta$,
the mixed partial derivative $\partial_x^\alpha\, \partial_y^\beta\,\widetilde{K}(x, y)$ is a linear
combination of finitely many expressions of the kind
\[
\bigl[ \partial_x^{\alpha'}\, \partial_y^{\beta'}\, K \bigr] (\psi(x), \psi(y))\cdot Q_{\alpha', \beta'}\,,
\]
where $\alpha'$, $\beta'$ are multi-indices with $1\le |\alpha'| \le |\alpha|$, $1 \le |\beta'| \le |\beta|$,
and $Q_{\alpha', \beta'}$ is a certain polynomial expression of partial derivatives of order at most
$\max(|\alpha|, |\beta|)$ of coordinate functions of $\psi$. In particular, if $K$ is $C^{k, k}$-smooth, then
so is $\widetilde K$. Since the maxima of higher order derivatives in the
definition of the norm $\| K \|_{L, Q, k}$ are multiplied by higher negative powers of $L$,
we conclude that for every compact $Q\subset V$ and $L\ge 1$,
\[
\| \widetilde{K} \|_{L, Q, k} \le C(\psi, Q, k)\, \| K \|_{L, \psi(Q), k}\,,
\]
where $ C(\psi, Q, k) $ is some factor depending on $\displaystyle \max_{|\gamma|\le k}\, \max_Q |\partial^\gamma \psi|$.

\medskip
Next, let $C^K_x$ be the matrix with the entries
$C^K_x(i, j) = \partial_{x_i}\, \partial_{y_j}\, K(x, y) $. Then
\[ \det C_{x}^{\widetilde K} = ( \det (D\psi)(x) )^2\, \det C_{\psi (x)}^K\,.\]

\medskip
One immediate consequence of these observations is that
\begin{itemize}\item
{\em the local controllability  of $K$
can be verified after any $C^2$-change of variables $\psi$,
and moreover, the corresponding constants at $x$
will change only by bounded factors depending on the first and second derivatives of $\psi$ and $\psi^{-1}$ at
$x$ and $\psi (x)$ respectively}.
\end{itemize}

\medskip
Now, let us see what the $C^2$-change of variable $\psi$ does to translation invariant scaling limits.
Let $z=\psi(z')\in U$. Assume that we have a sequence of kernels $K_L$ such that the corresponding
scaled kernels $K_{z, L}(u, v) = K_L(z+L^{-1}u, z+L^{-1}v)$ converge to $k(u-v)$, where $k$ is
a continuous function on $\bR^m$. Assume that for some $r>0$, there is a closed ball $\bar B = \bar B(z, r)\subset U$ such that
\[
\sup_{L>1}\, L^{-1} \max_{\bar B \times \bar B}\, (|\nabla_x K_L| + |\nabla_y K_L| ) = \mathfrak M < \infty\,.
\]
Let $u', v'\in\bR^m$. Then, for sufficiently large $L$, we have
\[
\bigl| \psi (z' + L^{-1}u') - \psi (z') - \tfrac1{L}(D\psi)(z') u' \bigr|  \le C(\psi)\, L^{-2} |u'|^2\,,
\]
and similarly for $v'$, where the constant $C(\psi)$ depends only on the bounds for the
second partial derivatives of $\psi$ in an arbitrarily small (but fixed) neighbourhood of $z'$.
Moreover, if $u'$ and $v'$ are fixed and $L$ is large, then the points
\[ \psi (z' + \tfrac1{L} u'), \quad  z + \tfrac1{L} (D\psi)(z') u'\,, \]
together with similar two points taken with $v'$ instead of
$u'$, belong to the ball $B$. So we obtain
\begin{multline*}
\bigl| \widetilde{K}(z'+\tfrac1{L}u', z'+\tfrac1{L}v') - K(z + \tfrac1{L} (D\psi)(z') u', z+ \tfrac1{L} (D\psi)(z') v') \bigr| \\[10pt]
\le L \cdot \mathfrak M \cdot C(\psi) L^{-2} (|u'|^2+|v'|^2) \to 0, \qquad {\rm as\ } L\to\infty\,.
\end{multline*}
Since $K(z + \frac1{L} (D\psi)(z') u', z + \tfrac1{L} (D\psi)(z') v')$ converge to $k((D\psi)(z') (u'-v'))$, we conclude that
$ \widetilde{K}(z'+L^{-1}u', z'+L^{-1}v') $ converge to $\widetilde{k}(u'-v')$, where
$\widetilde{k} (u') = k((D\psi)(z')(u'))$.

\medskip
Since a non-degenerate linear change of variable on the space  side corresponds
to a non-degenerate linear change of variable and renormalization on the Fourier side, the spectral measures
$\rho$ and $\widetilde\rho$, corresponding to $k$ and $\widetilde k$ respectively, do or do not have atoms
simultaneously.  This shows that
{\em the Gaussian parametric ensembles $f_L$ on $U$ and $ \widetilde{f}_L = f_L \circ \psi $
on $V=\psi^{-1}U$ are or aren't tame simultaneously}.
Finally, the corresponding limiting Gaussian functions $F_z$ and $\widetilde{F}_{z'}$
are related by $ \widetilde{F}_{z'} = F_z \circ (D\psi)(z') $, whence,
$\bar\nu_{\widetilde F_{z'}}(z') = |\det (D\psi)(z')|\, \bar\nu_{F_z} (z)$.

\subsection{Possibility to verify tameness in charts}\label{subsect:tameness-verific}

From the previous discussion it becomes clear why it suffices to check that
$ f_L \circ \pi_\alpha$ is tame on $U_\alpha$  for some atlas $\bA=(U_\alpha, \pi_\alpha)$
to be sure that
$f_L\circ\pi$ is tame on $U$ for any chart $\pi\colon U\to X$. Indeed, take any compact
$Q\subset \pi (U)$ and cover it by a finite union of open charts $ \bigcup_j \pi_{\alpha_j} (U_{\alpha_j})$.
Then we can choose compact sets $Q_j \subset Q \cap \pi_{\alpha_j} (U_{\alpha_j}) $ so that
$\bigcup_j Q_j = Q$. However, on each $Q_j$ the computations of all relevant quantities
in the charts $(U, \pi)$ and  $(U_{\alpha_j}, \pi_{\alpha_j})$ give essentially the same results
(up to bounded factors) because all partial derivatives of order $1$ and $2$ of the
transition mappings $\pi_{\alpha_j}^{-1} \circ \pi $ and $\pi^{-1}\circ \pi_{\alpha_j}$ are bounded
on $\pi^{-1}(Q_j)$ and $\pi_{\alpha_j}^{-1}(Q_j)$ respectively.

\medskip
If the atlas $\mathbb A$ has uniformly bounded distortions, our observations show that for every point $x\in X$,
{\em all} computations in {\em all} charts $(U, \pi)$ of $\mathbb A$ such that $x\in \pi (U)$
yield essentially the same results. Thus, for every point
$x\in X$, we can compute the relevant quantities in its own chart from $\mathbb A$ (the most convenient
one) without affecting the existence of uniform bounds for them, but, of course, affecting the best
possible values of those bounds.

\subsection{Completing the proof of Theorem~\ref{thm:manifold}}\label{subsect:proof-mflds}

Take two charts $\pi_1\colon U_1\to X$ and $\pi_2\colon U_2\to X$  and
consider the corresponding Gaussian parametric ensembles $f_{1, L}=f_L\circ \pi_1$ and
$f_{2, L}=f_L\circ \pi_2$ on $U_1$ and $U_2$ respectively.
For every $x\in \pi(U_1)\cap \pi(U_2) \subset X$, we have
\[
\bar\nu_1 (\pi_1^{-1}(x)) = \bigl| \det\bigl( [D(\pi_2^{-1} \pi_1)](\pi_1^{-1}(x)) \bigr) \bigr| \, \bar\nu_2 (\pi_2^{-1}(x))
\]
in the sense that if one side is defined, then so is the other and the equality holds.
Therefore, the push-forwards
$ (\pi_1)_* (\bar\nu_1\,\rd\vol) $ and $ (\pi_2)_* (\bar\nu_2\,\rd\vol) $ coincide on
$ \pi(U_1)\cap \pi(U_2) $, which allows us to define a Borel measure $n_\infty$ on $X$ unambiguously
and to justify the formula for its density with respect to any volume $\vol_X$ on $X$ compatible with the smooth
structure.

\medskip
The only thing that remains to do to establish Theorem~\ref{thm:main_thm} as stated, is to show that
\[
\lim_{L\to\infty} \cE^* \Bigl\{ \Bigl| L^{-m} \int_X \phi\, \rd n_L - \int_X \phi\, \rd n_\infty \Bigr| \Bigr\} = 0\,.
\]
The standard partition of unity argument allows us to reduce the problem to the case
when the support of the test function $\phi$ is contained in one chart $\pi (U)$.
Hence, the desired result
would be an immediate consequence of Theorem~\ref{thm:manifold} applied to the pull-back measures $\pi^* n_L$ and
the test-function $\phi\circ\pi$, if not for one minor nuisance: the pull-back to $U$ of a component counting measure
of $(f_L)$ on $X$ by the chart mapping $\pi$ may fail to be a component counting measure of $f_L\circ\pi$ on $U$
because the connected components of $Z(f_L)$ on $X$ that stretch outside $\pi (U)$ may get truncated
or split into several pieces when mapped to $U$ by $\pi^{-1}$. So the pull-back $\pi^*n_L$ may have
mass less than $1$ on some connected components of $f_L\circ \pi$ that stretch to the boundary of $U$.
We circumvent this difficulty by noticing that the closed support $\spt(\phi\circ\pi)$  is contained in $U$.
Thus, if we ``beef up'' the measure of each ``defective'' component $\gamma$ by adding an appropriate positive multiple
of a point mass at any point $u\in\gamma\setminus\spt(\phi\circ\pi)$, the pull-back $\pi^* n_L$ will turn
into a component counting measure
$n_L'$ but the total integral of $\phi\circ\pi$ will not be affected in any way. Now we can just apply Theorem~\ref{thm:main_thm}
to $n_L'$ instead of $\pi^* n_L$ and reach the desired conclusion. \hfill $\Box$

\begin{appendices}

\numberwithin{equation}{section}
\numberwithin{lemma}{section}
\numberwithin{theorem}{section}
\numberwithin{definition}{section}
\renewcommand{\thesection}{\Alph{section}}

\setcounter{subsection}{0}
\setcounter{section}{0}


\section{Smooth Gaussian functions}\label{App-A}

In this appendix, we collect well-known facts about smooth Gaussian
functions, which have been used throughout this paper.
Our smooth Gaussian functions will be defined on open subsets of $\bR^m$.
For a topological space $X$, by $\BB(X)$ we denote the Borel $\sigma$-algebra
generated by all open subsets of $X$.
As everywhere else in the paper, all Hilbert spaces are real and separable
and all Gaussian random variables have zero mean.

\subsection{The space $C^k(V)$}\label{A_subsect_C^k}

Let $V\subset\bR^m$ be an open set.
For $k\in\mathbb Z_+$, we denote by $C^k(V)$ the space of $C^k$-smooth functions on $V$.
The topology in $C^k(V)$ is generated by the seminorms\footnote{
The reader should be aware that the same notation was used in the main text for the seminorm
in $C^{k, k}(V\times V)$. This should not lead to a confusion because the functions measured
in these seminorms have different numbers of variables.
}
\[
\| g \|_{Q, k} = \max_Q\, \max_{|\alpha|\le k}\, \bigl| \partial^\alpha g \bigr|
\]
where $Q$ runs over all compact subsets of $V$. If $Q_n$ is an increasing sequence of
compact subsets of $V$ such that every compact set  $K\subset V$ is contained in each $Q_n$ with $n\ge n_0(K)$,
then the {\em countable}
family of the seminorms $\| g \|_{Q_n, k} $, $n=1, 2, $ \ldots \,, gives the same topology,
so $C^k(V)$ is metrizable. Since it is separable as well, every open set in $C^k(V)$
can be written as a countable union of ``standard neighbourhoods''
\[
\bB (Q, g_0, \e) = \bigl\{ g\in C^k(V)\colon \| g-g_0 \|_{Q, k} < \e \bigr\}.
\]

\medskip
We will need two simple lemmas.
\begin{lemma}\label{lemma:A_Borel1}
The Borel $\sigma$-algebra $\BB = \BB (C^k(V))$
coincides with the least $\sigma$-algebra
on $C^k(V)$ containing all ``intervals''
$ I(x; a, b)=\bigl\{ g\in C^k(V)\colon a\le g(x) < b \bigr\} $,
i.e., $\BB$ is generated by point evaluations $g\mapsto g(x)$.
\end{lemma}

\noindent{\em Proof of Lemma~\ref{lemma:A_Borel1}}: Denote
by $\BB'$ the least $\sigma$-algebra on $C^k(V)$ containing all intervals
$ I(x; a, b) $. We need to show that the $\sigma$-algebras $\BB$ and
$\BB'$ coincide. Since the mapping $C^k(V)\ni g \mapsto g(x) \in\bR$ is
continuous and, thereby, measurable, every interval $I(x; a, b)$ is Borel, that is
$\BB' \subset \BB$.

To show that $\BB \subset \BB'$, it suffices to check that
every standard neighourhood $\bB (Q, g_0, \e)$ belongs to $\BB'$,
or, which is the same, that the mapping $C^k(V)\ni g \mapsto \| g-g \|_{Q, k}$ is
$\BB$-measurable. Since for every fixed $x\in V$ and every multiindex
$\alpha$ with $|\alpha|\le k$, the mapping $g\mapsto \partial^\alpha g(x)$ can be represented
as a pointwise (on $C^k(V)$) limit of finite linear combinations of point
evaluations, it is measurable as well. It remains to note that
\[
\| g  - g_0 \|_{Q, k}
= \sup_{x\in Q'}\, \max_{|\alpha|\le k}\,
\bigl| \partial^\alpha g(x) - \partial^\alpha g_0(x) \bigr|,
\]
where $Q'$ is any countable dense (in $Q$) subset of $Q$. \mbox{}\hfill $\Box$

\begin{lemma}\label{lemma:A_Borel2}
$C^k(V)$ is a Borel subset of $C(V)$.
\end{lemma}

\noindent{\em Proof of Lemma~\ref{lemma:A_Borel2}}:
Take any function $\phi_1\in C_0^\infty(B)$, where $B$ is the unit ball in $\bR^m$,
put $\phi_j = j^m \phi (jx)$
and consider the sequence of continuous mappings
\[
C(V)\ni g \mapsto g * \phi_j \in C^k(V_{-1/j}).
\]
Note that $g\in C^k(V)$ if and only if $g*\phi_j$ converge in
$C^k$ uniformly on every compact set $Q\subset V$. Taking a countable exhaustion
$Q_n$ of $V$ and choosing $j(n)$ so that $Q_n \subset V_{-1/j(n)}$, we get the
representation
\[
C^k(V) = \bigcap_{q\ge 1}\, \bigcap_{n \ge 1}\, \bigcup_{j>j(n)}\, \bigcap_{s', s'' > j}
\bigl\{ g\in C(V)\colon \| g*\phi_{s'} - g* \phi_{s''} \|_{j,\, Q_n} < \tfrac1{q} \bigr\}.
\]
Clearly, the RHS is Borel in $C(V)$ (since each ``basic set'' on the RHS is open in $C(V)$).
\mbox{}\hfill $\Box$

\subsection{The definition and basic properties of $C^k$-smooth Gaussian
functions}\label{A_subsect_basic_defs}

\begin{definition}\label{def:A_smooth_random}
{\rm Let $ (\Omega, \mathfrak S, \cP) $ be a probability space.
The function $f\colon V \times \Omega \to \bR $ is a Gaussian function on $V$ if

\medskip\par\noindent {\rm (i)} for each $x\in V$, the mapping
$\omega\mapsto f(x, \omega)$ is measurable as a mapping from
$\bigl( \Omega, \mathfrak S \bigr)$ to $\bigl( \mathbb R, \BB(\mathbb R) \bigr)$;

\medskip\par\noindent {\rm (ii)} for each finite set of points $x_1, \ldots , x_n\in V$
and for each $c_1, \ldots , c_n\in\mathbb R$, the sum $\sum_j c_j f(x_j, \omega)$ is a Gaussian
random variable (maybe, degenerate).

\medskip\noindent Let $k\in \bZ_+$. The Gaussian function $f$ is called {\em $C^k$-smooth}
(or just $C^k$) if

\medskip\par\noindent {\rm (iii)} for almost every
$\omega\in\Omega$, the function $x\mapsto f(x, \omega)$
belongs to the space $C^k(V)$.
}
\end{definition}

\noindent Removing a subset of zero probability from $\Omega$,
we may (and will) just demand that the function $x\mapsto f(x, \omega)$
is in $C^k(V)$ for all $\omega\in\Omega$.

\medskip Every $ C^k$-Gaussian function $f$ generates two mappings
\[
\Omega\ni\omega \mapsto f(\, \cdot \,, \omega)\in C^k(V) \quad
{\rm and} \quad  V\ni x\to f(x, \, \cdot \,) \in L^2(\Omega, \mathcal P)\,.
\]
With some abuse of notation, we denote these mappings by the same letter
$f$.

\begin{lemma}\label{lemma:A_basic}
Suppose that $f$ is a $C^k$-smooth Gaussian function on $V$. Then

\smallskip\par\noindent{\rm (a)}
the mapping
$ f\colon \bigl( V\times\Omega, \BB(V) \times \mathfrak S) \to
\bigl( \mathbb R, \BB(\mathbb R) \bigr) $ is measurable;

\smallskip\par\noindent{\rm (b)}
the mapping
$f\colon \bigl( \Omega, \mathfrak S \bigr)\to \bigl( C^k(V), \BB(C^k(V)) \bigr) $
is measurable;

\smallskip\par\noindent{\rm (c)}
the mapping $f\colon V\to L^2 (\Omega, \mathcal P)$ is $C^k$-smooth.
\end{lemma}

\noindent{\em Proof of Lemma~\ref{lemma:A_basic}}:

\smallskip\par\noindent{\rm (a)}
We partition $V$ into countably many Borel sets $V_j$ of diameter
$\le 1/n$ each, fix an arbitrary collection of points $x_j\in V_j$,
and define a function $f_n\colon (V \times \Omega) \to \bR $ by
\[
f_n(x, \omega) = f(x_j, \omega) \qquad {\rm for\ } x\in V_j\,.
\]
The mappings $f_n\colon \bigl( V\times\Omega, \BB(V) \times \mathfrak S \bigr) \to \bigl( \bR, \BB(\bR) \bigr) $
are measurable and $\displaystyle f = \lim_{n\to\infty} f_n$ pointwise on $V\times\Omega$.

\medskip\par\noindent{\rm (b)}
It is an immediate consequence of Lemma~\ref{lemma:A_Borel1} combined with fact that,
for every $x\in V$, the mapping
$ \omega \mapsto f(x, \omega)$ is measurable.

\medskip\par\noindent{\rm (c)} Recall that if a sequence $\xi_n\colon \Omega\to\bR$
of Gaussian random variables converges pointwise to $\xi$, then $\xi$ is also a
Gaussian random variable. It follows that for every multiindex $\alpha$ with $|\alpha|\le k$, the
mapping $(x, \omega)\mapsto \partial^\alpha f(x, \omega)$ is a continuous Gaussian function.
Since the pointwise convergence of Gaussian random variables yields convergence in
in $L^2(\Omega, \cP)$, we see that the mapping
\[ V\ni x \mapsto \partial^\alpha f(x,\, \cdot\,) \in L^2(\Omega, \cP) \] is continuous
and gives the corresponding partial derivative of the mapping
$V\ni x \mapsto f(x, \,\cdot\,)$ considered as a function on $V$ with values
in $L^2(\Omega, \cP)$. \mbox{} \hfill $\Box$

\begin{definition}\label{def:A_equivalence}
{\rm
Let $f$ be a $C^k$ Gaussian function on $V$. Let
$\gamma_f \stackrel{\rm def}= f_* \mathcal P$ be the push-forward of the
probability measure $\mathcal P$ to $C^k(V)$ by $f$. We say that
two $C^k$ Gaussian functions $f_1$ and $f_2$ are {\em equivalent} if $\gamma_{f_1}=\gamma_{f_2}$.
We do not distinguish between equivalent Gaussian functions.
}
\end{definition}

In principle, we can
forget about the original probability space $ (\Omega, \mathfrak S, \cP)$
and consider the probability space $\bigl( C^k(V), \BB(C^k(V)), \gamma_f \bigr)$
and the mapping \[ V\times C^k(V)\ni (x, g) \mapsto g(x) \in \bR \] instead.
We can go one step further and remove any Borel subset of $\gamma_f$-measure
$0$ from $C^k(V)$ in this representation.

\subsection{Positive-definite kernels}\label{A_subsect_pos-def}

Let $f$ be a $C^k$ Gaussian function on $V$. Let
$K_f(t,s) \stackrel{\rm def}= \mathcal E \bigl\{ f(t) f(s) \bigr\}$
be the corresponding covariance kernel. It is a positive-definite symmetric
function\footnote{That is, the symmetric matrix
$ \bigl( K_f(x_i, x_j) \bigr)_{i, j =1}^n $ is positive definite for
every choice of $x_1, \ldots , x_n\in V$.} on $V\times V$.
The function $f$ is uniquely determined by $K_f$ up to equivalence.
Indeed, since a Gaussian distribution in $\bR^n$ is determined by its covariance matrix,
this fact is evident for the sets of the form
\[
S = \bigl\{ g\in C^k(V)\colon \bigl( g(x_1),\, \ldots , \, g(x_n) \bigr)\in B \bigr\}
\]
where $x_1, \, \ldots \,, x_n\in V$ and $B\in\BB(\bR^n)$. The general case follows immediately
because the fact that $\BB(C^k(V))$ is generated by point evaluations implies that every set
$S\in \BB(C^k(V))$ can be approximated by sets of such kind up to an arbitrary small $\gamma_f$-measure.

\medskip
Next, we observe that if $g$ is a continuous Gaussian function on $V$ with $K_g=K_f$,
then $g$ is equivalent, as a continuous Gaussian function, to the Gaussian function
$\widetilde{f}\colon \Omega \stackrel{f}\to C^k(V) \hookrightarrow C(V)$.
The function $\widetilde f$ generates a measure $\gamma_{\widetilde f}$ on $C(V)$:
\[
\gamma_{\widetilde f}(S) = \gamma_f (S\cap C^k(V))\,, \qquad  S\in\BB(C(V))\,.
\]
Furthermore,
\[
K_{\widetilde f}(x, y) = \cE \bigl\{ \widetilde f(x) \widetilde f(y) \bigr\}
= \cE \bigl\{ f(x) f(y) \bigr\} = K_f (x, y) = K_g (x, y)\,.
\]
Therefore, by the previous remark, $\gamma_{\widetilde f} = \gamma_g$.
In this situation, almost surely, $g$ is a $C^k$ Gaussian function. Indeed,
by Lemma~\ref{lemma:A_Borel2}, $C(V)\setminus C^k(V) \in\BB(C(V))$, whence,
\[
\gamma_g (C(V)\setminus C^k(V)) = \gamma_{\widetilde f} (C(V)\setminus C^k(V))
= \gamma_f (\emptyset) = 0\,.
\]
That is, {\em any continuous Gaussian function whose covariance kernel coincides with the
one of a $C^k$ Gaussian function, almost surely, is a $C^k$ function itself.}

\medskip Also observe that since the mapping \[ V\ni x \mapsto f(x, \, \cdot\,)\in L^2(\Omega, \cP) \] is
$C^k$,
the partial derivative $\partial_x^\alpha \partial_y^\beta K_f(x, y)$ exists and is continuous on
$V\times V$ for any multiindices $\alpha, \beta$ with $|\alpha|, |\beta|\le k$. Moreover,
\[
\partial^\alpha_x \partial^{\beta}_y K_f(x, y)
= \mathcal E \bigl\{ \partial^\alpha_x f(x)\, \partial^{\beta}_y f(y) \bigr\}.
\]

\subsection{From positive-definite kernels to reproducing kernel Hilbert spaces}
\label{A_subsect_RKHS}

In this section, we shall only assume that we are given
a continuous positive-definite symmetric kernel $K$ on $V\times V$.
We shall describe a canonical construction of the Hilbert space
$\cH = \cH(K)$ of continuous functions on $V$ such that $K$ is the {\em reproducing kernel}
in that space, that is, for every $g\in\cH$ and every $x\in V$, we have
$g(x) = \langle g, K_x \rangle_{\cH}$ where $K_x(y)=K(x,y)$.

\medskip
Consider the linear space $\mathcal L$ of all mappings
$\mathfrak h\colon V\to\bR$ such that $\mathfrak h(x)\ne 0$ only for
finitely many $x\in V$. Define the semi-definite scalar product on $\mathcal L$
by
\[
\langle \mathfrak h_1, \mathfrak h_2 \rangle
= \sum_{x, y\in V} K(x, y)\, \mathfrak h_1 (x) \mathfrak h_2(y)
\]
(this sum is actually finite). Since $K$ is positive-definite, we have
$\langle \mathfrak h, \mathfrak h \rangle \ge 0$ for every $\mathfrak h \in\mathcal L$.
Define the Hilbert seminorm on $\mathcal L$ by
$\| \mathfrak h \| = \sqrt{\langle h, h \rangle }$. Then $\langle \, \cdot \, , \, \cdot \, \rangle$ and
$\|\, \cdot \,\|$ become a nondegenerate scalar product and the associated Hilbert norm
on $\mathcal L / \mathcal L_0 $ where $\mathcal L_0$ is the linear subspace of $\mathcal L$
consisting of all $\mathfrak h\in\mathcal L$ with $\| \mathfrak h \| = 0$.
Let $H$ be the Hilbert space completion of the pre-Hilbert space
$\bigl( \mathcal L / \mathcal L_0, \langle \, \cdot \, , \, \cdot \, \rangle \bigr) $.

For $x\in V$, denote by $h_x$ the vector in $H$ corresponding to the function
\[
\mathfrak h_x(y) =
\begin{cases}
0, & y\ne x \\
1, & y=x.
\end{cases}
\]
Note that $\| h_x - h_y \|^2 = K(x, x) + K(y, y) - 2K(x, y) \to 0 $ as $y\to x$, so
the mapping $V\ni x \mapsto h_x \in H$ is continuous. Since
$\operatorname{span}\{h_x\colon x\in V \}$ is dense in $H$ by construction and since
for every countable dense subset $V'\subset V$, the set
$\{ h_x\colon x\!\in\!V'\}$ is dense in $\{ h_x\colon x\!\in\!V\}$,
$H$ is separable.

Now, define a linear map $\Phi\colon H\to C(V)$ by
$\Phi [h](x) = \langle h, h_x \rangle $, $h\in H$. If $\Phi[h]=0$, then
$\langle h, h_x \rangle = 0$ for all $x\in V$, so $h=0$.
Thus, we can identify $H$ with a linear subspace $\cH = \Phi (H)$ of $C(V)$.
Note also that \[ \Phi [h_x] (y) = \langle h_x, h_y \rangle = K(x, y) = K_x(y), \] so
$h_x$ is identified with $K_x$. Transferring the scalar product
from $H$ to $\cH$, we turn $\cH$ into a Hilbert space of continuous functions
on $V$ with the reproducing kernel $K$.

Observe, finally, that such a Hilbert space is unique. Indeed, if $\cH_1\subset C(V)$
is another Hilbert space of continuous functions with the same reproducing
kernel $K$, then the linear span $\cH_0$ of the functions $K_x$, $x\in V$,
is contained and dense in $\cH_1$ with respect to the Hilbert norm in $\cH_1$
(because if $g\in\cH_1$ is orthogonal to all $K_x$ in $\cH_1$, then
$ g(x) = \langle g, K_x \rangle_{\cH_1} = 0 $ for all $x\in V$, whence, $g=0$)
and for every pair of functions
\[
g_1 = \sum_{\tt finite} a_x K_x, \quad g_2 = \sum_{\tt finite} b_y K_y
\]
in $\cH_0$, we have
\[
\langle g_1, g_2 \rangle_{\cH_1} = \sum_{x, y} K(x, y) \, a_x b_y
= \langle g_1, g_2 \rangle_{\cH}\,.
\]
Thus the identity mapping $\cH_0 \to \cH_0$ can be extended to a bijective
isometry $\cH \to \cH_1$. Let now $g'\in\cH_1$ be the image of $g\in\cH$
under this isometry. Then
\[
g'(x) = \langle g', K_x \rangle_{\cH_1} = \langle g, K_x \rangle_{\cH} = g(x)\,,
\quad x\in V\,,
\]
so $\cH_1$ and $\cH$ consist of exactly the same functions on $V$ and are endowed
with the same scalar product.  \mbox{}\hfill $\Box$

\medskip
We end this section with a useful observation. Let $\{e_k\}$ be an arbitrary
orthonormal basis in $\cH$. For every $g\in\cH$, we put
$\widehat{g}(k) = \langle g, e_k \rangle_{\cH}$.
Then the Fourier series
$ \sum_k \widehat{g}(k) e_k$ converges to $g$ in $\cH$.
For every $y\in V$, we have
\begin{multline*}
\Bigl| g(y) - \sum_{1\le k \le N} \widehat{g}(k) e_k(y)  \Bigr|
= \Bigl| \Bigl\langle g - \sum_{1 \le k \le N} \widehat{g}(k) e_k, K_y \Bigr\rangle \Bigr| \\
\le \Bigl\| g - \sum_{1\le k \le N} \widehat{g}(k) e_k \Bigr\|_{\cH}\,
\| K_y \|_{\cH} \to 0 \quad {\rm as\ } N\to \infty\,.
\end{multline*}
Since $\| K_y \|_{\cH} = \sqrt{K(y, y)}$ is a continuous function of $y$ on $V$,
this yields the locally uniform convergence of the Fourier series $\sum_k \widehat{g}(k) e_k$ to $g$.

\smallskip Taking $g=K_x$ and observing that $\langle K_x, e_k \rangle = e_k(x)$, we
conclude that for every $x, y \in V$, we have
\[
\sum_k e_k (x) e_k(y) = K(x, y)\,.
\]

\subsection{Canonical series representation of continuous Gaussian functions}
\label{A_subsect_canon_series}

Let $H_0$ be any Gaussian subspace of $L^2(\Omega, \cP)$ and let
$ V\ni x \mapsto f_x\in H_0$
be a continuous mapping such that for every $x\in V$, the random variable
$f_x$ is Gaussian. The corresponding covariance
kernel $K(x, y) = \cE\{ f_x f_y \} = \langle f_x, f_y \rangle_{L^2(\Omega, \cP)} $
is also continuous. Let $H$ be the closed linear span of $\{ f_x \}_{x\in V}$
in $L^2(\Omega, \cP)$. It is a Gaussian subspace of $L^2(\Omega, \cP)$. For
$h\in H$, define $\Phi[h](x)= \langle h, f_x \rangle_{L^2(\Omega, \cP)}$.
Note that $\Phi[h]\in C(V)$ and $\Phi[h]=0$ if and only if $h=0$.
Also, $\Phi[f_x]=K(x, \,\cdot\,) = K_x$. Thus, $\cH = \{ \Phi[h]\colon h\in H \}$
is a linear subspace of $C(V)$ and if we endow it with the scalar product
$ \langle \Phi[h], \Phi[h'] \rangle_{\cH} = \langle h, h' \rangle_{L^2(\Omega, \cP)} $,
it will become a Hilbert space $\cH (K)$ of continuous functions with the reproducing kernel
$K$.

Now, take any orthonormal basis $\{ e_j \}$ in $\cH$ and choose $\xi_j\in H$ such that
$e_j = \Phi[\xi_j]$. Note that
\[
\langle \xi_i, \xi_j \rangle_{L^2(\Omega, \cP)} = \langle e_i, e_j \rangle_{\cH} =
\begin{cases}
0, & i\ne j \\
1, & i=j\,,
\end{cases}
\]
so $\xi_j$ are orthogonal and, thereby, independent standard Gaussian.
For every $x\in V$, we have
\[
f_x = \sum_j \langle f_x, \xi_j \rangle\, \xi_j
= \sum_j \Phi[\xi_j](x) \xi_j = \sum_j e_j(x)\xi_j\,.
\]
The upshot is that,
\begin{itemize}
\item
{\em given any Gaussian subspace $H_0\subset L^2(\Omega, \cP)$,
any continuous mapping $x\mapsto f_x$ from $V$ to $H_0$, and any orthonormal basis
$e_j$ in the reproducing kernel Hilbert space $\cH(K)$, where $K(x, y)=\cE\{f(x)f(y)\}$,
we can define independent standard Gaussian real variables on
$(\Omega, \mathfrak S, \cP)$ such that $ f_x = \sum_j \xi_j\, e_j (x) $ for all
$x\in V$}.
\end{itemize}

Assume now that we start with a continuous Gaussian function $f$ with some
underlying probability space $ (\Omega, \mathfrak S, \cP)$. Applying the above
construction to the induced mapping $x\mapsto f_x = f(x, \,\cdot\,)$, we get
\begin{equation}\label{eq:A_series0}
f(x, \omega) = \sum_j \xi_j(\omega)\, e_j(x)
\quad {\rm in\ } L^2(\Omega, \cP) \ {\rm for\ all\ } x\in V\,.
\end{equation}
Implementing $\xi_j$ as some everywhere defined functions on $\Omega$ and taking into
account that $L^2(\Omega, \cP)$-convergence yields convergence in probability, we have, in particular, that
\begin{equation}\label{eq:A_series}
{\rm for\ every \ } x\in V, \quad \sum_{j=1}^n \xi_j(\omega)\, e_j(x) \to f(x)
\quad {\rm in\ probability\ as\ } n\to\infty\,.
\end{equation}

\medskip
Now, put $ X_j = \xi_j e_j(x) $, $S_n = \sum_{j=1}^n X_j$ and note that for every compact $Q\subset V$,
the random variables $X_j$, $S_n$ and $S=f$ can be viewed as random vectors in the Banach space $C(Q)$. The random variables
$X_j$ are symmetric and independent, and~\eqref{eq:A_series} means that for every point evaluation
functional $z_x\in C(Q)^*$ given by $\langle z_x, g \rangle = g(x) $ for $x\in Q$,
we have $\langle z, S_n \rangle \to \langle z_x, S \rangle $ in probability.
By the classical Ito-Nisio theorem, which we will recall in the next section,
the series $\sum_j X_j$ converges to $S$ in $C(Q)$. Thus,
\begin{itemize}
\item
{\em the canonical series
representation~\eqref{eq:A_series0} actually converges in $C(V)$}.
\end{itemize}

\subsection{The Ito-Nisio theorem}\label{A_Ito-Nisio}

Let $\mathfrak X$ be a separable Banach space. An $\mathfrak X$-valued random variable
on a probability space $(\Omega, \mathfrak S, \cP)$ is just a measurable
mapping from $(\Omega, \mathfrak S, \cP)$ to $(\mathfrak X, \BB (\mathfrak X)\,) $.
Everywhere below, $X_j$ is a sequence of independent $\mathfrak X$-valued random
variables, $S$ is an $\mathfrak X$-valued random variable on the same probability
space, and $S_n = \sum_{j\le n} X_j$. We denote by $\| \cdot \|$ the norm in $\mathfrak X$,
and by $\langle  \cdot  ,  \cdot  \rangle$ the natural coupling of the dual space $\mathfrak X^*$ and $\mathfrak X$.

\medskip
First, we recall a classical

\medskip\par\noindent{\bf P. Lev\'y's lemma: }
{\em If $S_n$ converges to $S$ in probability, then $S_n$ converges to $S$
almost surely.}

\medskip\noindent{\em Proof of Lev\'y's lemma}: We will check that,  for almost every $\omega\in\Omega$,
$S_n(\omega)$ is a Cauchy sequence.
Take $\e\in (0, \tfrac12)$ and take $m$ so large that
\[
\cP \bigl\{ \| S_k - S \| \ge \tfrac12\, \e \bigr\} < \tfrac12\, \e\,,
\qquad k\ge m\,.
\]
Then take any positive integer $n>m$. For
$k=m, \ldots , n$, let
\[
A_k = \bigl\{\omega\in\Omega\colon \|S_\ell-S_m\| \le 2\e
\ {\rm for\ all\ } \ell = m,\, \ldots \,, k-1, \ {\rm but\ }
\| S_k - S_m \| > 2\e \bigr\}.
\]
Note that the events $A_k$ are disjoint,
and each $A_k$ is independent of $S_n - S_k$.
Also, if $\omega\in A_k$ and $\| S_n-S_k\|\le\e$, then
\[
\| S_n - S_m \| \ge \| S_k - S_m\| - \|S_n - S_k \| > 2\e - \e = \e\,.
\]
Furthermore, since
\[
\bigl\{ \| S_k - S_n \| > \e \bigr\} \subset
\bigl\{ \| S_k - S \| > \tfrac12\, \e \bigr\} \bigcup \bigl\{ \| S_n - S \| > \tfrac12\, \e \bigr\}\,,
\]
we have
\[
\cP \bigl\{ \| S_k - S_n \| \le \e \bigr\} = 1 - \cP \bigl\{ \| S_k - S_n \| > \e \bigr\} > 1-\e\,.
\]
Therefore,
\begin{multline*}
(1-\e)\, \cP\{ \max_{m\le k \le n} \| S_k - S_m \| > 2\e \}
= (1-\e)\, \sum_{k=m}^n P(A_k) \\
< \sum_{k=m}^n \cP \{ \| S_n - S_k \| \le \e \} \cdot \cP (A_k) \\
= \sum_{k=m}^n \cP \{ A_k  \ {\rm and\ } \| S_n-S_k \| \le \e   \} \\
\le \cP \{ \| S_n - S_m \|>\e \} \le \e\,.
\end{multline*}
Since $n$ and $\e$ are arbitrary, we see that, with
probability $1$, $S_n$ is Cauchy, so it converges almost surely. Clearly, the almost sure limit
and the limit in probability must be the same. \mbox{} \hfill $\Box$

\medskip We call a subset $Z\subset\mathfrak X^*$ {\em normalizing} if it is
countable and $ \| x \| = \sup\, \{ \langle z, x \rangle \colon z\in Z \} $
(then automatically $Z$ is contained in the unit ball of $\mathfrak X^*$).
Now, we can state the part of the Ito-Nisio theorem that we need:

\medskip\noindent{\bf Ito-Nisio theorem:} {\em Suppose that the random variables
$X_j$ are symmetric and that there exists a normalizing set $Z\subset\mathfrak X^*$
such that $\langle z, S_n \rangle \to \langle z, S \rangle $ in probability
for every $z\in Z$. Then $S_n\to S$ in $\mathfrak X$ almost surely.}

\medskip\noindent{\em Proof\,}: By P. Lev\'y's lemma, it is enough to show that
$S_n$ converges to $S$ in probability. First of all, note that the Borel $\sigma$-algebra
$\BB (\mathfrak X)$ coincides with the $\sigma$-algebra $\BB' (\mathfrak X)$
generated by the events
$\{x\colon \langle z, x \rangle \in [a, b) \}$, $z\in Z$, $a, b \in\bR$. Indeed,
for every $x_0\in\mathfrak X$, the mapping
$x\mapsto \| x-x_0 \| = \sup_Z \| \langle z, x \rangle - \langle z, x_0 \rangle \|$
is $\BB'(\mathfrak X)$-measurable. Thus, every open ball in $\mathfrak X$ is
$\BB'(\mathfrak X)$-measurable. Since $\mathfrak X$ is separable, every open
set in $\mathfrak X$ is $\BB'(\mathfrak X)$-measurable, so
$\BB (\mathfrak X)\subset \BB'(\mathfrak X)$. The inverse inclusion is obvious.

Next, we show that $S_n$ and $S-S_n$ are independent for every $n$. We need to check that
\[
\cP \{ S_n \in C_1, S-S_n \in C_2 \} =
\cP \{ S_n \in C_1 \} \cdot \cP\{ S-S_n \in C_2 \}
\]
for every $C_1, C_2\in\BB(\mathfrak X)$. Since $\BB (\mathfrak X) = \BB'(\mathfrak X)$,
it suffices to check this for the events of the form
\[
C=\bigl\{ \bigl( \langle z_1, x \rangle,\, \ldots \,, \langle z_q, x \rangle \bigr) \in B \bigr\},
\quad B\in\BB(\bR^q), \  z_1, \ldots , z_q \in Z\,,
\]
in which case it follows from the independence of $S_m-S_n$ and $S_n$ for $m>n$ and
the fact that, for $m\to\infty$,
\[
\bigl( \langle z_1, S_m-S_n \rangle,\, \ldots \,, \langle z_q, S_m-S_n \rangle \bigr)
\to \bigl( \langle z_1, S-S_n \rangle,\, \ldots \,, \langle z_q, S-S_n \rangle \bigr) \quad
{\rm in\ probability\,}.
\]

For a set $\mathfrak X' \subset \mathfrak X$, denote
\[
\mathfrak X'_{+\e} = \bigcup_{x\in\mathfrak X'} B(x, \e)\,.
\]
We claim that {\em for every finite set $\mathfrak X' \subset \mathfrak X$, there exists
a finite ``separating set'' of functionals $Z' \subset Z$ such that}
\[
\max_{z\in Z'} | \langle z, y' - y''\rangle | > \e\,,
\quad {\rm whenever\ } y', y'' \in \mathfrak X'_{+\e}
\ {\rm and\ } \| y'-y'' \|> 8\e\,.
\]
Indeed, consider all differences $x'-x''$ with $x', x''\in\mathfrak X'$ and for each
of them choose $z=z(x', x'')\in Z$ such that
$| \langle z, x'-x'' \rangle | \ge \tfrac12\, \| x' - x'' \|$.
Since $y', y''\in \mathfrak X'_{+\e}$, we can find $x', x''\in\mathfrak X'$
so that $\|x'-y'\|, \|x''-y''\|<\e$. Then $\|x'-x''\|>6\e$. Taking $z=z(x', x'')$, we get
\[
| \langle z, y' -y'' \rangle | > | \langle z, x' -x'' \rangle | -2\e
\ge 3\e - 2\e = \e\,,
\]
proving the claim.

\medskip
Now, comes the crux of the proof. Suppose that $A$ and $B$ are $\mathfrak X$-valued independent
random variables and $A$ is symmetric. Then, for every finite $\mathfrak X' \subset \mathfrak X$
and every $\e > 0$, we can write
\[
\cP \bigl\{ A\notin \bigl( \tfrac12 (\mathfrak X' - \mathfrak X') \bigr)_{+\varepsilon} \bigr\}
\le \cP \bigl\{ A + B \notin \mathfrak X'_{+\e }\bigr\}
+ \cP \bigl\{ -A + B \notin \mathfrak X'_{+\e }\bigr\}
= 2\, \cP \bigl\{ A + B \notin \mathfrak X'_{+\e }\bigr\}\,.
\]
The inequality here is due to the observation that
\[
{\rm if\ } a, b \in\mathfrak X \quad  {\rm and\ } a+b, -a + b \in \mathfrak X'_{+\e},
\quad {\rm then\ }
a\in\bigl( \tfrac12 (\mathfrak X' - \mathfrak X' ) \bigr)_{+\e}\,.
\]
The equality follows at once from the symmetry of $A$ and the independence of $A$ and $B$.

To finish the proof, we take $\e>0$ and let $x_1$, $x_2$, \ldots be a countable dense
set in $\mathfrak X$. Put $\mathfrak X'_N = \{x_1, \ldots , x_N\}$. Since
$(\mathfrak X'_N)_{+\e} \uparrow \mathfrak X$ as $N\to\infty$, we have
$\cP \{ S\notin (\mathfrak X'_N)_{+\e} \} < \e$ for large enough $N$.
We fix such $N$ and, to simplify notation, let $\mathfrak X' = \mathfrak X'_N$.
Since $S_n$ and $S-S_n$ are independent, and $S_n$ is symmetric, we can use them
as $A$ and $B$ in the argument above and get
$ \cP \bigl\{ S_n\notin \bigl( \tfrac12 (\mathfrak X' - \mathfrak X') \bigr)_{+\e} \bigr\} < 2\e $ for all $n$.
Let $Z'\subset Z$ be a finite separating set for
$\mathfrak X' \cup \frac12 (\mathfrak X' - \mathfrak X')$.
Then for every $n$,
\begin{multline*}
\cP \bigl\{ \| S_n - S \| > 8\e \bigr\} \\
\le \cP \bigl\{ S_n \notin \bigl( \tfrac12 (\mathfrak X' - \mathfrak X') \bigr)_{+\e} \bigr\}
+ \cP \bigl\{ S \notin \mathfrak X'_{+\e}\bigr\}
+ \cP \bigl\{ \max_{z\in Z'} | \langle z, S_n \rangle - \langle z, S \rangle | > \e \bigr\} \\
\le 3\e +
\sum_{z\in Z'} \cP \bigl\{ | \langle z, S_n \rangle - \langle z, S \rangle | > \e \bigr\}\,.
\end{multline*}
Since each term in the finite sum on the RHS tends to $0$ as $n\to\infty$ and $\e$ can be taken
as small as we want, the desired convergence in probability follows. \mbox{}\hfill $\Box$

\subsection{The local behavior of continuous Gaussian functions}\label{A_subsect_top_support}

Suppose that $f$ is a continuous Gaussian function on $V$ with the covariance
kernel $K$. As before, we denote by $\gamma_f$ the corresponding Gaussian
measure on $C(V)$.  The (closed) set $S(f)$ of functions $g\in C(V)$ for which
$\cP \{ f\in U \} = \gamma_f(U) >0 $ for every open neighbourhood $U$ of $g$ in $C(V)$
is called the {\em topological support of the measure} $\gamma_f$.

The following lemma gives a simple and useful description of the topological support of $\gamma_f$:
\begin{lemma}\label{lemma:A_top_suppport} $S(f) = \operatorname{Clos}_{C(V)} \cH(K)$.
\end{lemma}

\noindent{\em Proof of Lemma~\ref{lemma:A_top_suppport}}: First, we show that
for every $g\in\cH(K)$, every compact $Q\subset V$, and every $\e>0$, we have
$\cP\{ \| f-g \|_{C(Q)} < \e \} > 0$. We choose an orthonormal basis $\{ e_j \}$ in
$\cH(K)$ so that $g=te_1$ for some $t\in\bR$, and represent $f$ as $\sum_j \xi_j e_j$
where $\xi_j$ are independent Gaussian random variables. Since, by the Ito-Nisio theorem,
the series converges in $C(Q)$ with probability $1$, there exists $N=N(\e)$ such
that
\[
 \| \sum_{j>N }\xi_j e_j \|_{C(Q)} = \| f - \sum_{j\le N} \xi_j e_j \|_{C(Q)} < \tfrac12 \,\e
\]
with positive probability. Next, we choose $\eta$ so small that
\[
\eta \cdot \sum_{j \le N} \| e_j \|_{C(Q)} < \tfrac12\, \e\,.
\]
Now, suppose that
\[
\| \sum_{j>N }\xi_j e_j \|_{C(Q)} < \e\,, \quad \xi_1 \in (t-\eta, t+\eta)\,, \quad {\rm and} \quad \xi_2,\, \ldots \,, \xi_n\in(-\eta, \eta)\,.
\]
Then
\begin{multline*}
\| f-g \|_{C(Q)} \le |\xi_1-t|\,\| e_1 \|_{C(Q)} + \sum_{2\le j \le N} |\xi_j |\, \| e_j \|_{C(Q)} + \bigl\| \sum_{j>N} \xi_j e_j \bigr\|_{C(Q)} \\
\le \eta \sum_{j\le N} \| e_j\|_{C(Q)} + \bigl\| \sum_{j> N} \xi_j e_j \bigr\|_{C(Q)} < \e\,.
\end{multline*}
Hence,
\begin{multline*}
\cP \{ \| f-g \|_{C(Q)} < \e \} \ge \cP \bigl\{ \| \sum_{j > N} \xi_j e_j \|_{C(Q)}
< \tfrac12 \, \e \bigr\} \cdot \\
\cdot \cP \bigl\{ \xi_1 \in (t-\eta, t+\eta),  \ \xi_2, \ldots , \xi_N \in (-\eta, \eta) \bigr\} > 0\,.
\end{multline*}
Thus, $\cH(K)\subset S(f)$. Since $S(f)$ is closed in $C(V)$, we get
$ S(f) \supset \operatorname{Clos}_{C(V)} \cH(K) $.

\medskip To show the converse, assume that $g\in C(V)$ and
\[
\cP \{ \| f-g \|_{C(Q)} < \tfrac12\, \e \} = p > 0\,.
\]
We fix an orthonormal basis $\{e_j\}$ in $\cH(K)$ and choose $N$ so large that
\[
\cP \{ \| f - \sum_{j\le N} \xi_j e_j  \|_{C(Q)} > \tfrac12\, \e \} <  \tfrac12\, p\,.
\]
Then
\[
\cP \{ \| g - \sum_{j\le N} \xi_j e_j  \|_{C(Q)} < \e \} \ge p - \tfrac1{2}\, p > 0\,.
\]
Since
\[
\sum_{j\le N} \xi_j e_j \in \cH(K)\,,
\]
we conclude that the $\e$-neighbourhood of $g$ in $C(Q)$ intersects $\cH(K)$.
Since $\e$ and $Q$ are arbitrary, we see that
$ S(f) \subset \operatorname{Clos}_{C(V)} \cH(K) $.
\mbox{}\hfill $\Box$

\subsection{Fernique's theorem}\label{A_subsect_Fernique}

The next result we state was proven by Fernique and independently by Landau and Shepp.
It allows one to pass from some very weak estimates for various norms and semi-norms of
Gaussian functions to almost as strong bounds for tails as possible in principle.

\medskip\par\noindent{\bf Fernique's theorem}:
{\em Let $X$ be a random variable with values in a Banach space $\mathfrak X$, and let
$\{ \phi_j \} \subset \mathfrak X^*$ be an at most countable set of linear functionals on $\mathfrak X$
such that, for every choice of finitely many $\phi_j$'s, the joint distribution of $\phi_j(X)$ is Gaussian.
Suppose that, for some $\la>0$ and $\mu<\frac12$,
\begin{equation}\label{eq:Fernique1}
\cP\bigl\{ \sup_j |\phi_j(X)| \ge \la \bigr\} \le \mu\,.
\end{equation}
Then, for all $t\ge 1$,
\begin{equation}\label{eq:Fernique2}
\cP\bigl\{ \sup_j |\phi_j(X)| \ge \la t \bigr\} \le e^{-at^2}
\end{equation}
with a positive constants $a$ depending only on $\mu$.
}

\medskip Here, we present Fernique's original proof, which is short and elegant.
Landau and Shepp~\cite{LS} gave a different proof based on the Gaussian isoperimetry. The advantage of
the latter proof is that it does not need a priori assumption $\mu<\frac12$ and gives the
optimal RHS of~\eqref{eq:Fernique2}, which is $\Phi \bigl( t\Phi^{-1}(\mu) \bigr)$, where
\[
\Phi (s) = \sqrt{\frac2{\pi}}\, \int_s^\infty e^{-x^2/2}\, \rd x\,.
\]

\medskip\noindent{\em Proof of Fernique's theorem}:
Without loss of generality, we assume that $\la=1$.
Let $\Omega_n (t)$ be the event $\bigl\{ \sup_{1\le j \le n} |\phi_j (X)| > t \bigr\} $, $n\in\bN\cup \{\infty\}$.
We need to estimate $\cP\{\Omega_\infty (t)\}$  assuming that  $\cP\{\Omega_\infty (1)\} \le \mu $.
Since $\Omega_n (t) \subset \Omega_m (t)$ for $m\le n$, and $\Omega_\infty (t) = \bigcap_{n\ge 1} \Omega_n(t)$, it will
suffice to prove estimate~\eqref{eq:Fernique2} for every finite $n$ with constants $A$ and $a$ independent of $n$.

In what follows, we fix $n\in\mathbb N$, and put $\phi = \bigl( \phi_1(X),\, \ldots \,, \phi_n(X) \bigr)$. This is
a finite-dimensional Gaussian random vector. We let
$\displaystyle \| \phi \| \stackrel{\rm def}= \max_{1\le j \le n} |\phi_j(X)| $.
Fernique's proof is based on the following classical observation: {\em if $\psi$
is an $n$-dimensional Gaussian vector, which has the same distribution as $\phi$ and which
is independent of $\phi$, then $\frac1{\sqrt2} (\phi+\psi)$ and $\frac1{\sqrt2}(\phi-\psi)$ are two Gaussian
vectors, which have the same distribution as $\phi$ and which are independent of each other}.
The proof of this statement reduces to a routine verification of coincidence of all relevant
covariance matrices.

Now, take $t>0$ and $\tau>0$ and write
\begin{align*}
\cP \bigl\{ \| \psi \| \le \tau \bigr\} \cdot \cP \bigl\{ \| \phi \| > t \bigr\}
&= \cP \bigl\{ \| \tfrac1{\sqrt{2}} (\phi-\psi) \| \le \tau \bigr\} \cdot
\cP \bigl\{ \| \tfrac1{\sqrt{2}} (\phi+\psi) \| > t \bigr\} \\[10pt]
&= \cP \bigl\{ \| \tfrac1{\sqrt{2}} (\phi-\psi) \| \le \tau, \
\| \tfrac1{\sqrt{2}} (\phi+\psi) \| > t \bigr\} \\[10pt]
&\le \cP \bigl\{ \| \phi \| > \tfrac1{\sqrt{2}} (t-\tau), \
\| \psi \| > \tfrac1{\sqrt{2}} (t-\tau) \bigr\}\\[10pt]
&= \bigl( \cP\bigl\{ \| \phi \| > \tfrac1{\sqrt{2}} (t-\tau) \bigr\} \bigr)^2.
\end{align*}
Letting $\tau=1$ and recalling that $\cP\bigl\{ \| \psi \| \le 1 \bigr\} \ge 1-\mu$, we get
\[
\cP \bigl\{ \| \phi \| > t \bigr\} \le \frac1{1-\mu}  \left( \cP\bigl\{ \| \phi \| > \tfrac1{\sqrt{2}} (t-\tau) \bigr\} \right)^2.
\]
Put $p(t)=\cP \bigl\{ \| \phi \| > t \bigr\}$. This is a non-increasing function of $t$, which satisfies
\begin{align*}
p(t) &\le \frac1{1-\mu} p^2\bigl( \tfrac1{\sqrt{2}} (t-1) \bigr)\,, \qquad t\ge 1\,, \\
p(1) &\le \mu\,.
\end{align*}
Let
\[
t_k = \frac{(\sqrt{2})^{k+1}-1}{\sqrt{2}-1}\,, \qquad k\ge 0\,,
\]
that is, $t_0=1$, and $t_k = \frac1{\sqrt{2}} (t_{k+1}-1)$. Then, by induction on $k$, we have
\[
p(t_k) \le (1-\mu)\, \Bigl( \frac{\mu}{1-\mu} \Bigr)^{2^k}\,.
\]
Since
\[
t_{k+1}^2 = \Bigl[ \frac{(\sqrt{2})^{k+2}-1}{\sqrt{2}-1} \Bigr]^2
< \frac{2^{k+2}}{(\sqrt{2}-1)^2} = \frac4{(\sqrt{2}-1)^2}\, 2^k\,,
\]
we see that for $t_k \le t \le t_{k+1}$,
\[
p(t) \le p(t_k) < e^{-at^2} \quad {\rm with} \ a=\frac14\, (\sqrt{2}-1)^2\, \log\frac{\mu}{1-\mu} > 0\,,
\]
completing the proof. \hfill $\Box$

\subsection{Kolmogorov's theorem}\label{A_Kolmogorov}

Here, we formulate a version of the classical Kolmogorov's theorem for $C^{k, k}$ kernels.
Let $k\in\bN$ and let, as before, $V\subset\bR^m$ be an open set.

\begin{definition}\label{def:A_C^k}
{\rm We say that a symmetric function $K\colon V\times V \to \bR$
belongs to $C^{k, k}(V \times V)$ if all partial derivatives of $K$ including
at most $k$ differentiations in $x$ variables and at most $k$ differentiations
in $y$ variables exist and are continuous on $V\times V$ (in which case,
the order of differentiations does not matter and we can denote these derivatives
$\partial_x^\alpha \partial_y^\beta K(x, y)$ as usual).
}
\end{definition}

\medskip\noindent{\bf Kolmogorov's theorem:}
{\em Let $k\in\bN$. Suppose that $K\colon V\times V\to\bR$ is a positive definite symmetric function
of class $C^{k, k}(V\times V)$
and, in addition, that
\[
N_{V, k}(K) \stackrel{\rm def}=
\max_{|\alpha|, |\beta|\le k}\, \sup_{x, y\in V}\, \bigl| \partial_x^\alpha\, \partial_y^\beta K(x, y) \bigr| < \infty\,.
\]
Then there exists a (unique up to an equivalence) $C^{k-1}$ Gaussian function $f$ on $V$ with the covariance kernel $K$.

Moreover,
for every $\gamma\in (0, 1)$ and every closed ball $\bar B\subset V$,
we have
\[
\cE \bigl\{ \| f \|_{\bar B,\, k-1+\gamma} \bigr\} \le C(\bar B, V, k, \gamma)\, \sqrt{N_{V, k} (K)}\,.
\]
}

\medskip Note that since every compact set $Q\subset V$ can be covered by a finite union of
closed balls contained in $V$, the latter estimate immediately implies that, for any compact set $Q\subset V$,
\[
\cE \bigl\{ \| f \|_{Q,\, k-1} \bigr\} \le C(Q, V, k)\, \sqrt{N_{V, k} (K)}\,.
\]
The same is true for the H\"older norm, but the cover should be chosen carefully so that any two sufficiently close
points $x, y\in Q$ are covered by a single ball.
Then the resulting
bound on $Q$ depends on both the bounds on the balls and the geometry of the cover. We will never
need to estimate the H\"older norms on any compact set other than a ball, so we will not go into
the details here.

\medskip
It is also worth noting that
in the assumptions of Kolmogorov's theorem we use that $N_{V, k} (K) < \infty $ instead of the more natural for a function
defined on an open set assumption that $N_{Q, k} (K) < \infty $ for every compact set $Q\subset V$. This allows us to
reduce the number of nested compact sets we need to choose before doing any estimate. Of course, this replacement it is harmless.

\subsection{Proof of Kolmogorov's theorem}\label{A_subsect_proof_Kolmogorov}
To prove Kolmogorov's theorem we will use a ``convolution approach''. As far as high
order derivatives are concerned, this approach allows one to pass to the
limits in a family of covariance kernels easier than the more usual approach based on nets (see, for instance,~\cite[Section~3.1]{H}).

We split the proof into several steps.

\subsubsection{}\label{A-subsubsect:proof-Kolmogorov-prelimin}
As we have seen in~\ref{A_subsect_RKHS}, there exists a separable Hilbert space
$\cH$ and a continuous mapping $V\ni x \mapsto f_x \in \cH$ such that $K(x, y) =
\langle f_x, f_y \rangle$. Without loss of generality, we assume that $\cH$ is a
Gaussian subspace of $L^2(\Omega, \cP)$, where $(\Omega, \mathfrak S, \cP)$ is a
probability space. Our first task is to implement the mapping $x\mapsto f_x$
as a $\BB (V) \times \mathfrak S$-measurable function of $(x, \omega)$.

\medskip
We start with implementing each $f_x$ as an everywhere defined function
on $\Omega$. Then we pick a compact exhaustion $Q_n$ of $V$,
a sequence $\e_n>0$ with $\sum_n \e_n^2 < \infty$,
and choose $\rho_n>0$ so small that $\| f_x - f_y \|^2_{L^2(\Omega, \cP)}<\e_n^2$
for all $x\in Q_n$, $y\in V$ with $|x-y|<\rho_n$. We fix a countable partition
of $V$ into Borel sets $V_{j, n}$ of diameter less than $\rho_n$ each, choose
some point $x_{j, n}$ in every $V_{j, n}$ and put $f_n(x, \omega) = f_{x_{j, n}}(\omega)$
if $x\in V_j$. Then, for every $x\in Q_n$, $\| f_n - f_x \|^2_{L^2(\Omega, \cP)} < \e_n^2$,
so for each $t>0$, we have
\[
\cP \bigl\{ | f_n(x, \omega) - f_x(\omega) |> t \bigr\}
< t^{-2}\, \e_n^2\,.
\]
Since $\sum_n \e_n^2 < +\infty$ and $Q_n$ exhaust $V$, the functions $f_n(x, \, \cdot \,)$ converge to $f_x$
both $\cP$-almost surely and in $L^2(\Omega, \cP)$.

Let now $\displaystyle E = \{(x, \omega)\colon \lim_{n\to\infty} f_n (x, \omega)
\ {\rm exists\,}\} $. Since $f_n$ is $\BB(V)\times\mathfrak S$-measurable, so
is $E$. Also, for every $x\in V$, we have $\cP \{\omega\colon (x, \omega)\notin E\}=0$.
Thus,
\[
f(x, \omega) \stackrel{\rm def}=
\begin{cases}
\lim_{n\to\infty} f_n (x, \omega), & (x, \omega)\in E \\
0, & {\rm otherwise}
\end{cases}
\]
is a measurable representation of the mapping $x\mapsto f_x$.

\subsubsection{}\label{A-subsubsect:proof-Kolmogorov-prelimin2}
Denote $F_\omega (x) = f(x, \omega)$. By Fubini, for every compact $Q\subset V$,
\[
\cE \Bigl\{ \int_Q | F_\omega |^2\, \rd\vol \Bigr\}
= \int_Q \| f_x \|^2_{L^2(\Omega, \cP)}\, \rd\vol (x) \le \max_{x\in Q} K(x, x) \cdot \vol Q
\le N_{V, k}(K)^2\, \vol Q\,.
\]
Thus $F_\omega \in L^2_{\tt loc}(V)$ for every $\omega\in\Omega_1\subset\Omega$ with
$\cP(\Omega_1)=1$. Replacing $f(x, \omega)$ by $f(x, \omega) \done_{\Omega_1} (\omega)$,
we will assume that  $f(x, \omega)$ is such that $F_\omega\in L^2_{\tt loc}(V)$ for all
$\omega\in\Omega$.

\subsubsection{}\label{A-subsubsect:proof-kolmogorov-convolutions}
Next, we note that for every $\phi\in C_0^\infty (B(r))$, the convolution
in the $x$ variable
\[
(f *_x \phi)(x, \omega) \stackrel{\rm def}= (F_\omega * \phi)(x)
\]
is a $C^k$ (actually, $C^\infty$)
Gaussian function on $V_{-r}$ for every $\omega\in\Omega$. The only non-trivial part of this claim
is the Gaussian distribution property. To see it, observe that,
as an element of $\mathcal H\subset L^2(\Omega, \cP)$,
\[
(f *_x \phi )(x,\, \cdot\, ) = \int_{B(r)} f_{x+y}\, \phi (y)\, \rd\vol (y)\,.
\]
The integral on the RHS can be understood as the usual Riemann integral of a
continuous $L^2(\Omega, \cP)$-valued function, and hence, it
can be approximated in $L^2(\Omega, \cP)$ by finite Riemann sums
$\sum_j c_j f_{x+y_j} \in \cH$ and, therefore, lies in $\cH$ itself.
In what follows, we write $f * \phi$ instead of $f *_x \phi$ and view
$f * \phi$ as a random Gaussian function.

\subsubsection{}\label{A-subsubsect:proof-Kolmogorov-convolution-estimates}
We shall need a few estimates for $f * \phi$ and its derivatives
$\partial^\alpha (f * \phi) = f * \partial^\alpha \phi $ for $|\alpha| \le k-1$.
First of all, by Fubini,
\begin{multline*}
\cE \bigl\{ \bigl[ \partial^\alpha (f * \phi) (z) \bigr]^2 \bigr\}
= \iint_{B(r)\times B(r)} K(z+x, z+y) \partial^\alpha \phi (x) \partial^\alpha \phi (y)\,
\rd\vol (x)\, \rd\vol (y) \\
= \iint_{B(r)\times B(r)} \partial^\alpha_x \partial^\alpha_y K(z+x, z+y) \phi (x) \phi (y)\,
\rd\vol (x)\, \rd\vol (y)\,.
\end{multline*}
The expression on the right is trivially bounded by $\| \phi \|^2_{L^1}\, N_{V, k}(K)$.

If, in addition, the function $\phi$ has zero integral mean, we can
improve our trivial bound to $C r^2\, \| \phi \|^2_{L^1}\, N_{V, k}(K)$.
To see this, we put
\[
E_\alpha (z; x, y) = \partial_x^\alpha \partial_y^\alpha K(z+x, z+y)
- \partial_x^\alpha \partial_y^\alpha K(z, z+y) - \partial_x^\alpha \partial_y^\alpha K(z+x, z)
+ \partial_x^\alpha \partial_y^\alpha K(z, z)
\]
and note that by ``bilinear'' Lagrange's Mean-Value Theorem,
$ | E_\alpha (z; x, y) | \le C r^2\, N_{V, k}(K) $.
Then, writing
\[
\partial_x^\alpha \partial_y^\alpha K(z+x, z+y)  = - \partial_x^\alpha \partial_y^\alpha K(z, z)
+ \partial_x^\alpha \partial_y^\alpha K(z, z+y) + \partial_x^\alpha \partial_y^\alpha K(z+x, z)
+ E_\alpha (z; x, y)
\]
and integrating in $x$ and $y$ against $\phi(x)\phi(y)\, \rd\vol (x)\,\rd\vol (y)$, we obtain
\[
\cE \bigl\{ \bigl[ \partial^\alpha (f * \phi) (z) \bigr]^2 \bigr\}
= \iint_{B(r)\times B(r)} E_\alpha (z; x, y)
\,\rd\vol (x)\, \rd\vol (y) \le  C r^2\, \| \phi \|^2_{L^1}\, N_{V, k}(K)\,.
\]

\subsubsection{}\label{A-subsubsect:entropy-bound}
We shall also need the following ``entropy bound'':
\begin{lemma}[entropy bound]\label{lemma:A-entropy}
Let $r>0$. Let $g$ be a continuous Gaussian function on $V$ and
$\psi$ be any $C^\infty_0(B(r))$-function. Then $g*\psi$ is a continuous
Gaussian function on $V_{-r}$ and for every two compact sets $Q, Q'\subset V$
such that $Q_{+r}\subset Q'$, we have
\[
\cE \bigl\{ \| g * \psi \|_{C(Q)} \bigr\}
\le 5 \| \psi \|_{L^1} \cdot
\sqrt{1+\log\frac{\| \psi \|_{L^\infty} \vol Q'}{\| \psi \|_{L^1}}}\,
\cdot
\,\sqrt{\sup_{Q'} \cE \{|g|^2\}}\,.
\]
\end{lemma}

\noindent{\em Proof of the entropy bound}:
Without loss of generality, we assume that $\sup_{Q'} \cE \{ |g|^2 \} = 1$. Then
for every $x\in Q'$, $g(x)$ is a Gaussian random variable with
$\cE\{ g(x)^2 \} \le 1$. Hence,
\[
\cE e^{\frac14 g(x)^2} \le
\frac1{\sqrt{2\pi}}\, \int_{\bR} e^{\frac1{4} x^2}\, e^{-\frac1{2} x^2}
\, \rd x = \sqrt{2}\,.
\]
Take $\rho \ge \sqrt{2}$.
Noting that the function $\rho e^{-\frac14 \rho^2}$ decreases on $[\sqrt{2}, +\infty)$, we
estimate the convolution by
\begin{multline*}
(g * \psi)(x) = \int_{B(r)} g(x+y) \psi (y)\, \rd\vol (y) \\
= \int_{B(r) \cap \{|g|\le \rho\}} g(x+y) \psi (y)\, \rd\vol (y) +
\int_{B(r)\cap \{|g|>\rho\}} g(x+y) \psi (y)\, \rd\vol (y) \\
\le \rho \| \psi \|_{L^1} +
\rho e^{-\frac14 \rho^2} \| \psi \|_{L^\infty} \int_{Q'} e^{\frac14 g^2}\, \rd\vol\,,
\end{multline*}
so
\[
\cE \| g*\psi \|_{C(Q)}
\le \rho \bigl[ \| \psi \|_{L^1}
+ e^{-\frac14 \rho^2} \| \psi \|_{L^\infty} \sqrt{2} \vol Q'\bigr].
\]
Taking
\[
\rho = 2 \sqrt{1+\log\frac{\| \psi \|_{L^\infty} \vol Q'}{\| \psi \|_{L^1}}}
\]
(which is $\ge 2$ because $\| \psi \|_{L^1} \le \| \psi \|_{L^\infty} \vol B(r)
\le \| \psi \|_{L^\infty} \vol Q' $) and using that $2(1+\sqrt{2})<5$, we get
the desired bound. \mbox{}\hfill $\Box$

\subsubsection{}\label{A-subsubsect:proof-Kolmogorov-tw-expressions}
Now, we fix $\phi \ge 0$ in $ C_0^\infty (B(1)) $ with $\int \phi\, \rd\vol = 1$. For
$r>0$, let $\phi_r (x) = r^{-m} \phi (r^{-1}x)$ and note that
$\| \phi_r \|_{L^1} = 1 $ and $ \| \phi_r \|_{L^\infty} \le Cr^{-m} $ for all
$r>0$. Take a sequence $r_j = 2^{-j-1} $ and put
$f_j = f * \phi_{r_j} * \phi_{r_j}$. Then $f_j$ are $C^k$ Gaussian functions
on $V_{-2 r_j}$, and $f_j(x, \,\cdot\,) \to f_x$, as $j\to\infty$,
in $L^2(\Omega, \cP)$ for all $x\in V$.
Next, we fix a closed ball $\bar B = \bar B(x, r)\subset V$ and choose $j_0$ so large
that $\bar B(x, r+2r_{j_0})\subset V$.

Consider the series
\begin{equation}\label{eq:app-*}
f_{j_0} + \sum_{j\ge j_0} (f_{j+1}-f_j)\,.
\end{equation}
If we show that for every $\alpha$ with $|\alpha|\le k-1$, the expression
\[
\cE \| \partial^\alpha f_{j_0} \|_{C(\bar B)} + \sum_{j\ge j_0} \cE \| \partial^\alpha f_{j+1}
- \partial^\alpha f_{j} \|_{C(\bar B)}
\]
is bounded by $C(\bar B, j_0)\,\sqrt{N_{V, k}(K)}$,
and that for every $\alpha$ with $|\alpha|=k-1$ and every $\gamma\in (0, 1)$, the expression
\[
\cE \| \partial^\alpha f_{j_0} \|_{\bar B, \gamma} + \sum_{j\ge j_0} \cE \| \partial^\alpha f_{j+1}
- \partial^\alpha f_{j} \|_{\bar B, \gamma}
\]
is  bounded by   $C(\bar B, j_0, \gamma)\,\sqrt{N_{V, k}(K)}$, then we will
be done because then the series~\eqref{eq:app-*} will converge in $C^{k-1}(V)$ almost surely,
its sum will be a Gaussian function
$f$ with the covariance kernel $K$, and the desired bounds
for $\cE \| f \|_{\bar B, k-1+\gamma}$ will hold as well.

\subsubsection{}
For a multi-index $\alpha$ with $|\alpha|\le k$, we write $\partial^\alpha f_{j_0} = \partial^\alpha (f * \phi_{r_{j_0}}) * \phi_{r_{j_0}}$ and
note that the function $g = \partial^\alpha (f * \phi_{r_{j_0}})$ satisfies
$\cE \{ g(x)^2 \} \le N_{V, k}(K)$.
So Lemma~\ref{lemma:A-entropy} yields the bound $\cE \| \partial^\alpha f_{j_0 }\|_{C(\bar B)} \le C(\bar B, j_0)\, \sqrt{N_{V, k}(K)}$.
The interesting part is $\cE \| \partial^\alpha f_{j+1} - \partial^\alpha f_{j}\|_{C(\bar B)}$.
If $|\alpha|\le k-1$, writing
\[
\partial^\alpha f_{j+1}
- \partial^\alpha f_{j} = \partial^\alpha (f*(\phi_{r_{j+1}} - \phi_{r_{j}})) * (\phi_{r_{j+1}}+\phi_{r_{j+1}})
\]
applying the entropy bound with $g=\partial^\alpha (f*(\phi_{r_{j+1}} - \phi_{r_{j}}))$ and
$\psi = \phi_{r_{j+1}} + \phi_{r_{j}}$, and recalling that, by~\ref{A-subsubsect:proof-Kolmogorov-convolution-estimates},
$\cE\bigl\{ |g|^2 \bigr\} \le Cr_j^2\, N_{V, k}(K)$, we see that
\[
\cE \bigl\{ \| \partial^\alpha (f_{j+1}-f_j) \|_{C(\bar B)} \bigr\}
\le C\, r_j\,  \sqrt{1+\log (Cr_{j+1}^{-m})} \cdot \sqrt{N_{V, k}(K)}
\]
with some $C=C(\bar B)$. Since $\sum_j r_j\,  \sqrt{1+\log (Cr_{j+1}^{-m})} < \infty$, this takes care of the first of the series
in~\ref{A-subsubsect:proof-Kolmogorov-tw-expressions}
including the uniform norms of the derivatives of $f$ of order up to $k-1$.

\subsubsection{}
To get convergence of the series
\begin{equation}\label{eq:app-series}
\sum_{j\ge j_0} \cE \bigl\{ \| \partial^\alpha f_{j+1}
- \partial^\alpha f_{j} \|_{\bar B, \gamma} \bigr\}
\end{equation}
for a multi-index $\alpha$ with $|\alpha|=k-1$, we need the bound for
$ \cE \| \nabla \partial^\alpha f_{j+1} - \nabla \partial^\alpha f_{j} \|_{C(\bar B)} $.
Note that despite we still have convolutions with mean zero functions in the representation
of $ \nabla \partial^\alpha f_{j} $, we cannot use our trick from~\ref{A-subsubsect:proof-Kolmogorov-convolution-estimates}
because the kernel smoothness is totally exhausted.
Thus, we can use only the trivial estimate from~\ref{A-subsubsect:proof-Kolmogorov-convolution-estimates} without the factor $r_j$, and the entropy bounds yields
\[
\cE \bigl\{ \|  \nabla \partial^\alpha f_{j+1} - \nabla \partial^\alpha f_{j} \|_{C(\bar B)} \bigr\}
\le C \sqrt{1 + \log (Cr^{-m}_{j+1})} \cdot \sqrt{N_{V, k}(K)}\,.
\]
There is no hope to choose $r_j$ so that these terms will form a convergent series, so there is
no chance to show on this way that the $k-1$-st order derivatives are Lipschitz. Fortunately, we
do not need that much. All we really need is H\"older continuity.

Using a classical trick, we observe that for any function $h$
that is $C^1$ in some neighbourhood of $\bar B$, and for any
two points $x, y\in\bar B$, we have\footnote
{
We use the inequality $\min (a, b) \le a^{1-\gamma} b^\gamma $ valid for positive $a$ and $b$ and
for $\gamma\in (0, 1)$.
}
\[
| h(x) - h(y) | \le \min \bigl[ 2\| h \|_{C(\bar B)}, \| \nabla h \|_{C(\bar B)} |x-y| \bigr]
\le 2^{1-\gamma} \| h \|^{1-\gamma}_{C(\bar B)}\, \| \nabla h \|^\gamma_{C(\bar B)}|x-y|^\gamma
\]
By H\"older's inequality,
\begin{multline*}
\cE \Bigl\{ \|  \partial^\alpha f_{j+1} - \partial^\alpha f_{j} \|^{1-\gamma}_{C(\bar B)} \|
\cdot
 \|  \nabla\partial^\alpha f_{j+1} - \nabla\partial^\alpha f_{j} \|^{\gamma}_{C(\bar B)} \Bigr\} \\[5pt]
\le \Bigl( \cE  \|  \partial^\alpha f_{j+1} - \partial^\alpha f_{j} \|_{C(\bar B)} \Bigr)^{1-\gamma}
\cdot
 \Bigl( \cE  \| \nabla \partial^\alpha f_{j+1} - \nabla \partial^\alpha f_{j} \|_{C(\bar B)} \Bigr)^{\gamma}\,,
\end{multline*}
and, by the entropy bound, the RHS is
\[
\lesssim r_j^{1-\gamma} \bigl( 1 + \log ( Cr^{-m}_{j+1}  ) \bigr)\, \sqrt{N_{V, k}(K)}\,.
\]
Hence, the series~\eqref{eq:app-series} converges and the proof of Kolmogorov's theorem is complete. \mbox{}\hfill $\Box$

\subsection{Remarks to Kolmogorov's theorem}\label{A_subsect:remarks_Kolmogorov}

\subsubsection{}\label{A_subsubsect:remarks_Kolmogorov-Fernique}
Kolmogorov's theorem, as stated and proved, allows us to estimate $ \cE \| f \|_{\bar B, k-1+\gamma} $.
However, applying then Fernique's theorem, we immediately see that, in assumptions of Kolmogorov's theorem, we can estimate any moment
$ \cE \| f \|^p_{\bar B, k-1+\gamma} $ we want (that would be exactly as much as we use in this paper),
and even prove that the distribution tail
$  \cP \bigl\{ \| f \|_{\bar B, k-1+\gamma} > t \bigr\} $ at $t\to +\infty$ is Gaussian with controllable bounds.

Indeed, we take $\mathfrak X = C^{k-1}(\bar B)$, $X=f$, fix a countable dense set $B'\subset \bar B$, and put
\begin{align*}
\phi_{\alpha, x} (f) &= \partial^\alpha f(x), \qquad &|\alpha|\le k-1,\ x\in B'\,, \\
\phi_{\alpha, \gamma, x} (f) &= \frac{\partial^\alpha f(x) - \partial^\alpha f(y)}{|x-y|^\gamma}\,, \qquad &|\alpha|= k-1,\ x,y\in B'\,, x\ne y\,.
\end{align*}
Note that this is a countable system of linear functionals $\{ \phi_j \} \subset \mathfrak X^*$ satisfying the assumptions of Fernique's theorem,
and that $\| f \|_{\bar B, k-1+\gamma} = \sup_j |\phi_j(f)|$.
By Kolmogorov's theorem, there exists a positive constant $\la=\la(\bar B, V, k, \gamma)$ such that
\[
\cP \bigl\{ \| f \|_{\bar B, k-1+\gamma} > \la\cdot \sqrt{N_{V, k}(K)} \bigr\} < \tfrac14.
\]
Then Fernique's theorem tells us that
\[
\cP \bigl\{ \| f \|_{\bar B, k-1+\gamma} > t\, \la\cdot \sqrt{N_{V, k}(K)} \bigr\} <
e^{-at^2}\,, \qquad t\ge 1\,,
\]
whence,
\begin{equation}\label{eq:A-Kolmogorov-tails}
\cP \bigl\{ \| f \|_{\bar B, k-1+\gamma} > t \bigr\} <
C(B, V, k, \gamma) e^{-c(B, V, k, \gamma) t^2/N_{V, k}(K)}\,, \qquad t >0\,.
\end{equation}
In particular,
\[
\cE \bigl\{ \| f\|_{\bar B, k-1+\gamma}^p \bigr\}
\le C(B, V, k, \gamma)\, N_{V, k}^{p/2} (K)\,.
\]

\medskip
It is worth mentioning that one can also arrive at estimate~\eqref{eq:A-Kolmogorov-tails} directly after a certain modification of the proof of
Kolmogorov's theorem we gave.

\subsubsection{}
We have to distinguish between $C^k$ Gaussian functions on $U$ and Gaussian functions with $C^{k, k}(U\times U)$ covariance kernels:
the former are always the latter but, in general, not vice versa. However, by Kolmogorov's theorem, the continuous
Gaussian functions with $C^{k, k}(U\times U)$ covariance kernels fail to be in $C^k$ themselves just barely: they all
are in $C^{k-}(U) = \bigcap_{0<\gamma<1}\, C^{k-1+\gamma}(U)$.

\subsubsection{}
The ``convolution approach'' to Kolmogorov's theorem
allows one to approximate Gaussian
functions of finite smoothness by $C^\infty$ ones. This approximation can be used to establish some properties
of the kernel.

Using this idea, we will show now that every semi-norm $\|K\|_{Q, k}$ of a positive definite $C^{k, k}(U\times U)$ kernel
can be read from the ``diagonal'':
\[
\max_{|\alpha|, |\beta|\le k}\, \max_{x, y\in Q}\, \bigl| \partial^\alpha_x\, \partial^\beta_y\, K(x, y) \bigr|
=
\max_{|\alpha|\le k}\, \max_{x\in Q}\, \bigl| \partial^\alpha_x\, \partial^\alpha_y\, K(x, y)\big|_{y=x}\, \bigr|\,.
\]
Indeed, if $|\alpha|, |\beta|\le k-1$, then we can write
\begin{multline*}
\bigl| \partial^\alpha_x\, \partial^\beta_y\, K(x, y) \bigr|^2 = \bigl| \cE \bigl\{ \partial^\alpha f(x)\, \partial^\beta f(y) \bigr\} \bigr|^2 \\[10pt]
\le \cE \bigl\{ [ \partial^\alpha f(x)]^2 \bigr\}\, \cE \bigl\{ [ \partial^\beta f(y) ]^2 \bigr\}
= \bigl( \partial^\alpha_x\, \partial^\alpha_y\, K(x, y)\big|_{x=y} \bigr)\, \bigl( \partial^\beta_x\, \partial^\beta_y\, K(x, y)\big|_{y=x} \bigr)
\end{multline*}
for the $C^{k-}$ Gaussian function $f$ with the covariance kernel $K$, thus estimating the off-diagonal values by the square
root of the product of the two corresponding diagonal ones.
We cannot do the same estimate directly for the highest order derivatives, but we can consider the convolutions
$f * \phi$ that are infinitely smooth and get the inequality
\begin{equation}\label{eq:A_K_phi}
\bigl| \partial^\alpha_x\, \partial^\beta_y\, K_\phi (x, y) \bigr|^2
\le \bigl( \partial^\alpha_x\, \partial^\alpha_y\, K_\phi (x, y)\big|_{x=y} \bigr)\, \bigl( \partial^\beta_x\, \partial^\beta_y\, K_\phi (x, y)\big|_{y=x} \bigr)
\end{equation}
for the corresponding covariance kernels
\[
K_\phi (x, y) = \iint K(x+x', y+y') \phi(x')\phi(y')\, \rd\vol (x')\rd\vol (y')\,.
\]
Taking $\phi_1\in C_0(B(1))$ and $\phi (x) = \phi_r(x) = r^{-m}\phi_1(r^{-1}x)$, we can pass to
the limit
\[
\partial^\alpha_x\, \partial^\beta_y\, K_{\phi_r}(x, y) \to \partial^\alpha_x\, \partial^\beta_y\, K(x, y) \qquad \text{as\ } r\to 0\,,
\]
for any $|\alpha|, |\beta|\le k$, $x, y \in U$, we conclude that~\eqref{eq:A_K_phi} holds for $K$ as well.

\medskip Of course, here one can also work with the kernel directly, approximating the derivatives by finite difference ratios
and passing to the limit in some inequalities for long sums.

\subsubsection{}
The convolutions also facilitate convergence: if the kernels $K_\ell\in C^{k, k}(U\times U)$ are uniformly bounded on
compact subsets of $U\times U$ and converge pointwise to some kernel $K$ on $U\times U$, then $(K_\ell)_\phi \to K_\phi $
in $C^\infty (U_{-r}\times U_{-r})$ for any $\phi\in C^\infty_0(B(r))$. If we know, in addition, that for $|\alpha|, |\beta|\le k$,
the partial derivatives
$ \partial^\alpha_x\, \partial^\beta_y\, K_{\ell}(x, y)$ are uniformly locally bounded as well, we can use
elementary analysis to show that $K\in C^{k-1, k-1}(U\times U)$ and
$ \partial^\alpha_x\, \partial^\beta_y\, K_{\ell}(x, y) \to \partial^\alpha_x\, \partial^\beta_y\, K(x, y) $ for
$|\alpha|, |\beta|\le k-1$ uniformly on compact subsets of $U\times U$. However, in general, it is impossible to conclude
that $K\in C^{k,k}(U\times U)$. Surprisingly, this conclusion holds if the limiting kernel $K$ is translation invariant, i.e.,
$K(x, y) = \kappa (x-y)$ for some $\kappa\colon  \bR^m\to\bR$. This will be shown in the next section.

\subsection{Translation-invariant Gaussian functions}\label{A_subsect_tr-inv}

A continuous Gaussian function on $\bR^m$ is translation-invariant
if its covariance kernel $K(x, y)$ depends on $x-y$ only, i.e.,
$K(x, y) = \kappa (x-y) $ for some continuous positive definite $\kappa\colon \bR^m \to \bR$.
In this case, $\kappa$ can be written as a Fourier integral of some finite symmetric
positive Borel measure $\rho$ on $\bR^m$, i.e.,
\[
\kappa(x) = \int_{\bR^m} e^{2\pi {\rm i} (\la\cdot x)}\, \rd\rho(\la)\,.
\]
Consider the Hilbert space $L^2_{\tt H}(\rho)$ of all Hermitean ($h(-x)=\overline{h(x)}$)
functions $h\colon \bR^m\to\mathbb C$ with $\int |h|^2\, \rd\rho < \infty $. The standard
$L^2(\rho)$ scalar product $ \langle h_1, h_2 \rangle = \int h_1 \bar h_2\, \rd\rho$
is real on  $L^2_{\tt H}(\rho)$. Also, for every $x\in\bR^m$, the function
$f_x(\la) = e^{2\pi {\rm i}(\la\cdot x)}$ belongs to $L^2_{\tt H}(\rho)$ and
$\langle f_x, f_y \rangle = \kappa(x-y) $. Finally, the linear span of the functions $f_x$
is dense in $L^2_{\tt H}(\rho)$. Indeed, if $h\in L^2_{\tt H}(\rho)$, then
$\Phi [h] (x) = \langle h, f_x \rangle $ is the Fourier transform of the finite
Borel measure $h\,\rd\rho$. Hence, it vanishes identically only if $h=0$ $\rho$-a.e.\,.
Bringing all these observations together, we conclude that
\begin{itemize}
\item
{\em the Hilbert space $\cH(K)$ coincides with the
Fourier image $\mathcal F\, L^2_{\tt H}(\rho)$}.
\end{itemize}

\medskip
Now, we discuss the smoothness properties of translation invariant Gaussian functions and covariance
kernels. First of all, note that if $K(x, y)=\kappa(x-y)$, then
\[
\partial^\alpha_x\, \partial^\beta_y\, K (x, y)
= (-1)^{|\beta|} \bigl( \partial^{\alpha+\beta} \kappa \bigr)(x-y)\,.
\]
Thus, {\em $K$ is in $C^{k, k}(\bR^m\times\bR^m)$ if and only if $\kappa\in C^{2k}(\bR^m)$, that is, if and only if},
\begin{equation}\label{eq:app-2k_moment}
\int_{\bR^m} |\la|^{2k}\, \rd\rho(\la) < \infty\,.
\end{equation}

We end this section with a curious and quite useful observation:
\begin{itemize}\item
{\em if a sequence of positive definite
kernels $K_\ell\in C^{k, k}(U_\ell, U_\ell)$ with $U_\ell$ exhausting $\bR^m$ has a pointwise translation invariant limit
$ \kappa (x-y) $ and  $ \partial^\alpha_x\, \partial^\alpha_y\, K_{\ell}(x, y)\big|_{x=y=0} $ and stays bounded for $|\alpha|\le k$, then
$ \kappa\in C^{2k}(\bR^m) $}.
\end{itemize}

\medskip\noindent{\em Proof\,}: For $\phi\in C^\infty_0(\bR^m)$, $\phi(-x)=\phi (x)$ and put
$  K_\phi (x, y) = (\kappa * \phi * \phi)(x-y) $.
Since $\kappa= \cF \rho$ implies that $\kappa * \phi * \phi = \cF \rho_\phi$, where $\rd\rho_\phi = \widehat{\phi}\,^2\, \rd\rho$, we see
that
\begin{multline*}
(-1)^k\, \sum_{|\alpha|=k} \partial^{2\alpha} (\kappa * \phi * \phi) (0) =
(2\pi)^{2k}\, \sum_{|\alpha|=k} \int_{\bR^m} \la^{2\alpha_1}_1\, \ldots \, \la^{2\alpha_m}_m\, \rd\rho_\phi (\la) \\
= (2\pi)^{2k}\,  \int_{\bR^m} |\la|^{2k}\,\rd\rho_\phi (\la)\,.
\end{multline*}
If we know in advance that $\kappa\in C^{2k}(\bR^m)$, then the quantities $ \partial^{2\alpha} ( \kappa * \phi * \phi )(0) $,
$|\alpha|=k$, are uniformly bounded when
$\phi$ runs over even non-negative $C_0^\infty$ functions supported on a small ball centered at the origin and normalized by
$ \displaystyle \int_{\bR^m} \phi\,\rd\vol =1$.  Then, taking as before, $\phi_r (x) = r^{-m} \phi (r^{-1}x)$, letting $r\to 0$,
and applying Fatou's lemma, we get
\[
\int_{\bR^m} |\la|^{2k}\,\rd\rho(\la) \le \varliminf_{r\to 0}\, \int_{\bR^m}   |\la|^{2k}\,\rd\rho_{\phi_r}(\la) < \infty\,.
\]
Now, observe that the quantities  $ \partial^{2\alpha} (\kappa * \phi * \phi) (0) $ stay uniformly bounded even if $\kappa (x-y)$
is a pointwise limit of $C^{k,k}$ positive definite symmetric kernels $K_\ell (x, y)$ that are defined only in a neighbourhood
of the origin in $\bR^m \times\bR^m$ and that have uniformly bounded derivatives
$ \partial^\alpha_x\, \partial^\alpha_y\, K_{\ell}(x, y)\big|_{x=y=0} $. So, in this case, we still get
\[
\int_{\bR^m} |\la|^{2k}\,\rd\rho(\la) < \infty\,,
\]
completing the proof of our observation. \hfill $\Box\breve{}$

\section{Proof of the Fomin-Grenander-Maruyama theorem}\label{App-B}

Assuming that $\rho$ has no atoms, we need to show that if
$A\in\mathfrak S$ is a set satisfying $\gamma \bigl( (\tau_v A) \triangle A \bigr) = 0$ for every $v\in\bR^m$, then
$\gamma (A)$ is either $0$ or $1$. As before, we use the notation $ (\tau_v G)(u) = G(u+v)$, where $v\in\bR^m$ and $G\in X$.

\medskip
Since $\mathfrak S$ is generated by the intervals $I(u; a, b)$, given $\e>0$, we can take
finitely many points $u_1, \ldots , u_n\in\bR^m$ and a Borel set $B\subset \bR^n$ so that
$ \gamma \{ A \triangle P \} < \e $, where
\[
P = P(u_1, \ldots , u_n; B) \stackrel{\rm def}=
\bigl\{ G\in X\colon (G(u_1), \ldots , G(u_n))\in B\bigr\}\,.
\]
Without loss of generality, we may assume that the distribution of the
Gaussian vector $\bigl( G(u_1), \ldots , G(u_n) \bigr)$ is non-degenerate\footnote{
Otherwise one of the values, say, $G(u_n)$, is a linear combination
of other values with probability $1$. If
$ G(u_n) = \sum_{j=1}^{n-1} c_j G(u_j) $
is such a representation, then
\[
\gamma \bigl( \bigl\{ G\in X\colon (G(u_1), \ldots , G(u_n))\in B \bigr\}
\triangle \bigl\{ G\in X\colon (G(u_1), \ldots , G(u_{n-1}))\in B' \bigr\} \bigl) = 0\,,
\]
where
$ B' = \bigl\{ (t_1,  \ldots , t_{n-1})\in\bR^{n-1}\colon
\bigl( t_1, \ldots , t_{n-1}, \sum_{j=1}^{n-1} c_j t_j  \bigr) \in B \bigr\} $
is a Borel set in $\bR^{n-1}$, so we can remove the point $u_n$ from the
consideration at no cost.}.
In this case, we can write
\[
\gamma \bigl( P(u_1, \ldots , u_n; B) \bigr) = (2\pi)^{-n/2}\,
\bigl( \det \Lambda \bigr)^{-\frac12}\,
\int_B e^{-\frac12 (\Lambda^{-1} t \cdot t)}\, \rd\vol (t)\,,
\]
where $\Lambda = \bigl( k(u_i-u_j) \bigr)_{i,j = 1}^n$ is the covariance
matrix of the vector $\bigl( G(u_1), \ldots , G(u_n) \bigr)$. As before, we denote by
$k$ the Fourier integral  of the spectral measure $\rho$.

Since $\tau_v P = P(u_1+v, \ldots , u_n+v; B)$, we have
\[
P \cap \tau_v P
= P(u_1, \ldots , u_n, u_1+v, \ldots , u_n+v; B\times B)\,.
\]
Then
\[
\gamma \bigl( P \cap \tau_v P \bigr)
= (2\pi)^{-n}\, \bigl( \det \widetilde\Lambda \bigr)^{-\frac12}\,
\int_{B\times B} e^{-\frac12 (\widetilde\Lambda^{-1}(v) \widetilde t \cdot \widetilde t)}\,
\rd\vol (\,\widetilde t\,)\,
\]
where
\[
\widetilde\Lambda(v) =
\begin{pmatrix} \Lambda & \Theta (v) \\ \Theta^*(v) & \Lambda \end{pmatrix} \qquad \text{with}\  \Theta_{i, j}(v) = k(u_i - v - u_j)\,.
\]
Note that the matrix $\widetilde\Lambda (v)$ is
invertible and $\bigl( \widetilde\Lambda (v) \bigr)^{-1}$ is close to
$\begin{pmatrix} \Lambda^{-1} & 0 \\ 0 & \Lambda^{-1} \end{pmatrix}$ if $\| \Theta (v)  \|$
is small enough.

Next, we observe that we can
choose a sequence $v_\ell\in\bR^m$ so that $\| \Theta (v_\ell) \|\to 0$ as $\ell\to\infty$.
Indeed, letting  $ \Delta = \max_{i, j} |u_i-u_j| $, we have
\[
\frac1{\vol B(R)}\, \int_{B(R)} \sum_{i, j} k(u_i-v-u_j)^2\, \rd\vol (v)
\le \frac{n^2}{\vol B(R)}\, \int_{B(R+\Delta)} k^2\, \rd\vol\,,
\]
while by Wiener's theorem~\cite[VI.2.9]{Katz}, the absence of atoms in $\rho$ is equivalent to
\[
\lim_{R\to\infty} \frac1{\vol B(R)}\, \int_{B(R)} k^2\,\rd\vol = 0\,.
\]
Then, using the dominated convergence theorem, we conclude that
\[
\lim_{\ell\to\infty} \gamma \bigl( P \cap \tau_{v_\ell} P \bigr)
= \gamma \bigl( P \bigr)^2\,.
\]
Recalling that
$A \cap \tau_{v_\ell} A = A$ up to $\gamma$-measure $0$, we obtain
\[
\gamma (A) = \gamma (A \cap \tau_{v_\ell} A) \le \gamma (P \cap \tau_{v_\ell} P)
+ 2\e\, \stackrel{\ell\to\infty}\to \, \gamma (P)^2 + 2\e \le \bigl[ \gamma (A) \bigr]^2 + 2\e\,.
\]
Since $\e>0$ is arbitrary, we conclude that $\gamma (A) \le \gamma (A)^2$, whence
$\gamma (A) = 0$ or $\gamma (A) = 1$. \mbox{} \hfill $\Box$

\section{Condition $(\rho4)$}\label{App-C}

Here, we collect several observations that, in many instances, help
to verify condition $(\rho4)$. Recall that this condition asserts that
\begin{itemize}
\item
{\em there exist a finite compactly supported Hermitian
measure $\mu$ with $\spt(\mu) \subset \spt(\rho) $ and a bounded domain $D\subset\bR^m$ such that
$\cF\mu\big|_{\partial D}<0$ and $ (\cF\mu) (u_0) > 0 $ for some $u_0\in D$}.
\end{itemize}
Throughout this section, we assume that condition
$(\rho3)$ is satisfied, that is, that the measure $\rho$ is not supported on
a hyperplane in $\bR^m$.

\subsection{Quadratic hypersurface criterion}
{\em The support of any measure $\rho$ \underline{not} satisfying condition
$(\rho4)$ must be contained in a quadratic variety $A\la\cdot\la=b$, where $A$ is an $m\times m$ symmetric matrix
and $b\in\bR^m$}.

\medskip\noindent{\em Proof\,}: Suppose that $\spt(\rho)$  is not contained in any
quadratic variety of the above form. Then $\frac12 m(m+1) +1$-dimensional vectors
\[
v(\la) = \bigl\{ 1, \ \la(i)\la(j)\colon 1\le i \le j \le m \bigr\},
\quad \la\in\spt(\rho),
\]
span $\bR^{\frac12 m(m+1) +1}$ (here $\la (i)$ denotes the $i$-th coordinate of $\la$).
Then we can create two finite linear combinations of cosines:
\[
f(x) = \sum_{\la\in\spt(\rho)} a_\la \cos\bigl( 2\pi \la\cdot x \bigr), \quad
g(x) = \sum_{\la\in\spt(\rho)} b_\la \cos\bigl( 2\pi \la\cdot x \bigr),
\]
such that
\[
f(0) = 1, \quad (D^2 f)(0) = 0\,,
\]
and
\[
g(0) = 0, \quad (D^2 g)(0) = I\,,
\]
where $ D^2 f $ is the matrix with the entries $\partial^2_{x_i\, x_j} f$
and $I$ is the unit matrix. Note that we also automatically have
$Df(0)=Dg(0)=0$, Then the function $h=\e^2 f - g$ will satisfy
$h(0)=\e^2$ and $h(x)<0$ on $\bigl\{ |x|=2\e \bigr\}$, provided that $\e$ is
small enough. \hfill $\Box$

\subsection{Pjetro Majer's interior point criterion}\label{subsect:A-PietroMajer}
The next observation is due to Pietro Majer.

\medskip\noindent
{\em Let the interior of the convex hull of $\spt(\rho)$ contain
a point from $\spt (\rho)$. Then condition $(\rho 4)$ is satisfied}.

\smallskip\noindent In particular, condition $(\rho4)$ is satisfied when $0\in\spt(\rho)$.

\medskip\noindent{\em Proof\,}: Let $\upsilon$ be such a point. Since $\upsilon$ lies in the interior of
the convex hull of $\spt(\rho)$, there are $\la_1, \ldots , \la_n\in \spt(\rho)$ that span the whole
space $\bR^m$, such that
\[
\upsilon = \sum\nolimits_i t_i \la_i\,, \qquad
t_i \ge 0,\ \sum\nolimits_i t_i = \alpha < 1\,.
\]

Consider the function
\[
f(x) = \sum\nolimits_i b_i \cos\bigl( 2\pi \la_i \cdot x \bigr) - \cos (2\pi \upsilon \cdot x)
\]
with $b_i = \alpha t_i + n^{-1}(1-\alpha^2 + \e)$, where $\e>0$.  Then, for $x\to 0$,
\[
f(x) = \Bigr[ \sum\nolimits_i b_i -1 \Bigl] - 2\pi^2 \Bigr[ \sum\nolimits_i b_i (\la_i\cdot x)^2 - (\upsilon\cdot x)^2 \Bigl] + o(|x|^2)\,.
\]
In particular,
\[
f(0) = \sum\nolimits_i b_i -1 = \e >0\,.
\]
Next, we note that
\begin{multline*}
(\upsilon\cdot x)^2 = \Bigl( \sum\nolimits_i t_i \la_i \cdot x \Bigr)^2 = \Bigl( \sum\nolimits_i t_i^{1/2}\, t_i^{1/2} \la_i \cdot x \Bigr)^2 \\[7pt]
\le \Bigl( \sum\nolimits_i t_i \Bigr) \cdot \Bigl( \sum\nolimits_i t_i (\la_i\cdot x)^2 \Bigr) = \alpha\, \sum\nolimits_i t_i (\la_i\cdot x)^2\,.
\end{multline*}

Now, suppose that $x$ belongs to the non-degenerate ellipsoid
\[
E = \Bigl\{\, \sum\nolimits_i (\la_i \cdot x)^2 = \frac{\e n}{\pi^2(1-\alpha^2)}\, \Bigr\}\,.
\]
Since $\la_1, \ldots , \la_n$ span $\bR^n$, we have
\[
|x|^2 = O(\e)\,, \qquad \e\to 0, \, x\in E\,.
\]
Therefore, for $x\in E$ and $\e\to 0$, we have
\begin{multline*}
f(x) \le \e - 2\pi^2 \sum\nolimits_i (b_i - \alpha t_i) (\la_i \cdot x)^2 + o(\e) \\[7pt]
= \e - \frac{2\pi^2 (1-\alpha^2+\e)}n\, \sum\nolimits_i (\la_i \cdot x)^2 + o(\e) < \e - 2\e + o(\e) < 0\,,
\end{multline*}
completing the proof. \hfill $\Box$

\subsection{Analytic closure criterion}
Our last observation is that
\begin{itemize}
\item
the requirement
$\spt(\mu)\subset\spt(\rho)$ in condition ($\rho$4) can be relaxed to the requirement
$\spt(\mu)\subset\spt_{\,\tt r.a.}(\rho)$ where $\spt_{\,\tt r.a.}(\rho)$ is the intersection of all
real-analytic varieties containing $\spt(\rho)$.
\end{itemize}
Note that every quadratic variety is an analytic variety as well, so if $\spt(\rho)\subset V$ then $\spt_{\,\tt r.a.}(\rho)
\subset V$ too. Sometimes, $ \spt_{\,\tt r.a.}(\rho) $  is much larger that $\spt (\rho)$ and satisfy the assumption of~\ref{subsect:A-PietroMajer}
(or some other condition sufficient for establishing ($\rho$4) without $\spt(\rho)$ doing so).
For instance, suppose that $m=2$, $S\subset\bR^2$ is the unit circumference, and $\spt(\rho)\subset S$ is an infinite
set. Since infinite subsets of $S$ are uniqueness sets for real-analytic functions on $S$, we see that $\spt_{\, \tt r.a.}(\rho)=S$.
Then, taking $\mu=m_1$ (the Lebesgue measure on $S$), we conclude that condition $(\rho4)$ is satisfied.

\medskip\noindent{\em Proof\,}:
Let $Q\subset \bR^m$ be a compact set. Consider two linear subspaces of the space $C(Q)$ of {\em real-valued} continuous
functions on $Q$:
\[
X = \bigl\{ \mathcal F\mu\colon
\mu \text{\ is Hermitian, compactly supported}, \spt(\mu)\subset\spt(\rho) \bigr\}
\]
and
\[
X_{\, \tt r.a.} = \bigl\{ \mathcal F\mu\colon
\mu \text{\ is Hermitian, compactly supported}, \spt(\mu)\subset\spt_{\,\tt r.a.}(\rho) \bigr\}.
\]
We need to show that the $C(Q)$-closure of $X$ contains $X_{\, \tt r.a.}$. We will be using a simple duality argument.
Suppose that a signed measure $\nu$ supported by $Q$ annihilates $X$, that is,
\[
\int_Q \bigl( \mathcal F \mu \bigr)\, \rd\nu = 0  \qquad \text{for all admissible\ } \mu\,.
\]
Taking $\mu = \frac12 (\delta_\la + \delta_{-\la})$, $\la\in\spt(\rho)$, we find that the cosine-transform
\[
\bigl( \mathcal C\nu \bigr)(\la) = \int_Q \cos\bigl( 2\pi \la\cdot x \bigr)\,\rd\nu (x)
\]
vanishes on $\spt(\rho)$. However, $\mathcal C\nu$ is an entire function, which is real on $\bR^m$.
Hence, if it vanishes on $\spt(\rho)$, it must also vanish on $\spt_{\, \tt r.a.}(\rho)$. Therefore,
the measure $\nu$ annihilates the subspace $X_{\, \tt r.a.}$ as well. \mbox{} \hfill $\Box$

\end{appendices}

\end{document}